%% file: apcfe.tex
\documentclass[sn-mathphys-num,Numbered,pdflatex]{sn-jnl}
\usepackage{lmodern}
\usepackage{graphicx} 
\usepackage{amsmath}
\usepackage{amssymb}
\usepackage{amsthm}
\usepackage{mathtools}
\usepackage[sort]{cleveref}
\usepackage{subcaption}
\usepackage{url}
\usepackage{xspace}
\usepackage{textcomp}
\usepackage{manyfoot}
\usepackage{booktabs}
\usepackage{epsfig}
\usepackage{layouts}
\usepackage{nicefrac}
\usepackage{siunitx}
\usepackage{enumerate}
\usepackage{comment}
\newcommand{\eps}{\varepsilon}
\newcommand{\Ass}{\ensuremath{A^{7\times 7}}}
\newcommand{\R}{\mathbb R}
\newcommand{\C}{\mathbb C}
\renewcommand{\i}{\mathrm{i}}
\newcommand{\Eigen}{\texttt{Eigen3}\xspace}
\newcommand{\Id}{\mathrm{Id}}
\renewcommand{\Re}{\mathrm{Re}}
\newcommand{\im}{\mathrm{imag}}
\newcommand{\re}{\mathrm{real}}
\newcommand{\Nmax}{{N_{\mathrm{max}}}}
\newcommand{\Pmax}{{P_{\mathrm{max}}}}
\date{\today}
\begin{document}
\title[Arbitrary precision computation]{Arbitrary precision computation of hydrodynamic stability eigenvalues}
\author*[1]{\fnm{Patrick} \sur{Dondl}}\email{patrick.dondl@mathematik.uni-freiburg.de}
\author[2]{\fnm{Ludwig} \sur{Striet}}\email{ludwig.striet@mathematik.uni-freiburg.de}
\author[3]{\fnm{Brian} \sur{Straughan}}\email{brian.straughan@durham.ac.uk}
\affil*[1]{\orgdiv{Abteilung f\"ur Mathematik}, \orgname{Albert - Ludwigs - Universit\"at Freiburg}, 
\orgaddress{\street{Hermann Herder Strasse},
\city{Freiburg}, \postcode{79104}, \country{Germany}}}
\affil[2]{\orgdiv{Abteilung f\"ur Mathematik}, \orgname{Albert - Ludwigs - Universit\"at Freiburg}, 
\orgaddress{\street{Hermann Herder Strasse},
\city{Freiburg}, \postcode{79104}, \country{Germany}}}
\affil[3]{\orgdiv{Department of Mathematics}, \orgname{University of Durham}, \orgaddress{\street{Stockton Road},
\city{Durham}, \postcode{DH1 3LE}, \country{U.K.}}}
\abstract{
We show that by using higher order precision arithmetic, i.e., using floating point types with more significant bits than standard \texttt{double} precision numbers, one may accurately compute eigenvalues for non-normal matrices arising in hydrodynamic stability problems. The basic principle is illustrated by a classical example of two $7\times 7$ matrices for which it is well known that eigenvalue computations fail when using standard double precision arithmetic. We then present an implementation of the Chebyshev tau--QZ method  allowing the use of a large number of Chebyshev polynomials together with arbitrary precision arithmetic. This is used to compute the behavior of the spectra for Couette and Poiseuille flow at high Reynolds number. An experimental convergence analysis finally makes it evident that high order precision is required to obtain accurate results.}
\maketitle
\section{Introduction}
\setcounter{equation}{0}
The problem of instability of parallel shear flows in fluid mechanics is one which has increasingly occupied current research. In classical linear theory of hydrodynamic instability the leading eigenvalue, i.e., the one with largest real part, is traditionally examined to deduce whether a flow is stable or not. The shortcomings of this approach are well known in some areas of fluid mechanics where, for example, in Couette flow linear instability theory never predicts instability. A similar outcome arises for Poiseuille flow, albeit with a finite Reynolds number. However, transition to instability does occur in real life, and much below the critical Reynolds numbers predicted by linear theory. The review by
\citet{Kerswell:2018} explains this concept in an extremely lucid manner and also describes how other approaches have modified the 
classical technique. For example, the idea of transient growth, cf.\ \citet{ButlerFarrell:1992}, shows that by considering
eigenvalues other than simply the one with largest real part, one may find modes have very rapid growth over a relatively short interval of time. One difficulty with such an approach is that numerical calculation of accurate values of some eigenvalues
in the spectrum away from the leading one is not reliable, at least in standard \texttt{double} precision arithmetic. Thus, to have any faith in such an approach, one 
requires a numerical method which will guarantee a prescribed accuracy for all eigenvalues, including those for which, due to non self-adjointness, the eigenvectors are nearly linearly dependent. The purpose of this article is to present a method whereby one may calculate an
eigenvalue to the required accuracy by employing arbitrary precision arithmetic. 

We stress that for parallel shear flow the theory has been substantially developed in alternative directions, in particular, by 
nonlinearity to nonmodal analysis. We again refer to the review by \citet{Kerswell:2018}. Other developments concern growth to a 
localized nonlinear optimal, see \citet{Pringle:2015}, or a three-dimensional statistical stability analysis, see 
\citet{MarkeviciuteKerswell:2023}.  For such shear flow problems such as Couette or Poiseuille flow we here present an
arbitrary precision development of the Chebyshev tau--QZ algorithm, cf.\
\citet{Orszag:1971}, \citet{MolerStewart:1971}, \citet{Straughan1996}, which
additionally allows one to use a large number of Chebyshev polynomials.

While we concentrate on describing accurate calculation of the spectra for Couette and Poiseuille flows in a linearly viscous fluid
we point out that there has been a tremendous amount of recent activity in the calculation of such spectra involving non-normal matrices
in other problems of real life. For example, parallel shear flow with a viscoelastic liquid, \citet{ShankarShivakumara:2023};
magnetohydrodynamic channel flow of a Navier-Stokes-Voigt fluid, \citet{Kavitha:2024}; Poiseuille flow in a non-ideal fluid, 
\citet{Zheng:2024};
shear flows in porous media, \citet{Arjun:2025}, \citet{ShivarajKumar:2024}, \citet{Samanta:2022};
shear and related flows for a fluid overlying various types of porous material, an area with many real life applications, see
\citet{DondlStraughan:2025}, including blood flow in an artery, \citet{Samanta:2020,Samanta:2023}, \citet{Sharma:2025} or
Faraday instabilities, \citet{Samanta:2020prsa}. 

A complete stability analysis  in each of these problems will require accurate calculation of a large part of the spectrum and as the
matrices are non-normal, a truly accurate calculation of the spectrum will necessarily require the ideas of this paper where many
spectral polynomials are involved, but crucially arbitrary precision is also required. The above lists examples in fluid mechanics,
but recent work has also identified the need for the calculation of spectra of non-normal matrices in other areas of applied mathematics.
Notably in transmitted microwave radiation, \citet{DavyGenack:2019}; and also in solid mechanics, particularly for wave motion in
quasicrystals, \citet{Benzi:2025}. Thus, we believe the ideas of this paper and subsequent development from our codes will be invaluable in
future work on these areas.

In the next section we commence with two $7\times 7$ matrices of \citet{Godunov:1992} and of \citet{TrefethenEmbree:2005}, which have eigenvalues
of 0,1,1,$\pm$2,$\pm$4,
and
0,$\pm$1,$\pm$2,$\pm$4,
respectively.
In connection with the \citet{Godunov:1992} example \citet{Brown:2010} write, ...
``{\sl no matter what software is used, numerical computations yield a set of complex eigenvalues, ... with imaginary parts as large as 
8.33, which are nowhere near the true eigenvalues,"} cf.\ also \citet{HinrichsenPritchard:2005}. They infer that this is due to the 
non-normality of the $7\times 7$ matrix so that the eigenvalues are subject to unavoidable round off errors and consequently the eigenvalue
computations cannot be trusted. We show in the next section that we may resolve this numerical accuracy problem by employing higher
precision arithmetic. We give examples of the computation of eigenvalues for plane Poiseuille flow and Couette flow at high Reynolds
numbers where double precision calculations are totally inaccurate and we show how this may be overcome by a suitable combination of
Chebyshev polynomials and high precision arithmetic. The convergence of the solution is checked by computing the Hausdorff distance between
sets of eigenvalues. We then show that in many physical regimes the eigenvalues in the zone where they are close to being linearly 
dependent, depend critically on the number of polynomials \emph{and} on the higher precision arithmetic. 

Our multi precision implementation of the QZ algorithm is an adaption of
the Real QZ algorithm provided by Alexey Korepanov with \Eigen that works with
complex input data.
Furthermore, we use the shift strategy explained in \cite{MolerStewart:1971}. We remark
that many improvements to the algorithm presented there have been presented in literature
in recent years that could help to further speed up convergence of the method.

By consistently templating for the scalar types in our \texttt{C++} implementation
of the algorithm, it is easily possible to replace built-in \texttt{double} variables
by multi precision floating point numbers. We use the wrapper \cite{mpfrcpp} to
the \cite{mpfrpaper} multi precision floating point library in our implementations.

\section{A Classical Example}
We consider the $7\times7$ matrix
\begin{equation*}
A_s^{7\times7} \coloneqq L^{-1} \tilde A_s^{7\times 7} L
\end{equation*}
with 
\begin{equation*}
\resizebox{\textwidth}{!}{$
\tilde A_s^{7x7} = \begin{pmatrix}
1 & 2048 &  256 & 128 &   64 &   32 &   16 \\ 
0 &   -2 & 1024 & 512 &  256 &  128 &   32 \\
0 &    0 &    4 & 512 & 1024 &  256 &   64 \\
0 &    0 &    0 &   0 &  512 &  512 &  128 \\
0 &    0 &    0 &   0 &   -4 & 1024 &  256 \\
0 &    0 &    0 &   0 &    0 &    2 & 2048 \\
0 &    0 &    0 &   0 &    0 &    0 &    s
\end{pmatrix}\text{, }
L = \begin{pmatrix}
1\; & \;0\; & \;0\; & \;0\; & \;0\; & \;0\; & \;0 \\ 
  		0\; & \;1\; & \;0\; & \;0\; & \;0\; & \;0\; & \;0 \\ 
1\; & \;0\; & \;1\; & \;0\; & \;0\; & \;0\; & \;0 \\ 
0\; & \;0\; & \;0\; & \;1\; & \;0\; & \;0\; & \;0 \\ 
0\; & \;0\; & \;1\; & \;0\; & \;1\; & \;0\; & \;0 \\ 
1\; & \;0\; & \;0\; & \;0\; & \;0\; & \;1\; & \;0 \\ 
0\; & \;1\; & \;1\; & \;0\; & \;1\; & \;0\; & \;1
\end{pmatrix}$}
\end{equation*}
which, going back to \cite{Godunov:1992} and pointed out in \cite[p. 547]{Hinrichsen2005},
has the spectrum
\begin{equation}
\label{eq:spectrum_Ass_p1}
\sigma(A_1^{7\times 7}) = \{ -4, -2, 0, 1, 1, 2, 4 \}
\end{equation}
in the case $s = 1$ and as pointed out in \cite[p. 489]{TrefethenEmbree:2005}, has the spectrum
\begin{equation}
\label{eq:spectrum_Ass_m1}
\sigma(A_{-1}^{7\times 7}) = \{ -4, -2, -1, 0, 1, 2, 4 \}
\end{equation}
in the case $s = -1$. It is reported in the aforementioned works that computing
the spectra of the matrix $\Ass_{s}$, $s\in\{-1,1\}$, with any standard numeric software
in usual double precision floating arithmetic fails. Indeed, computing the spectra
using standard \texttt{double} precision in the \texttt{EigenSolver} from the
\texttt{Eigen::Eigenvalues} module of the \Eigen library \cite{Eigen3}, which
computes the eigenvalues using a Schur decomposition gives, rounded to three decimal places, the
spectra
\begin{equation*}
\sigma_d(\Ass_{-1}) = \{ -5.446\pm2.095, -1.355\pm 4.859\i, 3.738\pm 3.904\i, 6.127 \}
\end{equation*}
and 
\begin{equation*}
\sigma_d(\Ass_{1}) = \{ -5.342\pm 2.149\i, -1.056\pm 5.053\i,4.139\pm 4.108\i, 6.517\}.
\end{equation*}
Computing the spectra using the \texttt{eig}-function of the free software \texttt{GNU Octave}
gives
\begin{equation*}
\sigma_d(\Ass_{-1}) = \{ 6.630 \pm 2.812\i, 1.590 \pm 6.328\i, -4.562 \pm 5.010\i, -7.316  \}
\end{equation*}
and 
\begin{equation*}
\sigma_d(\Ass_{1}) = \{ 6.788 \pm 2.787\i, 1.868 \pm 6.223\i, -4.187 \pm 4.865\i, -6.938\}.
\end{equation*}
All of the above obviously do not at all resemble the correct spectra
$\sigma\left(\Ass_{-1}\right)$ and $\sigma\left(\Ass_{1}\right)$.
However, the situation changes if we compute the results using multi precision
floating point arithmetic. To show this, we compute the eigenvalues of the aforementioned
$7\times 7$ matrix with higher precision using the \texttt{mpfr} library \cite{mpfrpaper} 
with the \texttt{C++} interface \texttt{MPFR C++} \cite{mpfrcpp} for which \Eigen
has built-in support. Here and throughout this article, we refer with $P$ to
the number of significant bits and with
\begin{equation*}
\eps_P = 2\cdot 2^{-P}
\end{equation*}
to the \emph{interval machine precision}. We refer with
$\sigma_P\left(\Ass_s\right)$ to the spectrum of $\Ass_s$ that we compute
in floating point arithmetic with $P$ significant bits.

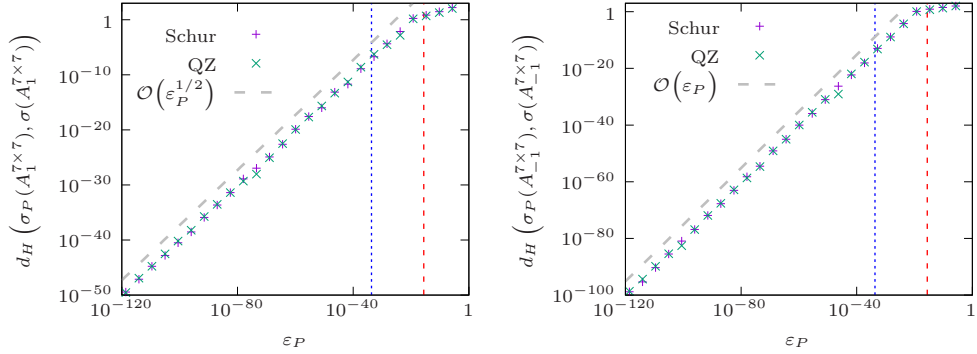
\begin{figure}
\begin{center}
\resizebox{\textwidth}{!}{
\input{apcfe-paper-figure-convana_7x7.txt}
}
\end{center}
\caption{The Hausdorff distances of the spectra $\sigma_P(\Ass_{-1})$ (left) and $\sigma_P(\Ass_{1})$ (right)
to the correct spectra  $\sigma(\Ass_{-1})$ (left) and $\sigma(\Ass_{1})$ (right) for different numbers of significant
bits $P$ or respective machine precisions $\eps_P$. The distance decreases like $\mathcal O\big(\eps_P^{\nicefrac{1}{2}}\big)$ in 
the case $s = 1$ and like $\mathcal O\left(\eps_P\right)$ in the case $s = -1$.
The dashed red and blue lines indicate the machine precision $\eps_{\texttt{double}} = 2.2204\cdot 10^{-16}$
for standard double precision with 53 significant bits and $\eps_{\texttt{ext}} = 1.9259\cdot 10^{-34}$ for
extended precision with 113 significant bits, respectively.}
\label{fig:convana_7x7}
\end{figure}

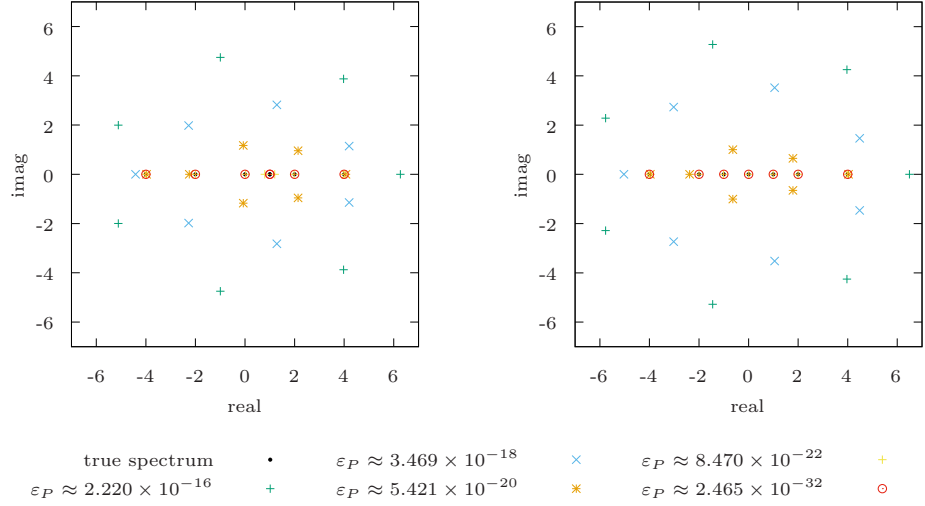
\begin{figure}
\begin{center}
\resizebox{\textwidth}{!}{
\input{apcfe-paper-figure-7x7_spects.txt}
}
\end{center}
\caption{The spectra $\sigma_P(\Ass_{1})$ (left) and $\sigma_P(\Ass_{-1})$ (right) that we compute 
for different numbers of significant bits $P$. For high numbers of $P$, we clearly obtain
the correct spectra given in \cref{eq:spectrum_Ass_m1,eq:spectrum_Ass_p1}.}
\label{fig:spectrum_7x7}
\end{figure}

We compute the spectra of $\Ass_s$ for $s\in\{-1,1\}$ using (i) the Eigenvalue solver
provided by \Eigen version 3.4 and (ii) our implementation
of the QZ-algorithm where we set $B = \Id\in\R^{7\times 7}$. Afterwards, we compute
the Hausdorff distances\footnote{The \emph{Hausdorff distance} between two non-empty sets $\mathcal A,\mathcal B\subset\C$
is defined as $
d_H(\mathcal A, \mathcal B) = \max\left\{
\max\{D(a, \mathcal B)\,:\,a\in \mathcal A\},\max\{D(b, \mathcal A)\,:\,b\in \mathcal B\}
\right\}$ where $D(a, \mathcal B) = \min\{|a-b|\,:\,b\in \mathcal B\}$
and $D(b, \mathcal A)$ is defined analogously.} $d_H(\sigma(\Ass_s),\sigma_P(\Ass_s))$ of the correct
spectrum and the spectrum we compute numerically using $P$ significant
bits. 
We show the spectra that we compute using different precisions in \cref{fig:spectrum_7x7}
and the distances in \cref{fig:convana_7x7}. We observe that
using floating point numbers with more significant bits decreases the error
and gives valuable results for the computation of eigenvalues. The distance
decreases like $\mathcal O\big(\eps_P^{\nicefrac{1}{2}}\big)$ in the case $s = 1$ and like
$\mathcal O\left(\eps_P\right)$ in the case $s = -1$. More generally, we observe
that the rate is $-1/2$ if one ore more eigenvalues occur twice and $-1$ if all
the eigenvalues occur only once.

The promising results for this artificial example motivate the use of high-precision
floating point arithmetic to compute spectra of hydrodynamic stability problems.

\section{Chebyshev-tau method and the Orr-Sommerfeld equation}
We briefly summarize the Chebyshev-tau method as we use it to compute the eigenvalues
of the hydrodynamic problems as described in the introduction here. For details,
we refer to \cite{Straughan1996}.

We begin the  discussion with the Orr-Sommerfeld equation which reads
\begin{equation}
\label{eq:orr-sommerfeld}
(D^2-a^2)^2 \phi = \i a \Re (U-c)(D^2 - a^2) \phi - \i a \Re U'' \phi,\quad z\in(-1,1).
\end{equation}
Here, $D$ is the derivative w.r.t the variable $z$ and $\Re$, $a$ and $c$ are the
Reynolds number, wavenumber and the (sought) eigenvalue, respectively.
We aim to solve \cref{eq:orr-sommerfeld}
with $U = (1-z^2)$ in the case of Poiseuille flow and $U = z$ in the case of Couette
flow, both subject to the boundary conditions
\begin{equation}
\label{eq:orr-sommerfeld-bc}
\phi = D\phi = 0\text{ for } z\pm 1.
\end{equation}
We aim to solve \cref{eq:orr-sommerfeld} using the $D^2$ and the $D^4$ method.
The $D^2$ method is reported to be a good compromise between matrix size and growth
of its entries \cite{Straughan1996}. However, the growth of the entries is less
of a problem here, as we use high-precision numerical computations and the actual
computing (wall-)time is influenced more by the size of the matrix than by the number of significant bits. This makes for a compelling argument to use the $D^4$-method, where the matrix size is halved compared to the $D^2$ method.
Still, we mainly focus the discussion in this work to the $D^2$ method for
the sake of comparability to other works in this field.

\subsection{The \texorpdfstring{$D^2$}{D2}-method}
\label{sub:D2_method}
In the $D^2$-method, we write \cref{eq:orr-sommerfeld} as two equations containing
the operators $L_1$, $L_2$ defined via
\begin{equation}
\label{eq:orr-sommerfeld-D2}
\begin{aligned}
L_1(\phi,\chi) &\equiv (D^2-a^2)\phi - \chi = 0  \\
L_2(\phi,\chi) &\equiv (D^2-a^2)\chi - \i a \Re (U - c)\chi + \i a \Re U''\phi = 0.
\end{aligned}
\end{equation}
We expand $\phi$ and $\chi$ in Chebyshev polynomials via
\begin{equation*}
\phi = \sum\limits_{i=0}^{N+2}\phi_i T_i(z),\quad \chi = \sum\limits_{i=0}^{N+2} \chi_i T_i(z)
\end{equation*}
and solve the equations
\begin{equation}
\label{eq:boundary_conditions}
L_1(\phi,\chi) = \tau_1 T_{N+1} + \tau_2 T_{N+2},\quad L_2(\phi,\chi) = \tau_3 T_{N+1} + \tau_4 T_{N+2}
\end{equation}
Multiplying by $T_i, i = 0,\ldots,N$ yields $2(N+1)$ equations for the coefficients
$\phi_i,\chi_i$ \cite{Straughan1996}. We write the boundary conditions in \cref{eq:orr-sommerfeld-bc}
as
\begin{equation*}
\sum\limits_{\substack{i=0\\i\text{ even}}}^{N+2} \phi_i = 0,\quad
\sum\limits_{\substack{i=1\\i\text{ odd}}}^{N+2} \phi_i = 0,\quad
\sum\limits_{\substack{i=1\\i\text{ odd}}}^{N+2} i^2\phi_i = 0,\quad
\sum\limits_{\substack{i=0\\i\text{ even}}}^{N+2} i^2\phi_i = 0.
\end{equation*}
Computing the eigenvalues $c$ from \cref{eq:orr-sommerfeld-D2} then reduces  to solving
the generalized eigenvalue problem $A x = c B x$ with $x = (\phi,\chi)$ and
$A^{N} \in \C^{2(N+3)\times 2(N+3)}, B^{N} \in \C^{2(N+3)\times 2(N+3)}$ where,
in case of Poiseuille flow,
\begin{equation*}
\label{eq:def_A_pois}
A_{\text{real}}^N = \begin{pmatrix}
D^2 - a^2\Id & -\Id \\
\mathrm{BC1} & 0\cdots 0 \\
\mathrm{BC2} & 0\cdots 0 \\
\mathbf 0 & D^2 - a^2\Id \\
\mathrm{BC3} & 0\cdots 0 \\
\mathrm{BC4} & 0\cdots 0 \\
\end{pmatrix},\quad
A_\text{imag}^N = \begin{pmatrix}
\mathbf 0 & \mathbf 0 \\
0\cdots0 & 0\cdots0 \\
0\cdots0 & 0\cdots0 \\
-2a\,\Re\,\Id & a\,\Re\,(\Pi^{z^2} - \Id)\\
0\cdots0 & 0\cdots0 \\
0\cdots0 & 0\cdots0 \\
\end{pmatrix},
\end{equation*}
and
\begin{equation}
\label{eq:def_B_pois}
B_{\text{real}}^N = \mathbf 0,\quad B_\text{imag}^N = \begin{pmatrix}
\mathbf 0 & \mathbf 0 \\
0\cdots 0 & 0\cdots0\\
0\cdots 0 & 0\cdots0 \\
\mathbf 0 & -a\,\Re\,\Id \\
0\cdots 0 & 0\cdots0\\
0\cdots 0 & 0\cdots0.
\end{pmatrix}
\end{equation}
The row vectors $\mathrm{BC1}-\mathrm{BC4}$ express the boundary conditions \cref{eq:boundary_conditions}. The matrix $\Id\in \R^{(N+1)\times (N+3)}$ is the matrix whose left 
block is the $(N+1)\times(N+1)$ identity matrix and that is appropriately filled with $0$ on the
right side.
The matrix $\Pi^{z^2}\in\R^{(N+1)\times (N+3)}$ is the matrix that represents multiplication by $z^2$ and is obtained by writing 
$z^2 = \frac{1 + T_2(z)}{2}$ and taking the inner products $(T_i,z^2\phi)$. Its non-zero
entries are, using $1$-based indices for matrix-entries,
\begin{equation*}
\Pi^{z^2}_{1,1} = \nicefrac{1}{2},\;
\Pi^{z^2}_{1,3} = \nicefrac{1}{4},\;\Pi^{z^2}_{2,2} = \nicefrac{3}{4}, \;\Pi^{z^2}_{2,4} = \nicefrac{1}{4},\;
\Pi^{z^2}_{3,1} = \Pi^{z^2}_{3,3} = \nicefrac{1}{2}, \;\Pi^{z^2}_{3,5} = \nicefrac{1}{4}
\end{equation*}
and
\begin{equation*}
\Pi^{z^2}_{i,i} = \nicefrac{1}{2},\quad \Pi^{z^2}_{i,i-2} = \Pi^{z^2}_{i,i+2} = \nicefrac{1}{4} \;\text{ for } \;i\geq 4.
\end{equation*} 

For Couette flow, we have the same matrices $A^N, B^N$ except that
\begin{equation*}
\label{eq:def_A_coue}
A_\text{imag}^N = \begin{pmatrix}
\mathbf 0 & \mathbf 0 \\
0\cdots0 & 0\cdots0 \\
0\cdots0 & 0\cdots0 \\
\mathbf 0 & -a\,\Re\,\Pi^z\\
0\cdots0 & 0\cdots0 \\
0\cdots0 & 0\cdots0 \\
\end{pmatrix}
\end{equation*}
and that the matrix $\Pi^z\in\R^{(N+1)\times (N+3)}$ is the
matrix that represents multiplication by $z$. Writing $z = T_1(z)$ it is obtained
by taking the inner products $(T_i, z)$. Its non-zero entries are, again using
$1$-based indices,
\begin{equation*}
\Pi^z_{1,1} = \nicefrac{1}{2}, \Pi^z_{2,1} = 1,
\Pi^z_{2,3} = \nicefrac{1}{2}\text{ and } \Pi^z_{i,i-1} = \Pi^z_{i,i+1} = \nicefrac{1}{2}\text{ for }i > 2.
\end{equation*}
The matrix $D^2\in\R^{(N+1)\times(N+3)}$ represents the second derivative operator
$\nicefrac{\mathrm d^2}{\mathrm dz^2}$ on the space of Chebyshev polynomials. It is defined
via
\begin{equation*}
D^2 \phi = \sum\limits_{i=0}^{N+2} \phi_i^{(2)}T_i(z)
\end{equation*}
with
\begin{equation*}
\phi_i^{(2)} = \frac{1}{b_i}\sum\limits_{\substack{p=i+2\\p+i\text{even}}}^{N+2} p(p^2-i^2) \phi_p,
\end{equation*}
and $b_i = 2$ if $i = 0$, $b_i = 1$ otherwise. To avoid any confusion, we display
the top-left blocks of the respective matrices in the following. 

\subsection{The \texorpdfstring{$D^4$}{D4}-method}
\label{sub:D4_method}
Different to the $D^2$ method discussed above, the $D^4$ method abstains from
splitting \cref{eq:orr-sommerfeld} into the two equations \eqref{eq:orr-sommerfeld-D2}. Instead we directly discretize by expanding
\begin{equation*}
\phi = \sum\limits_{i=0}^{N+4}\phi_i T_i(z).
\end{equation*}
Then, we define the operator $L$ via
\begin{equation*}
L\phi = D^4\phi -2a^2D^2\phi - \i a \Re (U-c)(D^2 - a^2)\phi + (a^4 + \i a\Re U'')\phi
\end{equation*}
and solve 
\begin{equation*}
L\phi \equiv \tau_1 T_{N+1} + \tau_2 T_{N+2} + \tau_3 T_{N+3} + \tau_4 T_{N+4},
\end{equation*}
equipped with boundary conditions
\begin{equation*}
\sum\limits_{\substack{i=0\\i\text{ even}}}^{N+3} \phi_i = 0,\quad
\sum\limits_{\substack{i=1\\i\text{ odd}}}^{N+4} \phi_i = 0,\quad
\sum\limits_{\substack{i=1\\i\text{ odd}}}^{N+4} i^2\phi_i = 0,\quad
\sum\limits_{\substack{i=0\\i\text{ even}}}^{N+3} i^2\phi_i = 0.
\end{equation*}
The matrix $D^4$ represents the fourth derivative operator and is defined through the 
identities
\begin{equation*}
D^4\phi = \sum\limits_{i=0}^{N+4}\phi_i^{(4)} T_i(z)
\end{equation*}
with
\begin{equation*}
\phi_i^{(4)} = \frac{1}{24b_i}\sum\limits_{\substack{p=i+4\\p+i\text{ even}}}^{p=N+4}
p\left(
p^2(p^2-4)^2 - 3p^4i^2 + 3p^2i^4 - i^2(i^2 - 4)^2
\right)\phi_p
\end{equation*}
and $b_i = 2$ if $i = 0$, $b_i = 1$ otherwise. The matrix $D^2$ has to be formed
to a matrix of appropriate size accordingly. We abstain from showing the concrete
values of the matrix $D^4$ here do their strong growth.

Computing the eigenvalues $c$ reduces then to solving the generalized eigenvalue
problem $Ax = cBx$ where
\begin{equation*}
\resizebox{\textwidth}{!}{$
A_{\text{real}} = \begin{pmatrix}
D^4 - 2 a^2D^2 + a^4\Id
\end{pmatrix},
A_{\text{imag}} = 
\begin{pmatrix}
-a\,\Re (D^2 - a^2\Id) + a\,\Re\, \Pi\,(D^2 -a^2\Id) - 2a\,\Re\,\Id \\
\text{BC1} \\
\text{BC2} \\
\text{BC3} \\
\text{BC4} \\
\end{pmatrix},
$}
\end{equation*}
with $\Pi = \Pi^{z^2}$ or $\Pi = \Pi^z$, as above in \cref{sub:D2_method}, and
\begin{equation*}
B_{\text{real}} = \mathbf 0,\quad
B_{\text{imag}} = \begin{pmatrix}
-a\,\Re\,(D^2-a^2\Id) \\
0\cdots 0\\
0\cdots 0\\
0\cdots 0\\
0\cdots 0
\end{pmatrix}.
\end{equation*}
In this case, $\Id$ is the $(N+1)\times(N+5)$ matrix that has the $(N+1)\times(N+1)$
identity matrix in its top-left corner and is zero otherwise.
We refer, again, to \cite{Straughan1996} for more details on the $D^4$ method.
As stated there, a major
drawback of the $D^4$ method is the growth of the entries of the matrices which
is $\mathcal{O}\left(N^{7}\right)$ and makes the numerical treatment challenging due to extremely bad conditioning for large $N$. However, this drawback can be mitigated by the use of high precision floating point numbers, where the benefit of smaller system size -- which is half that of the $D^2$-method -- may outweigh the worse conditioning. 

\section{Numerical experiments}
\subsection{The \texorpdfstring{$D^2$}{D2} method for Plane Poiseuille Flow}
In this section, we analyze the output of the $D^2$ Chebyshev tau method for plane Poiseuille
flow as described in \cref{sub:D2_method}.
We refer with $\sigma_P(A^N, B^N)$ to the resulting spectrum of the matrix pair
that we compute using our implementation of the QZ algorithm using floating point
arithmetic with $P$ significant bits. 

\begin{figure}
\begin{center}
\input{apcfe-paper-figure-Straughan-paper-repl_Straugh_Fig1.txt}
\end{center}
\caption{The spectrum for plane Poiseuille flow with $\Re = 10^4$, $a = 1$,
computed in double precision with $N=200$ polynomials
for the even- and odd eigenfunctions. This reproduces the result in Fig.\ 1 of \cite{Straughan1996}.\label{fig:Straughan_Fig1}}
\end{figure}
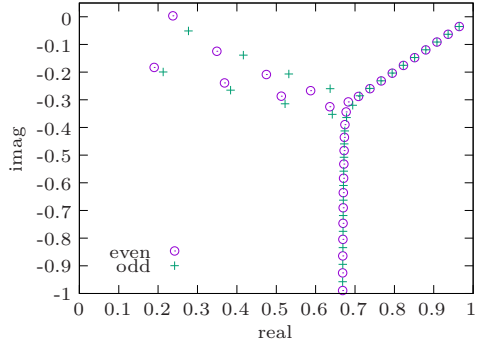

As a first numerical test, we reproduce the 
results obtained by Dongarra, Straughan, and Walker in \cite{Straughan1996}. The obtained eigenvalues  belonging to the odd- and even eigenfunctions of
the Orr-Sommerfeld equation for Poiseuille flow with $Re = 10^4$, $a = 1$ can be found in
\cref{fig:Straughan_Fig1}.
In our further numerical experiments, we analyze the effects of varying $N$ and $P$ on the accuracy with which eigenvalues of the Poiseuille and Couette flow problems can be computed.
The questions we pose are
\begin{enumerate}[(i)]
\item how the number of polynomials $N$ affects the computed spectra and
\item how the number of significant bits $P$ used in the computations affects
the computed spectra.
\end{enumerate}

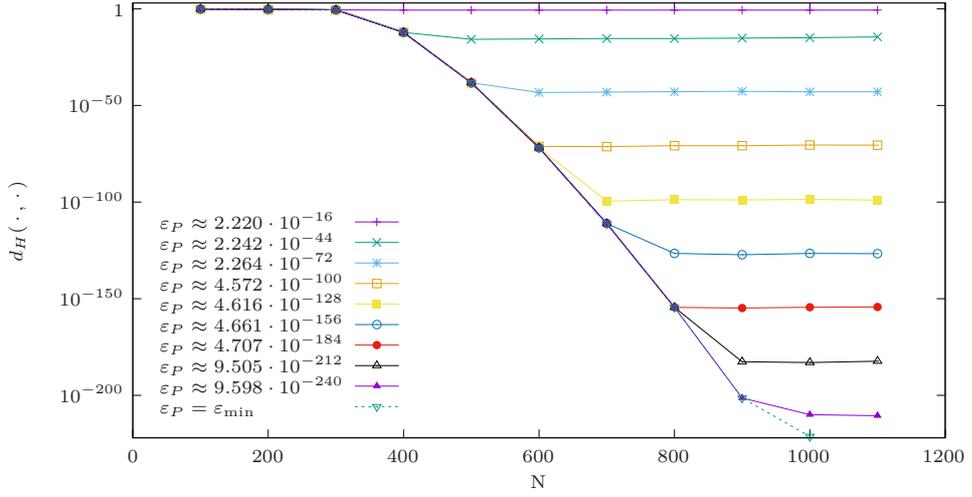
\begin{figure}
\begin{center}
\input{apcfe-paper-figure-dists_NPmax_Re1e5_N_dH.txt}
\end{center}
\caption{\mbox{\textit{Poiseuille flow, $\Re = 10^{5}$}}. The distances 
as described in \cref{eq:def_distances} for Poiseuille flow with $Re = 10^{5}$ and different
values of $P$. We note that that increasing the number of polynomials does
not lead to a decrease of the error if the number of significant bits is not
chosen high enough.}
\label{fig:convana_Re1e5_N_dH}
\end{figure}

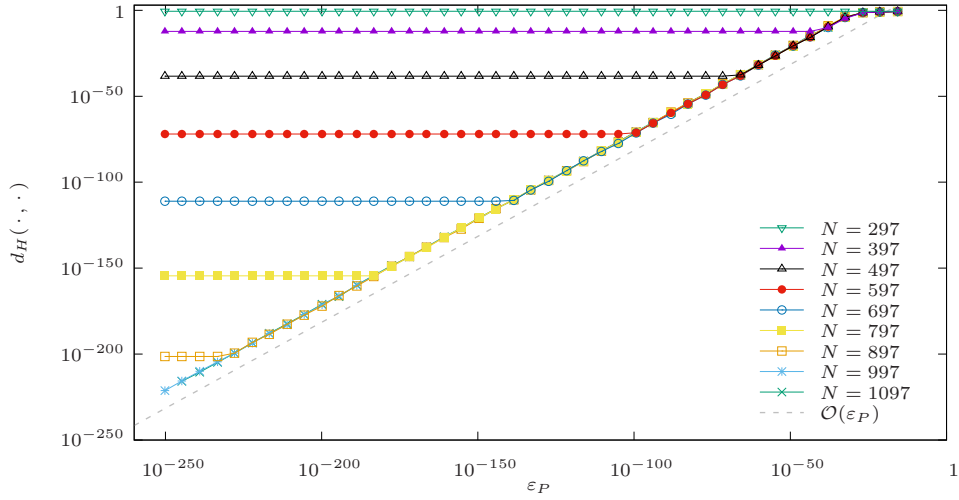
\begin{figure}
\begin{center}
\input{apcfe-paper-figure-dists_NPmax_Re1e5_P_dH.txt}
\end{center}
\caption{\textit{Poiseuille flow, $\Re = 10^5$}. The distances 
as described in \cref{eq:def_distances} for Poiseuille flow with $Re = 10^5$ and different
values of $N$. We note that that increasing the number of significant bits
$P$ (or decreasing $\eps_P$ respectively) does not lead to a decrease of the error
if the number of polynomials is not chosen high enough.}
\label{fig:convana_Re1e5_P_dH}
\end{figure}

To this end, we compute the spectra for fixed Reynolds numbers $Re = 10^5$, $Re = 2\cdot 10^5$, $Re = 5\cdot 10^5$ and different values of $P \leq \Pmax = 832$ and $N\leq \Nmax$. We vary 
the values of $\Nmax$ between $1100$ and $1500$ depending on the Reynolds number.
As usual, we only consider eigenvalues in the domain
\begin{equation*}
Q\coloneqq \{z\in\C\,|\, -1\leq \im(z)\leq 0,\, 0\leq \re(z)\leq 1 \}
\end{equation*}
and compute for
each of the chosen Reynolds numbers and all $N$, $P$, in the sense of an experimental convergence analysis, the Hausdorff distances
\begin{equation}
\label{eq:def_distances}
d_H(\sigma_P(A^N, B^N)\cap Q,\, \sigma_\Pmax(A^{\Nmax}, B^{\Nmax})\cap Q).
\end{equation}
The respective distances for fixed values of $P$ and varying $N$ can be found in
\cref{fig:convana_Re1e5_N_dH,fig:convana_Re2e5_N_dH,fig:convana_Re5e5_N_dH}.
Our results clearly suggest that fixing $\eps_P$ and increasing the
number of polynomials improves the results only up to a certain point. Then, a further
increase of the number of of polynomials only decreases the distance if the $\eps_P$
is decreased (or, equivalently, if $P$ is increased) as well.
Conversely, the results for fixed values of $N$ and varying values of $P$ in 
\cref{fig:convana_Re1e5_P_dH,fig:convana_Re2e5_P_dH,fig:convana_Re5e5_P_dH}.
show that also decreasing $\eps_P$ decreases the Hausdorff distance to the reference spectrum only if the number of polynomials is increased as well.
The aforementioned figures show that a minimal number of polynomials $N$ and precision $P$
that depend on the Reynolds Number $\Re$ is needed to obtain a distance of the spectrum
computed with $N$, $P$ to the spectrum computed with $\Nmax,\Pmax$ that is below
double precision. We show the respective value for $N$ and $P$ to obtain a given accuracy for the spectrum in \cref{tab:need_NP_pois}. Note that, for the Reynolds numbers considered here,  even a modest accuracy of 10\% can \emph{not} be achieved using standard \texttt{double} or, for Reynolds Numbers $5.0\cdot 10^5$, \texttt{extended} precision arithmetic.

\begin{table}
\captionsetup{width=\textwidth}
\caption{The required minimal values for the number of polynomials $N$ and significant bits $P$ to obtain an accuracy of the computed spectra below \texttt{single}, \texttt{double}
and \texttt{extended} precision. As reference we use the spectra computed with the maximal number of polynomials $N_\text{max}$ and maximal number of significant bits $P_\text{max}$.	\label{tab:need_NP_pois}}
\begin{tabular}{r|cc|cc|cc|cc} \multicolumn{9}{c}{}\\
    accuracy & \multicolumn{2}{c|}{\texttt{10\%}} & \multicolumn{2}{c|}{\texttt{single}} & \multicolumn{2}{c|}{\texttt{double}} & \multicolumn{2}{c}{\texttt{extended}} \\
					& $N$ & $P$ & $N$ & $P$ & $N$ & $P$ \\
$\Re = 1.0\cdot 10^5$ & 400 &  90 & 400 & 127 & 500 & 146 &  500 & 220 \\
$\Re = 2.0\cdot 10^5$ & 500 & 109 & 600 & 164 & 600 & 201 &  700 & 257 \\
$\Re = 5.0\cdot 10^5$ & 800 & 164 & 900 & 239 & 900 & 276 & 1000 & 331 
\end{tabular}
\end{table}

\begin{figure}
\begin{center}
\input{apcfe-paper-figure-dists_NPmax_Re2e5_N_dH.txt}
\end{center}
\caption{\textit{Poiseuille flow, $\Re = 2\cdot10^{5}$}. The distances 
as described in \cref{eq:def_distances} for Poiseuille flow with $Re = 2\cdot10^5$ and different
values of $P$. We note that that increasing the number of polynomials does
not lead to a decrease of the error if the number of significant bits is not
chosen high enough.}
\label{fig:convana_Re2e5_N_dH}
\end{figure}

\begin{figure}
\begin{center}
\input{apcfe-paper-figure-dists_NPmax_Re5e5_N_dH.txt}
\end{center}
\caption{\textit{Poiseuille flow, $\Re = 5\cdot10^{5}$}. The distances as
described in \cref{eq:def_distances} for Poiseuille flow with $Re =
5\cdot10^5$ and different values of $P$. We note that that increasing the
number of polynomials does not lead to a decrease of the error if the number
of significant bits is not chosen high enough.}
\label{fig:convana_Re5e5_N_dH}
\end{figure}

\begin{figure}
\begin{center}
\input{apcfe-paper-figure-dists_NPmax_Re2e5_P_dH.txt}
\end{center}
\caption{\textit{Poiseuille flow, $\Re = 2\cdot 10^5$}. The distances
as described in \cref{eq:def_distances} for Poiseuille flow with $Re = 2\cdot10^5$ and different
values of $N$. We note that that increasing the number of significant bits
$P$ (or decreasing $\eps_P$ respectively) does not lead to a decrease of the error
if the number of polynomials is not chosen high enough.}
\label{fig:convana_Re2e5_P_dH}
\end{figure}

\begin{figure}
\begin{center}
\input{apcfe-paper-figure-dists_NPmax_Re5e5_P_dH.txt}
\end{center}
\caption{\textit{Poiseuille flow, $\Re = 5\cdot 10^5$}. The distances 
as described in \cref{eq:def_distances} for Poiseuille flow with $Re = 5\cdot10^5$ and different
values of $N$. We note that that increasing the number of significant bits
$P$ (or decreasing $\eps_P$ respectively) does not lead to a decrease of the error
if the number of polynomials is not chosen high enough.}
\label{fig:convana_Re5e5_P_dH}
\end{figure}

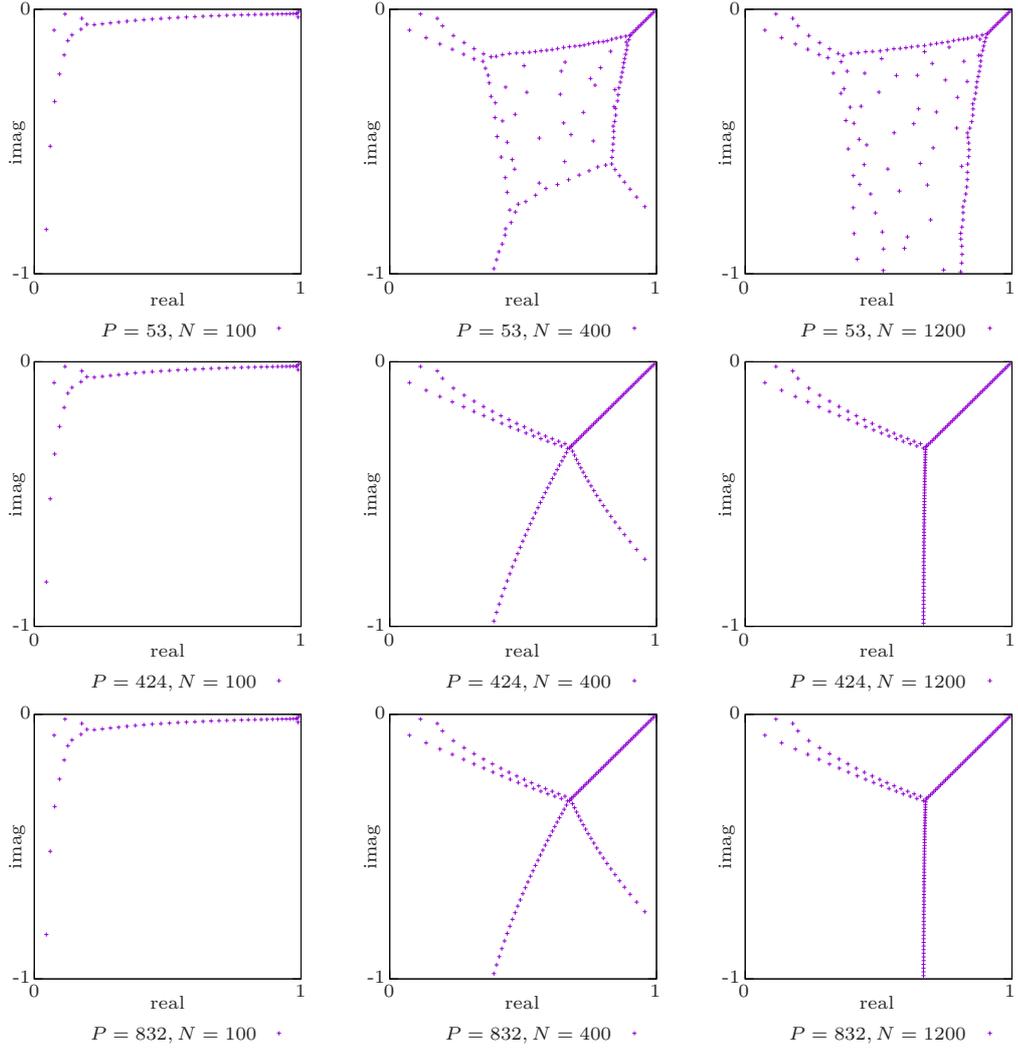
\begin{figure}
\begin{center}
\input{apcfe-paper-figure-out_pois_2e5_spectra.txt}
        \mbox{}\\
\end{center}
\caption{\textit{Poiseuille flow, $\Re = 2\cdot10^5$}. The Eigenvalues for
Poiseuille flow in the case $Re = 2\cdot10^5$ for different values of $N$ and
$P$.}
\end{figure}

\begin{figure}
\begin{center}
\input{apcfe-paper-figure-out_pois_5e5_spectra.txt}
        \mbox{}\\
\end{center}
\caption{\textit{Poiseuille flow, $\Re = 5\cdot10^5$}. The Eigenvalues for
Poiseuille flow in the case $Re = 5\cdot10^5$ for different values of $N$ and
$P$.}
\end{figure}
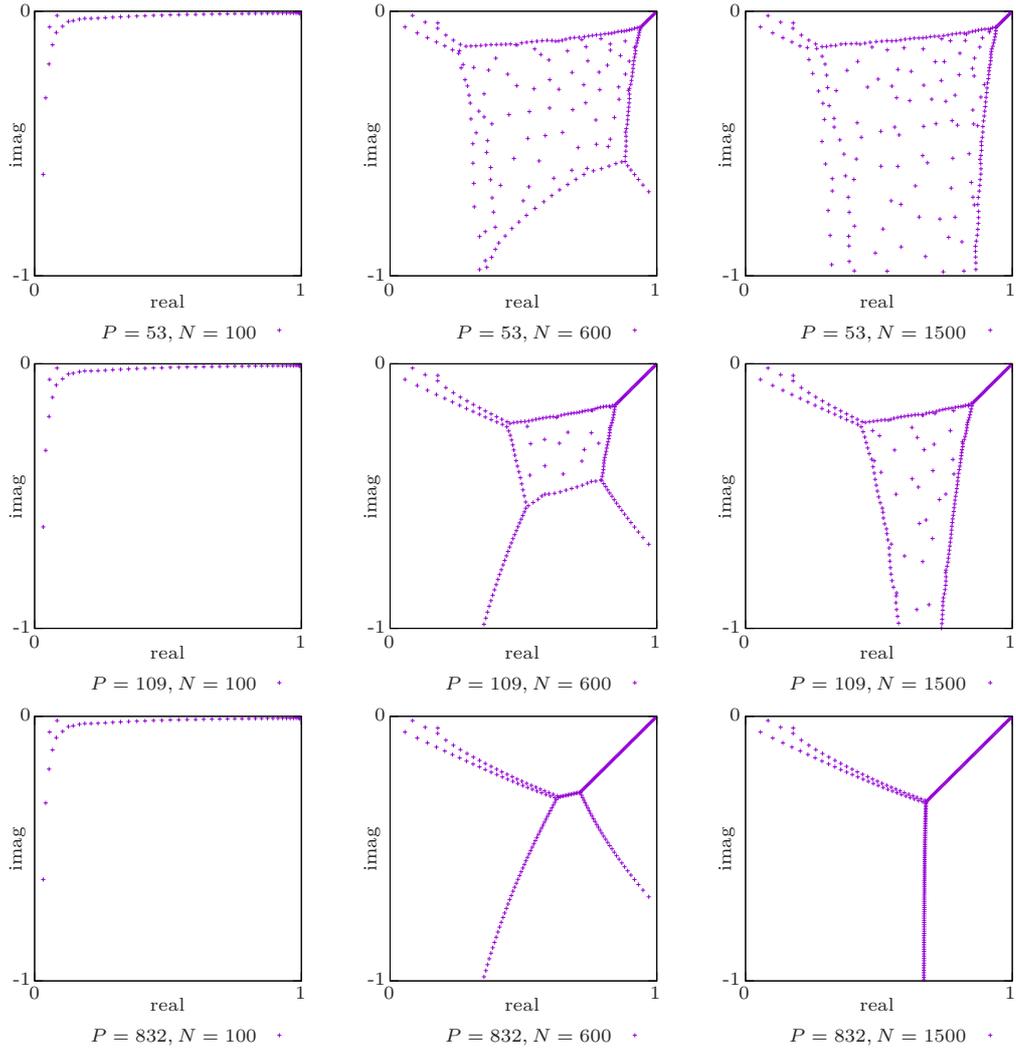

\subsection{The \texorpdfstring{$D^2$}{D2} method for Couette Flow}
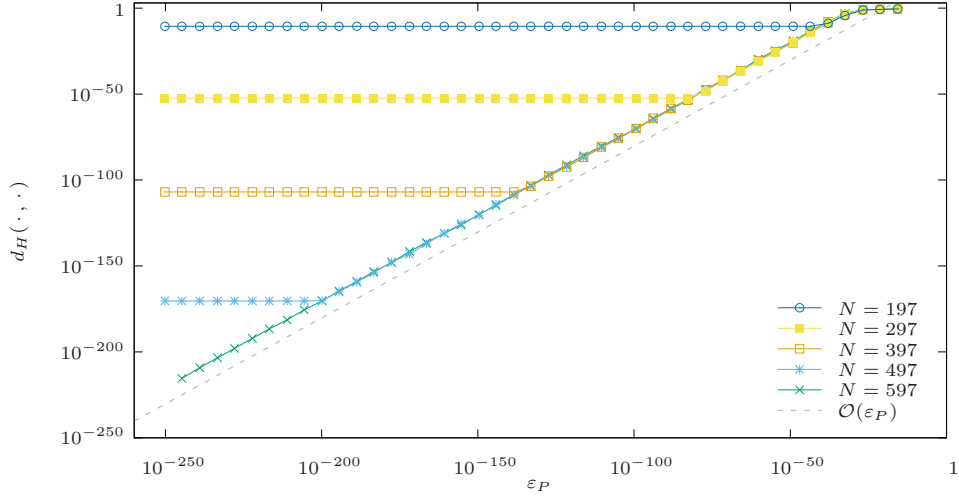
\begin{figure}
\begin{center}
\input{apcfe-paper-figure-dists_coue_NPmax_Re13000_P_dH.txt}
\end{center}
\caption{\textit{Couette flow, $\Re = 13,000$}. The distances 
as described in \cref{eq:def_distances} for Couette flow with $Re = 13,000$ and different
values of $N$. We note that that increasing the number of significant bits
$P$ (or decreasing $\eps_P$ respectively) does not lead to a decrease of the error
if the number of polynomials is not chosen high enough.}
\label{fig:convana_coue_Re13000_P_dH}
\end{figure}

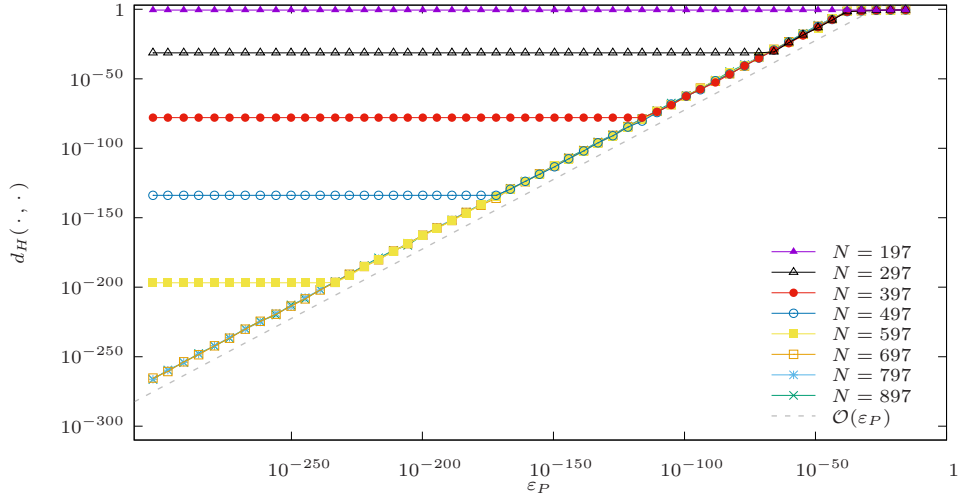
\begin{figure}
\begin{center}
\input{apcfe-paper-figure-dists_coue_NPmax_Re20000_P_dH.txt}
\end{center}
\caption{\textit{Couette flow, $\Re = 20,000$}. The distances 
as described in \cref{eq:def_distances} for Couette flow with $Re = 20,000$ and different
values of $N$. We note that that increasing the number of significant bits
$P$ (or decreasing $\eps_P$ respectively) does not lead to a decrease of the error
if the number of polynomials is not chosen high enough.}
\label{fig:convana_coue_Re20000_P_dH}
\end{figure}

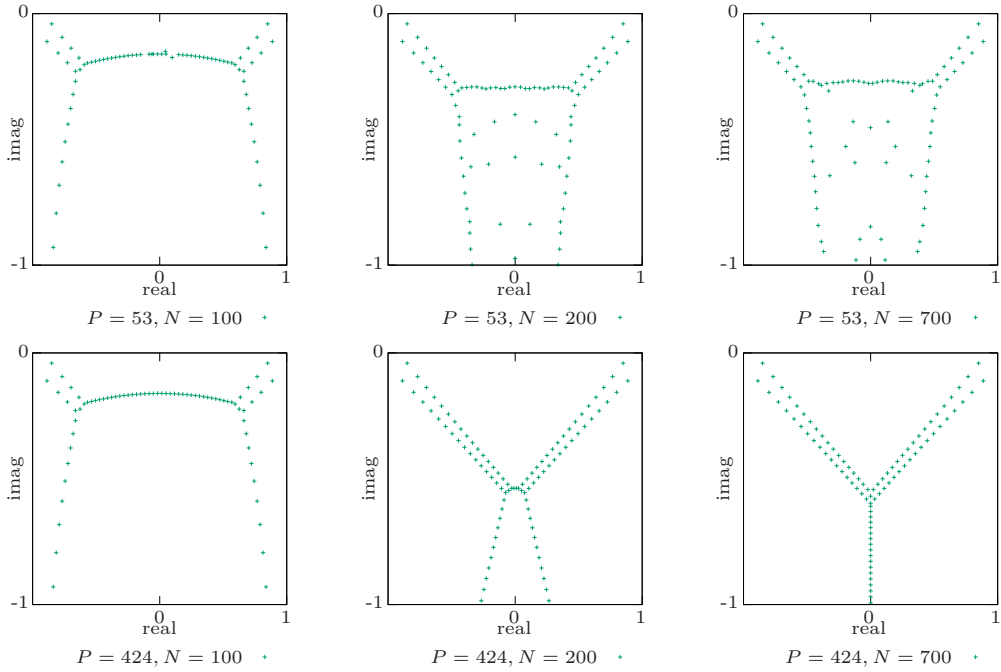
\begin{figure}
\begin{center}
\input{apcfe-paper-figure-out_coue_20000_spectra.txt}
        \mbox{}\\
\end{center}
\caption{\textit{Couette flow, $\Re = 20,000$}. The Eigenvalues for Couette flow in the case $Re = 20,000$ for
different values of $N$ and $P$.
}
\label{fig:coue_Re20000_spectrum}
\end{figure}

We also analyze the effects of the number of polynomials $N$ and the number of significant bits $P$ used in the Chebyshev tau method computation
for Couette flow as described in \cref{sub:D2_method}. Again, we thus compute the
spectra $\sigma_P(A^N, B^N)$ of the matrices as described there and compare them to
the spectra $\sigma_\Pmax(A^\Nmax, B^\Nmax)$ computed with maximal number of
polynomials $\Nmax$ and maximal number of significant bits $\Pmax$. In case of Couette
flow, we consider the eigenvalues restricted to the domain
\begin{equation*}
Q\coloneqq \{z\in\C\,|\, -1\leq \im(z)\leq 0,\, -1\leq \re(z)\leq 1 \}.
\end{equation*}
We show
the computed Hausdorff-distances for $\Re = 13,000$ in \cref{fig:convana_coue_Re13000_N_dH,fig:convana_coue_Re20000_N_dH}
and for $\Re = 20,000$ in \cref{fig:convana_coue_Re13000_P_dH,fig:convana_coue_Re20000_P_dH}.
Again, we notice that both $N$ and $P$ have to be chosen high enough. Finally,
we show the spectra that we compute for $\Re = 20,000$ and different values of $N$
and $P$ in \cref{fig:coue_Re20000_spectrum}.

\begin{figure}
\begin{center}
\input{apcfe-paper-figure-dists_coue_NPmax_Re13000_N_dH.txt}
\end{center}
\caption{\textit{Couette flow, $\Re = 13,000$}. The distances as described in \cref{eq:def_distances} for Couette flow with $Re = 13,000$ and different
values of $P$. We note that that increasing the number of polynomials does
not lead to a decrease of the error if the number of significant bits is not
chosen high enough.}
\label{fig:convana_coue_Re13000_N_dH}
\end{figure}
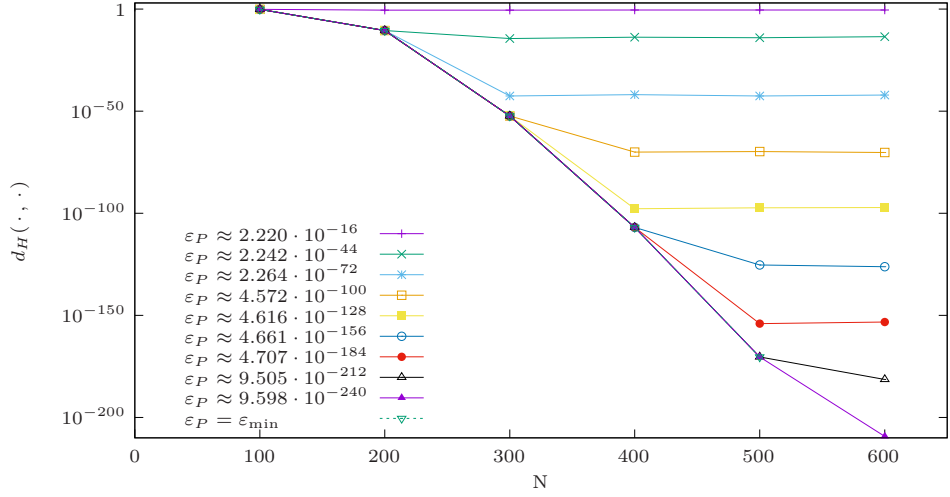

\begin{figure}
\begin{center}
\input{apcfe-paper-figure-dists_coue_NPmax_Re20000_N_dH.txt}
\end{center}
\caption{\textit{Couette flow, $\Re = 13,000$}. The distances 
as described in \cref{eq:def_distances} for Couette flow with $Re = 20,000$ and different
values of $P$. We note that that increasing the number of polynomials does
not lead to a decrease of the error if the number of significant bits is not
chosen high enough.}
\label{fig:convana_coue_Re20000_N_dH}
\end{figure}

\subsection{Comparison of the \texorpdfstring{$D^2$}{D2} and the \texorpdfstring{$D^4$}{D4} method}
\begin{figure}
\begin{center}
\input{apcfe-paper-figure-dists_comp_D2D4_P_832.txt}
\end{center}
\caption{Comparison of the results we obtain using the $D^2$ method and the $D^4$ method
for Poiseuille flow with $\Re = 10^{5}$ and $P = 832$. We use the spectrum computed using the $D^2$ method with $\Pmax,\Nmax$ as reference. The $D^4$ method is now equivalent in accuracy to the $D^2$-method. Furthermore the $D^4$ method is still faster for the 
same values of $N$.}
\label{fig:comp_D2D4_P_832}
\end{figure}
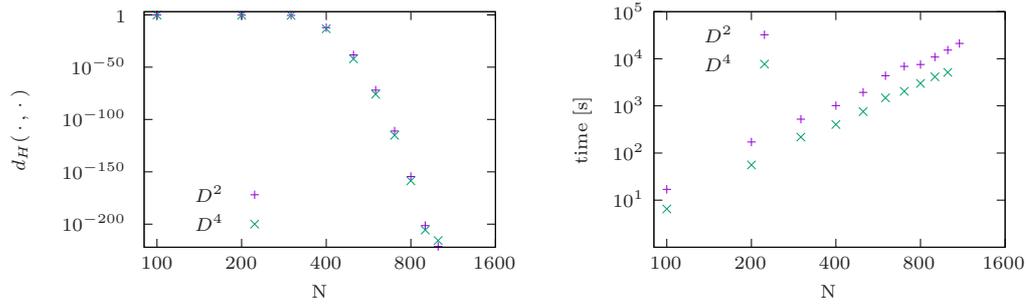

As discussed in \cref{sub:D4_method}, the use of high precision floating point
numbers may make the $D^4$ method feasible again. In this section, we compare the results we obtain using the $D^2$ method and the $D^4$ method.
Thus, we compute the Eigenvalues for Poiseuille flow using both methods for 
 fixed Reynolds number $\Re = 10^5$ and use different values of $N$ and $P$.
We use the spectrum that we obtained using the $D^2$ method with $\Nmax = 1000$
and $\Pmax = 832$ as the reference and compute the Hausdorff distances as described
before for the results we obtain using $D^2$ method and $D^4$ method.
We show the results for selected values of $P$ in
\cref{fig:comp_D2D4_P_53,fig:comp_D2D4_P_832} and see that
the error, measured as usual by the Hausdorff distance to the reference computation, is lower for the same values of $N$ when we use the $D^2$ method
and small values of $P$ such as, e.g., $P = 53$. However, the situation changes
if we use high values of $P$ such as $P = 832$. Furthermore, the actual wall time
needed to compute the results is significantly lower in the case of the $D^4$ method, as the matrix size is only half that of the $D^2$ method. 

\begin{figure}
\begin{center}
\input{apcfe-paper-figure-dists_comp_D2D4_P_53.txt}
\end{center}
\caption{Comparison of the results we obtain using the $D^2$ method and the $D^4$ method
for Poiseuille flow with $\Re = 10^{5}$ and $P = 53$. We use the spectrum computed using the $D^2$ method with $\Pmax,\Nmax$ as reference. The $D^2$ method
yields much better accuracy this low value of $P$ (equivalent to standard \texttt{double} precision). The $D^4$ method is faster for the 
same values of $N$ as the respective matrix sizes are halved.}
\label{fig:comp_D2D4_P_53}
\end{figure}
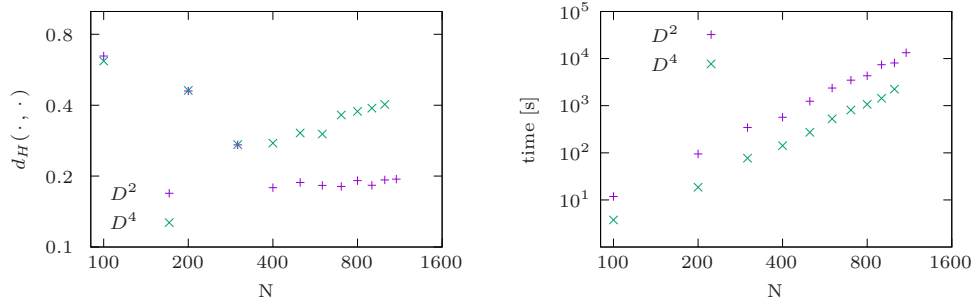

\section*{Funding}
The work of BS was supported by an Emeritus Fellowship of the
Leverhulme Trust, EM-2019-022/9.
\bibliography{bibliography,thompson}
\end{document}

%% file: apcfe-paper-figure-convana_7x7.txt
\begingroup
  \makeatletter
  \providecommand\color[2][]{%
    \GenericError{(gnuplot) \space\space\space\@spaces}{%
      Package color not loaded in conjunction with
      terminal option `colourtext'%
    }{See the gnuplot documentation for explanation.%
    }{Either use 'blacktext' in gnuplot or load the package
      color.sty in LaTeX.}%
    \renewcommand\color[2][]{}%
  }%
  \providecommand\includegraphics[2][]{%
    \GenericError{(gnuplot) \space\space\space\@spaces}{%
      Package graphicx or graphics not loaded%
    }{See the gnuplot documentation for explanation.%
    }{The gnuplot epslatex terminal needs graphicx.sty or graphics.sty.}%
    \renewcommand\includegraphics[2][]{}%
  }%
  \providecommand\rotatebox[2]{#2}%
  \@ifundefined{ifGPcolor}{%
    \newif\ifGPcolor
    \GPcolorfalse
  }{}%
  \@ifundefined{ifGPblacktext}{%
    \newif\ifGPblacktext
    \GPblacktexttrue
  }{}%
  \let\gplgaddtomacro\g@addto@macro
  \gdef\gplbacktext{}%
  \gdef\gplfronttext{}%
  \makeatother
  \ifGPblacktext
    \def\colorrgb#1{}%
    \def\colorgray#1{}%
  \else
    \ifGPcolor
      \def\colorrgb#1{\color[rgb]{#1}}%
      \def\colorgray#1{\color[gray]{#1}}%
      \expandafter\def\csname LTw\endcsname{\color{white}}%
      \expandafter\def\csname LTb\endcsname{\color{black}}%
      \expandafter\def\csname LTa\endcsname{\color{black}}%
      \expandafter\def\csname LT0\endcsname{\color[rgb]{1,0,0}}%
      \expandafter\def\csname LT1\endcsname{\color[rgb]{0,1,0}}%
      \expandafter\def\csname LT2\endcsname{\color[rgb]{0,0,1}}%
      \expandafter\def\csname LT3\endcsname{\color[rgb]{1,0,1}}%
      \expandafter\def\csname LT4\endcsname{\color[rgb]{0,1,1}}%
      \expandafter\def\csname LT5\endcsname{\color[rgb]{1,1,0}}%
      \expandafter\def\csname LT6\endcsname{\color[rgb]{0,0,0}}%
      \expandafter\def\csname LT7\endcsname{\color[rgb]{1,0.3,0}}%
      \expandafter\def\csname LT8\endcsname{\color[rgb]{0.5,0.5,0.5}}%
    \else
      \def\colorrgb#1{\color{black}}%
      \def\colorgray#1{\color[gray]{#1}}%
      \expandafter\def\csname LTw\endcsname{\color{white}}%
      \expandafter\def\csname LTb\endcsname{\color{black}}%
      \expandafter\def\csname LTa\endcsname{\color{black}}%
      \expandafter\def\csname LT0\endcsname{\color{black}}%
      \expandafter\def\csname LT1\endcsname{\color{black}}%
      \expandafter\def\csname LT2\endcsname{\color{black}}%
      \expandafter\def\csname LT3\endcsname{\color{black}}%
      \expandafter\def\csname LT4\endcsname{\color{black}}%
      \expandafter\def\csname LT5\endcsname{\color{black}}%
      \expandafter\def\csname LT6\endcsname{\color{black}}%
      \expandafter\def\csname LT7\endcsname{\color{black}}%
      \expandafter\def\csname LT8\endcsname{\color{black}}%
    \fi
  \fi
    \setlength{\unitlength}{0.0500bp}%
    \ifx\gptboxheight\undefined%
      \newlength{\gptboxheight}%
      \newlength{\gptboxwidth}%
      \newsavebox{\gptboxtext}%
    \fi%
    \setlength{\fboxrule}{0.5pt}%
    \setlength{\fboxsep}{1pt}%
    \definecolor{tbcol}{rgb}{1,1,1}%
\begin{picture}(7370.00,2834.00)%
    \gplgaddtomacro\gplbacktext{%
      \csname LTb\endcsname
      \put(744,2651){\makebox(0,0)[r]{\strut{}\footnotesize{$1$}}}%
      \put(744,566){\makebox(0,0)[r]{\strut{}\footnotesize{$10^{-50}$}}}%
      \put(744,983){\makebox(0,0)[r]{\strut{}\footnotesize{$10^{-40}$}}}%
      \put(744,1400){\makebox(0,0)[r]{\strut{}\footnotesize{$10^{-30}$}}}%
      \put(744,1817){\makebox(0,0)[r]{\strut{}\footnotesize{$10^{-20}$}}}%
      \put(744,2234){\makebox(0,0)[r]{\strut{}\footnotesize{$10^{-10}$}}}%
      \put(3426,456){\makebox(0,0){\strut{}\footnotesize{$1$}}}%
      \put(810,456){\makebox(0,0){\strut{}\footnotesize{$10^{-120}$}}}%
      \put(1682,456){\makebox(0,0){\strut{}\footnotesize{$10^{-80}$}}}%
      \put(2554,456){\makebox(0,0){\strut{}\footnotesize{$10^{-40}$}}}%
    }%
    \gplgaddtomacro\gplfronttext{%
      \csname LTb\endcsname
      \put(1527,2541){\makebox(0,0)[r]{\strut{}\footnotesize{Schur}}}%
      \csname LTb\endcsname
      \put(1527,2321){\makebox(0,0)[r]{\strut{}\footnotesize{QZ}}}%
      \csname LTb\endcsname
      \put(1527,2101){\makebox(0,0)[r]{\strut{}\footnotesize{$\mathcal O\big(\eps_P^{1/2}\big)$}}}%
      \csname LTb\endcsname
      \put(73,1671){\rotatebox{-270.00}{\makebox(0,0){\strut{}\footnotesize{$d_H\left(\sigma_P(A^{7\times 7}_1), \sigma(A^{7\times 7}_1)\right)$}}}}%
      \put(2118,236){\makebox(0,0){\strut{}\footnotesize{$\eps_P$}}}%
    }%
    \gplgaddtomacro\gplbacktext{%
      \csname LTb\endcsname
      \put(4540,2712){\makebox(0,0)[r]{\strut{}\footnotesize{$1$}}}%
      \put(4540,566){\makebox(0,0)[r]{\strut{}\footnotesize{$10^{-100}$}}}%
      \put(4540,995){\makebox(0,0)[r]{\strut{}\footnotesize{$10^{-80}$}}}%
      \put(4540,1424){\makebox(0,0)[r]{\strut{}\footnotesize{$10^{-60}$}}}%
      \put(4540,1853){\makebox(0,0)[r]{\strut{}\footnotesize{$10^{-40}$}}}%
      \put(4540,2283){\makebox(0,0)[r]{\strut{}\footnotesize{$10^{-20}$}}}%
      \put(7221,456){\makebox(0,0){\strut{}\footnotesize{$1$}}}%
      \put(4606,456){\makebox(0,0){\strut{}\footnotesize{$10^{-120}$}}}%
      \put(5478,456){\makebox(0,0){\strut{}\footnotesize{$10^{-80}$}}}%
      \put(6349,456){\makebox(0,0){\strut{}\footnotesize{$10^{-40}$}}}%
    }%
    \gplgaddtomacro\gplfronttext{%
      \csname LTb\endcsname
      \put(5323,2602){\makebox(0,0)[r]{\strut{}\footnotesize{Schur}}}%
      \csname LTb\endcsname
      \put(5323,2382){\makebox(0,0)[r]{\strut{}\footnotesize{QZ}}}%
      \csname LTb\endcsname
      \put(5323,2162){\makebox(0,0)[r]{\strut{}\footnotesize{$\mathcal O\big(\eps_P\big)$}}}%
      \csname LTb\endcsname
      \put(3869,1671){\rotatebox{-270.00}{\makebox(0,0){\strut{}\footnotesize{$d_H\left(\sigma_P(A_{-1}^{7\times 7}), \sigma(A_{-1}^{7\times 7})\right)$}}}}%
      \put(5913,236){\makebox(0,0){\strut{}\footnotesize{$\eps_P$}}}%
    }%
    \gplbacktext
    \put(0,0){\includegraphics[width={368.50bp},height={141.70bp}]{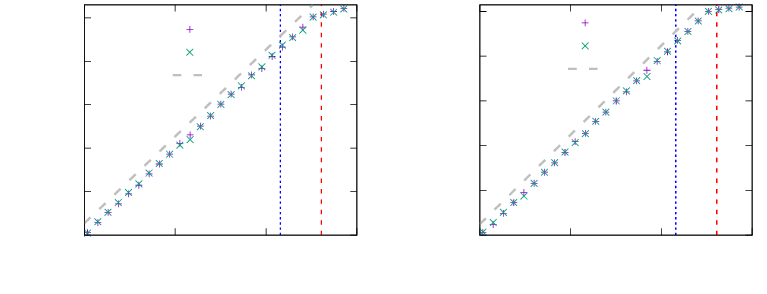}}%
    \gplfronttext
  \end{picture}%
\endgroup

%% file: apcfe-paper-figure-7x7_spects.txt
\begingroup
  \makeatletter
  \providecommand\color[2][]{%
    \GenericError{(gnuplot) \space\space\space\@spaces}{%
      Package color not loaded in conjunction with
      terminal option `colourtext'%
    }{See the gnuplot documentation for explanation.%
    }{Either use 'blacktext' in gnuplot or load the package
      color.sty in LaTeX.}%
    \renewcommand\color[2][]{}%
  }%
  \providecommand\includegraphics[2][]{%
    \GenericError{(gnuplot) \space\space\space\@spaces}{%
      Package graphicx or graphics not loaded%
    }{See the gnuplot documentation for explanation.%
    }{The gnuplot epslatex terminal needs graphicx.sty or graphics.sty.}%
    \renewcommand\includegraphics[2][]{}%
  }%
  \providecommand\rotatebox[2]{#2}%
  \@ifundefined{ifGPcolor}{%
    \newif\ifGPcolor
    \GPcolorfalse
  }{}%
  \@ifundefined{ifGPblacktext}{%
    \newif\ifGPblacktext
    \GPblacktexttrue
  }{}%
  \let\gplgaddtomacro\g@addto@macro
  \gdef\gplbacktext{}%
  \gdef\gplfronttext{}%
  \makeatother
  \ifGPblacktext
    \def\colorrgb#1{}%
    \def\colorgray#1{}%
  \else
    \ifGPcolor
      \def\colorrgb#1{\color[rgb]{#1}}%
      \def\colorgray#1{\color[gray]{#1}}%
      \expandafter\def\csname LTw\endcsname{\color{white}}%
      \expandafter\def\csname LTb\endcsname{\color{black}}%
      \expandafter\def\csname LTa\endcsname{\color{black}}%
      \expandafter\def\csname LT0\endcsname{\color[rgb]{1,0,0}}%
      \expandafter\def\csname LT1\endcsname{\color[rgb]{0,1,0}}%
      \expandafter\def\csname LT2\endcsname{\color[rgb]{0,0,1}}%
      \expandafter\def\csname LT3\endcsname{\color[rgb]{1,0,1}}%
      \expandafter\def\csname LT4\endcsname{\color[rgb]{0,1,1}}%
      \expandafter\def\csname LT5\endcsname{\color[rgb]{1,1,0}}%
      \expandafter\def\csname LT6\endcsname{\color[rgb]{0,0,0}}%
      \expandafter\def\csname LT7\endcsname{\color[rgb]{1,0.3,0}}%
      \expandafter\def\csname LT8\endcsname{\color[rgb]{0.5,0.5,0.5}}%
    \else
      \def\colorrgb#1{\color{black}}%
      \def\colorgray#1{\color[gray]{#1}}%
      \expandafter\def\csname LTw\endcsname{\color{white}}%
      \expandafter\def\csname LTb\endcsname{\color{black}}%
      \expandafter\def\csname LTa\endcsname{\color{black}}%
      \expandafter\def\csname LT0\endcsname{\color{black}}%
      \expandafter\def\csname LT1\endcsname{\color{black}}%
      \expandafter\def\csname LT2\endcsname{\color{black}}%
      \expandafter\def\csname LT3\endcsname{\color{black}}%
      \expandafter\def\csname LT4\endcsname{\color{black}}%
      \expandafter\def\csname LT5\endcsname{\color{black}}%
      \expandafter\def\csname LT6\endcsname{\color{black}}%
      \expandafter\def\csname LT7\endcsname{\color{black}}%
      \expandafter\def\csname LT8\endcsname{\color{black}}%
    \fi
  \fi
    \setlength{\unitlength}{0.0500bp}%
    \ifx\gptboxheight\undefined%
      \newlength{\gptboxheight}%
      \newlength{\gptboxwidth}%
      \newsavebox{\gptboxtext}%
    \fi%
    \setlength{\fboxrule}{0.5pt}%
    \setlength{\fboxsep}{1pt}%
    \definecolor{tbcol}{rgb}{1,1,1}%
\begin{picture}(7370.00,4534.00)%
    \gplgaddtomacro\gplbacktext{%
      \csname LTb\endcsname
      \put(678,1547){\makebox(0,0)[r]{\strut{}\footnotesize{-6}}}%
      \put(678,1921){\makebox(0,0)[r]{\strut{}\footnotesize{-4}}}%
      \put(678,2294){\makebox(0,0)[r]{\strut{}\footnotesize{-2}}}%
      \put(678,2668){\makebox(0,0)[r]{\strut{}\footnotesize{0}}}%
      \put(678,3042){\makebox(0,0)[r]{\strut{}\footnotesize{2}}}%
      \put(678,3415){\makebox(0,0)[r]{\strut{}\footnotesize{4}}}%
      \put(678,3789){\makebox(0,0)[r]{\strut{}\footnotesize{6}}}%
      \put(997,1140){\makebox(0,0){\strut{}\footnotesize{-6}}}%
      \put(1371,1140){\makebox(0,0){\strut{}\footnotesize{-4}}}%
      \put(1744,1140){\makebox(0,0){\strut{}\footnotesize{-2}}}%
      \put(2118,1140){\makebox(0,0){\strut{}\footnotesize{0}}}%
      \put(2492,1140){\makebox(0,0){\strut{}\footnotesize{2}}}%
      \put(2865,1140){\makebox(0,0){\strut{}\footnotesize{4}}}%
      \put(3239,1140){\makebox(0,0){\strut{}\footnotesize{6}}}%
    }%
    \gplgaddtomacro\gplfronttext{%
      \csname LTb\endcsname
      \put(1881,503){\makebox(0,0)[r]{\strut{}\footnotesize{true spectrum}}}%
      \csname LTb\endcsname
      \put(1881,283){\makebox(0,0)[r]{\strut{}\footnotesize{$\eps_P \approx 2.220 \times 10^{-16}$}}}%
      \csname LTb\endcsname
      \put(4188,503){\makebox(0,0)[r]{\strut{}\footnotesize{$\eps_P \approx 3.469 \times 10^{-18}$}}}%
      \csname LTb\endcsname
      \put(4188,283){\makebox(0,0)[r]{\strut{}\footnotesize{$\eps_P \approx 5.421 \times 10^{-20}$}}}%
      \csname LTb\endcsname
      \put(6495,503){\makebox(0,0)[r]{\strut{}\footnotesize{$\eps_P \approx 8.470 \times 10^{-22}$}}}%
      \csname LTb\endcsname
      \put(6495,283){\makebox(0,0)[r]{\strut{}\footnotesize{$\eps_P \approx 2.465 \times 10^{-32}$}}}%
      \csname LTb\endcsname
      \put(403,2668){\rotatebox{-270.00}{\makebox(0,0){\strut{}\footnotesize{imag}}}}%
      \put(2118,920){\makebox(0,0){\strut{}\footnotesize{real}}}%
    }%
    \gplgaddtomacro\gplbacktext{%
      \csname LTb\endcsname
      \put(4474,1547){\makebox(0,0)[r]{\strut{}\footnotesize{-6}}}%
      \put(4474,1920){\makebox(0,0)[r]{\strut{}\footnotesize{-4}}}%
      \put(4474,2294){\makebox(0,0)[r]{\strut{}\footnotesize{-2}}}%
      \put(4474,2668){\makebox(0,0)[r]{\strut{}\footnotesize{0}}}%
      \put(4474,3041){\makebox(0,0)[r]{\strut{}\footnotesize{2}}}%
      \put(4474,3415){\makebox(0,0)[r]{\strut{}\footnotesize{4}}}%
      \put(4474,3788){\makebox(0,0)[r]{\strut{}\footnotesize{6}}}%
      \put(4793,1140){\makebox(0,0){\strut{}\footnotesize{-6}}}%
      \put(5166,1140){\makebox(0,0){\strut{}\footnotesize{-4}}}%
      \put(5540,1140){\makebox(0,0){\strut{}\footnotesize{-2}}}%
      \put(5914,1140){\makebox(0,0){\strut{}\footnotesize{0}}}%
      \put(6287,1140){\makebox(0,0){\strut{}\footnotesize{2}}}%
      \put(6661,1140){\makebox(0,0){\strut{}\footnotesize{4}}}%
      \put(7034,1140){\makebox(0,0){\strut{}\footnotesize{6}}}%
    }%
    \gplgaddtomacro\gplfronttext{%
      \csname LTb\endcsname
      \put(4199,2667){\rotatebox{-270.00}{\makebox(0,0){\strut{}\footnotesize{imag}}}}%
      \put(5913,920){\makebox(0,0){\strut{}\footnotesize{real}}}%
    }%
    \gplbacktext
    \put(0,0){\includegraphics[width={368.50bp},height={226.70bp}]{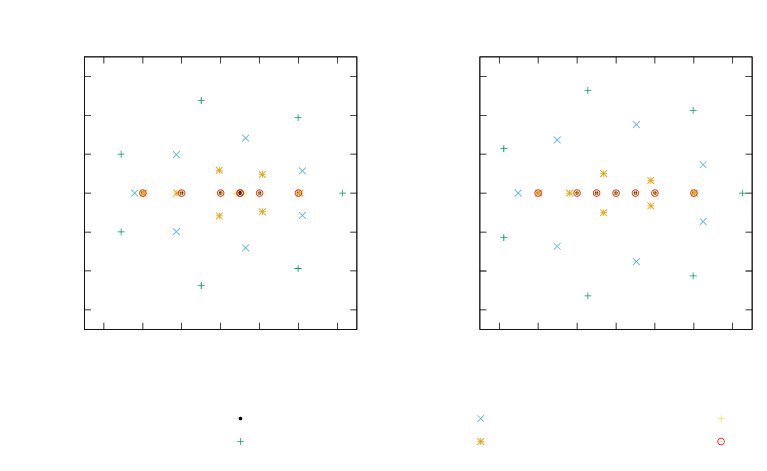}}%
    \gplfronttext
  \end{picture}%
\endgroup

%% file: apcfe-paper-figure-Straughan-paper-repl_Straugh_Fig1.txt
\begingroup
  \makeatletter
  \providecommand\color[2][]{%
    \GenericError{(gnuplot) \space\space\space\@spaces}{%
      Package color not loaded in conjunction with
      terminal option `colourtext'%
    }{See the gnuplot documentation for explanation.%
    }{Either use 'blacktext' in gnuplot or load the package
      color.sty in LaTeX.}%
    \renewcommand\color[2][]{}%
  }%
  \providecommand\includegraphics[2][]{%
    \GenericError{(gnuplot) \space\space\space\@spaces}{%
      Package graphicx or graphics not loaded%
    }{See the gnuplot documentation for explanation.%
    }{The gnuplot epslatex terminal needs graphicx.sty or graphics.sty.}%
    \renewcommand\includegraphics[2][]{}%
  }%
  \providecommand\rotatebox[2]{#2}%
  \@ifundefined{ifGPcolor}{%
    \newif\ifGPcolor
    \GPcolorfalse
  }{}%
  \@ifundefined{ifGPblacktext}{%
    \newif\ifGPblacktext
    \GPblacktexttrue
  }{}%
  \let\gplgaddtomacro\g@addto@macro
  \gdef\gplbacktext{}%
  \gdef\gplfronttext{}%
  \makeatother
  \ifGPblacktext
    \def\colorrgb#1{}%
    \def\colorgray#1{}%
  \else
    \ifGPcolor
      \def\colorrgb#1{\color[rgb]{#1}}%
      \def\colorgray#1{\color[gray]{#1}}%
      \expandafter\def\csname LTw\endcsname{\color{white}}%
      \expandafter\def\csname LTb\endcsname{\color{black}}%
      \expandafter\def\csname LTa\endcsname{\color{black}}%
      \expandafter\def\csname LT0\endcsname{\color[rgb]{1,0,0}}%
      \expandafter\def\csname LT1\endcsname{\color[rgb]{0,1,0}}%
      \expandafter\def\csname LT2\endcsname{\color[rgb]{0,0,1}}%
      \expandafter\def\csname LT3\endcsname{\color[rgb]{1,0,1}}%
      \expandafter\def\csname LT4\endcsname{\color[rgb]{0,1,1}}%
      \expandafter\def\csname LT5\endcsname{\color[rgb]{1,1,0}}%
      \expandafter\def\csname LT6\endcsname{\color[rgb]{0,0,0}}%
      \expandafter\def\csname LT7\endcsname{\color[rgb]{1,0.3,0}}%
      \expandafter\def\csname LT8\endcsname{\color[rgb]{0.5,0.5,0.5}}%
    \else
      \def\colorrgb#1{\color{black}}%
      \def\colorgray#1{\color[gray]{#1}}%
      \expandafter\def\csname LTw\endcsname{\color{white}}%
      \expandafter\def\csname LTb\endcsname{\color{black}}%
      \expandafter\def\csname LTa\endcsname{\color{black}}%
      \expandafter\def\csname LT0\endcsname{\color{black}}%
      \expandafter\def\csname LT1\endcsname{\color{black}}%
      \expandafter\def\csname LT2\endcsname{\color{black}}%
      \expandafter\def\csname LT3\endcsname{\color{black}}%
      \expandafter\def\csname LT4\endcsname{\color{black}}%
      \expandafter\def\csname LT5\endcsname{\color{black}}%
      \expandafter\def\csname LT6\endcsname{\color{black}}%
      \expandafter\def\csname LT7\endcsname{\color{black}}%
      \expandafter\def\csname LT8\endcsname{\color{black}}%
    \fi
  \fi
    \setlength{\unitlength}{0.0500bp}%
    \ifx\gptboxheight\undefined%
      \newlength{\gptboxheight}%
      \newlength{\gptboxwidth}%
      \newsavebox{\gptboxtext}%
    \fi%
    \setlength{\fboxrule}{0.5pt}%
    \setlength{\fboxsep}{1pt}%
    \definecolor{tbcol}{rgb}{1,1,1}%
\begin{picture}(4534.00,2834.00)%
    \gplgaddtomacro\gplbacktext{%
      \csname LTb\endcsname
      \put(1157,440){\makebox(0,0)[r]{\strut{}\footnotesize{-1}}}%
      \put(1157,647){\makebox(0,0)[r]{\strut{}\footnotesize{-0.9}}}%
      \put(1157,854){\makebox(0,0)[r]{\strut{}\footnotesize{-0.8}}}%
      \put(1157,1061){\makebox(0,0)[r]{\strut{}\footnotesize{-0.7}}}%
      \put(1157,1268){\makebox(0,0)[r]{\strut{}\footnotesize{-0.6}}}%
      \put(1157,1475){\makebox(0,0)[r]{\strut{}\footnotesize{-0.5}}}%
      \put(1157,1682){\makebox(0,0)[r]{\strut{}\footnotesize{-0.4}}}%
      \put(1157,1889){\makebox(0,0)[r]{\strut{}\footnotesize{-0.3}}}%
      \put(1157,2096){\makebox(0,0)[r]{\strut{}\footnotesize{-0.2}}}%
      \put(1157,2303){\makebox(0,0)[r]{\strut{}\footnotesize{-0.1}}}%
      \put(1157,2510){\makebox(0,0)[r]{\strut{}\footnotesize{0}}}%
      \put(1210,330){\makebox(0,0){\strut{}\footnotesize{0}}}%
      \put(1503,330){\makebox(0,0){\strut{}\footnotesize{0.1}}}%
      \put(1795,330){\makebox(0,0){\strut{}\footnotesize{0.2}}}%
      \put(2088,330){\makebox(0,0){\strut{}\footnotesize{0.3}}}%
      \put(2381,330){\makebox(0,0){\strut{}\footnotesize{0.4}}}%
      \put(2674,330){\makebox(0,0){\strut{}\footnotesize{0.5}}}%
      \put(2966,330){\makebox(0,0){\strut{}\footnotesize{0.6}}}%
      \put(3259,330){\makebox(0,0){\strut{}\footnotesize{0.7}}}%
      \put(3552,330){\makebox(0,0){\strut{}\footnotesize{0.8}}}%
      \put(3844,330){\makebox(0,0){\strut{}\footnotesize{0.9}}}%
      \put(4137,330){\makebox(0,0){\strut{}\footnotesize{1}}}%
    }%
    \gplgaddtomacro\gplfronttext{%
      \csname LTb\endcsname
      \put(737,1526){\rotatebox{-270}{\makebox(0,0){\strut{}\footnotesize{imag}}}}%
      \put(2673,154){\makebox(0,0){\strut{}\footnotesize{real}}}%
      \csname LTb\endcsname
      \put(1748,758){\makebox(0,0)[r]{\strut{}\footnotesize{even}}}%
      \csname LTb\endcsname
      \put(1748,648){\makebox(0,0)[r]{\strut{}\footnotesize{odd}}}%
    }%
    \gplbacktext
    \put(0,0){\includegraphics[width={226.70bp},height={141.70bp}]{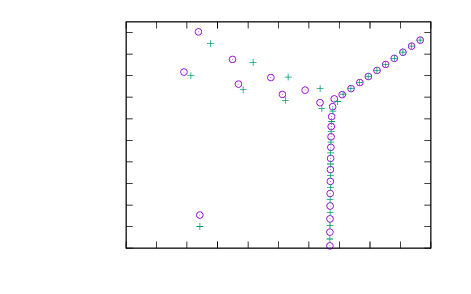}}%
    \gplfronttext
  \end{picture}%
\endgroup

%% file: apcfe-paper-figure-dists_NPmax_Re1e5_N_dH.txt
\begingroup
  \makeatletter
  \providecommand\color[2][]{%
    \GenericError{(gnuplot) \space\space\space\@spaces}{%
      Package color not loaded in conjunction with
      terminal option `colourtext'%
    }{See the gnuplot documentation for explanation.%
    }{Either use 'blacktext' in gnuplot or load the package
      color.sty in LaTeX.}%
    \renewcommand\color[2][]{}%
  }%
  \providecommand\includegraphics[2][]{%
    \GenericError{(gnuplot) \space\space\space\@spaces}{%
      Package graphicx or graphics not loaded%
    }{See the gnuplot documentation for explanation.%
    }{The gnuplot epslatex terminal needs graphicx.sty or graphics.sty.}%
    \renewcommand\includegraphics[2][]{}%
  }%
  \providecommand\rotatebox[2]{#2}%
  \@ifundefined{ifGPcolor}{%
    \newif\ifGPcolor
    \GPcolorfalse
  }{}%
  \@ifundefined{ifGPblacktext}{%
    \newif\ifGPblacktext
    \GPblacktexttrue
  }{}%
  \let\gplgaddtomacro\g@addto@macro
  \gdef\gplbacktext{}%
  \gdef\gplfronttext{}%
  \makeatother
  \ifGPblacktext
    \def\colorrgb#1{}%
    \def\colorgray#1{}%
  \else
    \ifGPcolor
      \def\colorrgb#1{\color[rgb]{#1}}%
      \def\colorgray#1{\color[gray]{#1}}%
      \expandafter\def\csname LTw\endcsname{\color{white}}%
      \expandafter\def\csname LTb\endcsname{\color{black}}%
      \expandafter\def\csname LTa\endcsname{\color{black}}%
      \expandafter\def\csname LT0\endcsname{\color[rgb]{1,0,0}}%
      \expandafter\def\csname LT1\endcsname{\color[rgb]{0,1,0}}%
      \expandafter\def\csname LT2\endcsname{\color[rgb]{0,0,1}}%
      \expandafter\def\csname LT3\endcsname{\color[rgb]{1,0,1}}%
      \expandafter\def\csname LT4\endcsname{\color[rgb]{0,1,1}}%
      \expandafter\def\csname LT5\endcsname{\color[rgb]{1,1,0}}%
      \expandafter\def\csname LT6\endcsname{\color[rgb]{0,0,0}}%
      \expandafter\def\csname LT7\endcsname{\color[rgb]{1,0.3,0}}%
      \expandafter\def\csname LT8\endcsname{\color[rgb]{0.5,0.5,0.5}}%
    \else
      \def\colorrgb#1{\color{black}}%
      \def\colorgray#1{\color[gray]{#1}}%
      \expandafter\def\csname LTw\endcsname{\color{white}}%
      \expandafter\def\csname LTb\endcsname{\color{black}}%
      \expandafter\def\csname LTa\endcsname{\color{black}}%
      \expandafter\def\csname LT0\endcsname{\color{black}}%
      \expandafter\def\csname LT1\endcsname{\color{black}}%
      \expandafter\def\csname LT2\endcsname{\color{black}}%
      \expandafter\def\csname LT3\endcsname{\color{black}}%
      \expandafter\def\csname LT4\endcsname{\color{black}}%
      \expandafter\def\csname LT5\endcsname{\color{black}}%
      \expandafter\def\csname LT6\endcsname{\color{black}}%
      \expandafter\def\csname LT7\endcsname{\color{black}}%
      \expandafter\def\csname LT8\endcsname{\color{black}}%
    \fi
  \fi
    \setlength{\unitlength}{0.0500bp}%
    \ifx\gptboxheight\undefined%
      \newlength{\gptboxheight}%
      \newlength{\gptboxwidth}%
      \newsavebox{\gptboxtext}%
    \fi%
    \setlength{\fboxrule}{0.5pt}%
    \setlength{\fboxsep}{1pt}%
    \definecolor{tbcol}{rgb}{1,1,1}%
\begin{picture}(7370.00,3968.00)%
    \gplgaddtomacro\gplbacktext{%
      \csname LTb\endcsname
      \put(990,3703){\makebox(0,0)[r]{\strut{}\footnotesize{$1$}}}%
      \put(990,803){\makebox(0,0)[r]{\strut{}\footnotesize{$10^{-200}$}}}%
      \put(990,1528){\makebox(0,0)[r]{\strut{}\footnotesize{$10^{-150}$}}}%
      \put(990,2253){\makebox(0,0)[r]{\strut{}\footnotesize{$10^{-100}$}}}%
      \put(990,2978){\makebox(0,0)[r]{\strut{}\footnotesize{$10^{-50}$}}}%
      \put(1056,352){\makebox(0,0){\strut{}\footnotesize{$0$}}}%
      \put(2064,352){\makebox(0,0){\strut{}\footnotesize{$200$}}}%
      \put(3072,352){\makebox(0,0){\strut{}\footnotesize{$400$}}}%
      \put(4081,352){\makebox(0,0){\strut{}\footnotesize{$600$}}}%
      \put(5089,352){\makebox(0,0){\strut{}\footnotesize{$800$}}}%
      \put(6097,352){\makebox(0,0){\strut{}\footnotesize{$1000$}}}%
      \put(7105,352){\makebox(0,0){\strut{}\footnotesize{$1200$}}}%
    }%
    \gplgaddtomacro\gplfronttext{%
      \csname LTb\endcsname
      \put(1258,2104){\makebox(0,0)[l]{\strut{}\footnotesize{$\eps_P \approx 2.220\cdot 10^{-16}$}}}%
      \csname LTb\endcsname
      \put(1258,1950){\makebox(0,0)[l]{\strut{}\footnotesize{$\eps_P \approx 2.242\cdot 10^{-44}$}}}%
      \csname LTb\endcsname
      \put(1258,1796){\makebox(0,0)[l]{\strut{}\footnotesize{$\eps_P \approx 2.264\cdot 10^{-72}$}}}%
      \csname LTb\endcsname
      \put(1258,1642){\makebox(0,0)[l]{\strut{}\footnotesize{$\eps_P \approx 4.572\cdot 10^{-100}$}}}%
      \csname LTb\endcsname
      \put(1258,1488){\makebox(0,0)[l]{\strut{}\footnotesize{$\eps_P \approx 4.616\cdot 10^{-128}$}}}%
      \csname LTb\endcsname
      \put(1258,1334){\makebox(0,0)[l]{\strut{}\footnotesize{$\eps_P \approx 4.661\cdot 10^{-156}$}}}%
      \csname LTb\endcsname
      \put(1258,1180){\makebox(0,0)[l]{\strut{}\footnotesize{$\eps_P \approx 4.707\cdot 10^{-184}$}}}%
      \csname LTb\endcsname
      \put(1258,1026){\makebox(0,0)[l]{\strut{}\footnotesize{$\eps_P \approx 9.505\cdot 10^{-212}$}}}%
      \csname LTb\endcsname
      \put(1258,872){\makebox(0,0)[l]{\strut{}\footnotesize{$\eps_P \approx 9.598\cdot 10^{-240}$}}}%
      \csname LTb\endcsname
      \put(1258,718){\makebox(0,0)[l]{\strut{}\footnotesize{$\eps_P = \eps_{\text{min}}$}}}%
      \csname LTb\endcsname
      \put(187,2115){\rotatebox{-270.00}{\makebox(0,0){\strut{}\footnotesize{$d_H(\,\cdot\,,\,\cdot\,)$}}}}%
      \put(4080,154){\makebox(0,0){\strut{}\footnotesize{N}}}%
    }%
    \gplbacktext
    \put(0,0){\includegraphics[width={368.50bp},height={198.40bp}]{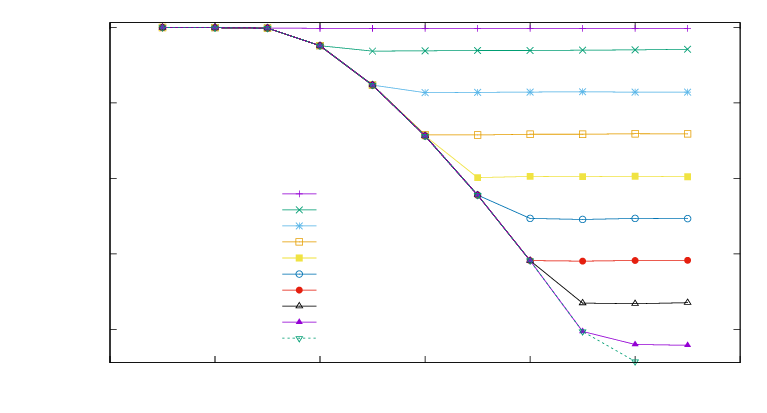}}%
    \gplfronttext
  \end{picture}%
\endgroup

%% file: apcfe-paper-figure-dists_NPmax_Re2e5_N_dH.txt
\begingroup
  \makeatletter
  \providecommand\color[2][]{%
    \GenericError{(gnuplot) \space\space\space\@spaces}{%
      Package color not loaded in conjunction with
      terminal option `colourtext'%
    }{See the gnuplot documentation for explanation.%
    }{Either use 'blacktext' in gnuplot or load the package
      color.sty in LaTeX.}%
    \renewcommand\color[2][]{}%
  }%
  \providecommand\includegraphics[2][]{%
    \GenericError{(gnuplot) \space\space\space\@spaces}{%
      Package graphicx or graphics not loaded%
    }{See the gnuplot documentation for explanation.%
    }{The gnuplot epslatex terminal needs graphicx.sty or graphics.sty.}%
    \renewcommand\includegraphics[2][]{}%
  }%
  \providecommand\rotatebox[2]{#2}%
  \@ifundefined{ifGPcolor}{%
    \newif\ifGPcolor
    \GPcolorfalse
  }{}%
  \@ifundefined{ifGPblacktext}{%
    \newif\ifGPblacktext
    \GPblacktexttrue
  }{}%
  \let\gplgaddtomacro\g@addto@macro
  \gdef\gplbacktext{}%
  \gdef\gplfronttext{}%
  \makeatother
  \ifGPblacktext
    \def\colorrgb#1{}%
    \def\colorgray#1{}%
  \else
    \ifGPcolor
      \def\colorrgb#1{\color[rgb]{#1}}%
      \def\colorgray#1{\color[gray]{#1}}%
      \expandafter\def\csname LTw\endcsname{\color{white}}%
      \expandafter\def\csname LTb\endcsname{\color{black}}%
      \expandafter\def\csname LTa\endcsname{\color{black}}%
      \expandafter\def\csname LT0\endcsname{\color[rgb]{1,0,0}}%
      \expandafter\def\csname LT1\endcsname{\color[rgb]{0,1,0}}%
      \expandafter\def\csname LT2\endcsname{\color[rgb]{0,0,1}}%
      \expandafter\def\csname LT3\endcsname{\color[rgb]{1,0,1}}%
      \expandafter\def\csname LT4\endcsname{\color[rgb]{0,1,1}}%
      \expandafter\def\csname LT5\endcsname{\color[rgb]{1,1,0}}%
      \expandafter\def\csname LT6\endcsname{\color[rgb]{0,0,0}}%
      \expandafter\def\csname LT7\endcsname{\color[rgb]{1,0.3,0}}%
      \expandafter\def\csname LT8\endcsname{\color[rgb]{0.5,0.5,0.5}}%
    \else
      \def\colorrgb#1{\color{black}}%
      \def\colorgray#1{\color[gray]{#1}}%
      \expandafter\def\csname LTw\endcsname{\color{white}}%
      \expandafter\def\csname LTb\endcsname{\color{black}}%
      \expandafter\def\csname LTa\endcsname{\color{black}}%
      \expandafter\def\csname LT0\endcsname{\color{black}}%
      \expandafter\def\csname LT1\endcsname{\color{black}}%
      \expandafter\def\csname LT2\endcsname{\color{black}}%
      \expandafter\def\csname LT3\endcsname{\color{black}}%
      \expandafter\def\csname LT4\endcsname{\color{black}}%
      \expandafter\def\csname LT5\endcsname{\color{black}}%
      \expandafter\def\csname LT6\endcsname{\color{black}}%
      \expandafter\def\csname LT7\endcsname{\color{black}}%
      \expandafter\def\csname LT8\endcsname{\color{black}}%
    \fi
  \fi
    \setlength{\unitlength}{0.0500bp}%
    \ifx\gptboxheight\undefined%
      \newlength{\gptboxheight}%
      \newlength{\gptboxwidth}%
      \newsavebox{\gptboxtext}%
    \fi%
    \setlength{\fboxrule}{0.5pt}%
    \setlength{\fboxsep}{1pt}%
    \definecolor{tbcol}{rgb}{1,1,1}%
\begin{picture}(7370.00,3968.00)%
    \gplgaddtomacro\gplbacktext{%
      \csname LTb\endcsname
      \put(990,3699){\makebox(0,0)[r]{\strut{}\footnotesize{$1$}}}%
      \put(990,516){\makebox(0,0)[r]{\strut{}\footnotesize{$10^{-200}$}}}%
      \put(990,1312){\makebox(0,0)[r]{\strut{}\footnotesize{$10^{-150}$}}}%
      \put(990,2108){\makebox(0,0)[r]{\strut{}\footnotesize{$10^{-100}$}}}%
      \put(990,2903){\makebox(0,0)[r]{\strut{}\footnotesize{$10^{-50}$}}}%
      \put(1056,352){\makebox(0,0){\strut{}\footnotesize{$0$}}}%
      \put(1987,352){\makebox(0,0){\strut{}\footnotesize{$200$}}}%
      \put(2917,352){\makebox(0,0){\strut{}\footnotesize{$400$}}}%
      \put(3848,352){\makebox(0,0){\strut{}\footnotesize{$600$}}}%
      \put(4778,352){\makebox(0,0){\strut{}\footnotesize{$800$}}}%
      \put(5709,352){\makebox(0,0){\strut{}\footnotesize{$1000$}}}%
      \put(6640,352){\makebox(0,0){\strut{}\footnotesize{$1200$}}}%
    }%
    \gplgaddtomacro\gplfronttext{%
      \csname LTb\endcsname
      \put(1242,2190){\makebox(0,0)[l]{\strut{}\footnotesize{$\eps_P \approx 2.220\cdot 10^{-16}$}}}%
      \csname LTb\endcsname
      \put(1242,2036){\makebox(0,0)[l]{\strut{}\footnotesize{$\eps_P \approx 2.242\cdot 10^{-44}$}}}%
      \csname LTb\endcsname
      \put(1242,1882){\makebox(0,0)[l]{\strut{}\footnotesize{$\eps_P \approx 2.264\cdot 10^{-72}$}}}%
      \csname LTb\endcsname
      \put(1242,1728){\makebox(0,0)[l]{\strut{}\footnotesize{$\eps_P \approx 4.572\cdot 10^{-100}$}}}%
      \csname LTb\endcsname
      \put(1242,1574){\makebox(0,0)[l]{\strut{}\footnotesize{$\eps_P \approx 4.616\cdot 10^{-128}$}}}%
      \csname LTb\endcsname
      \put(1242,1420){\makebox(0,0)[l]{\strut{}\footnotesize{$\eps_P \approx 4.661\cdot 10^{-156}$}}}%
      \csname LTb\endcsname
      \put(1242,1266){\makebox(0,0)[l]{\strut{}\footnotesize{$\eps_P \approx 4.707\cdot 10^{-184}$}}}%
      \csname LTb\endcsname
      \put(1242,1112){\makebox(0,0)[l]{\strut{}\footnotesize{$\eps_P \approx 9.505\cdot 10^{-212}$}}}%
      \csname LTb\endcsname
      \put(1242,958){\makebox(0,0)[l]{\strut{}\footnotesize{$\eps_P \approx 9.598\cdot 10^{-240}$}}}%
      \csname LTb\endcsname
      \put(1242,804){\makebox(0,0)[l]{\strut{}\footnotesize{$\eps_P = \eps_{\text{min}}$}}}%
      \csname LTb\endcsname
      \put(187,2115){\rotatebox{-270.00}{\makebox(0,0){\strut{}\footnotesize{$d_H(\,\cdot\,,\,\cdot\,)$}}}}%
      \put(4080,154){\makebox(0,0){\strut{}\footnotesize{N}}}%
    }%
    \gplbacktext
    \put(0,0){\includegraphics[width={368.50bp},height={198.40bp}]{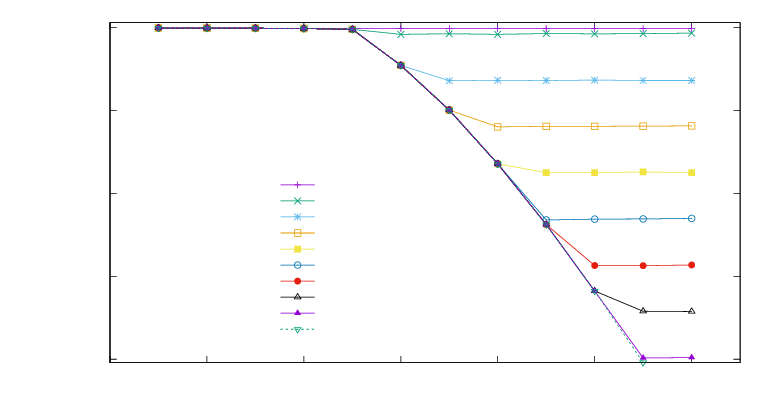}}%
    \gplfronttext
  \end{picture}%
\endgroup

%% file: apcfe-paper-figure-dists_NPmax_Re5e5_N_dH.txt
\begingroup
  \makeatletter
  \providecommand\color[2][]{%
    \GenericError{(gnuplot) \space\space\space\@spaces}{%
      Package color not loaded in conjunction with
      terminal option `colourtext'%
    }{See the gnuplot documentation for explanation.%
    }{Either use 'blacktext' in gnuplot or load the package
      color.sty in LaTeX.}%
    \renewcommand\color[2][]{}%
  }%
  \providecommand\includegraphics[2][]{%
    \GenericError{(gnuplot) \space\space\space\@spaces}{%
      Package graphicx or graphics not loaded%
    }{See the gnuplot documentation for explanation.%
    }{The gnuplot epslatex terminal needs graphicx.sty or graphics.sty.}%
    \renewcommand\includegraphics[2][]{}%
  }%
  \providecommand\rotatebox[2]{#2}%
  \@ifundefined{ifGPcolor}{%
    \newif\ifGPcolor
    \GPcolorfalse
  }{}%
  \@ifundefined{ifGPblacktext}{%
    \newif\ifGPblacktext
    \GPblacktexttrue
  }{}%
  \let\gplgaddtomacro\g@addto@macro
  \gdef\gplbacktext{}%
  \gdef\gplfronttext{}%
  \makeatother
  \ifGPblacktext
    \def\colorrgb#1{}%
    \def\colorgray#1{}%
  \else
    \ifGPcolor
      \def\colorrgb#1{\color[rgb]{#1}}%
      \def\colorgray#1{\color[gray]{#1}}%
      \expandafter\def\csname LTw\endcsname{\color{white}}%
      \expandafter\def\csname LTb\endcsname{\color{black}}%
      \expandafter\def\csname LTa\endcsname{\color{black}}%
      \expandafter\def\csname LT0\endcsname{\color[rgb]{1,0,0}}%
      \expandafter\def\csname LT1\endcsname{\color[rgb]{0,1,0}}%
      \expandafter\def\csname LT2\endcsname{\color[rgb]{0,0,1}}%
      \expandafter\def\csname LT3\endcsname{\color[rgb]{1,0,1}}%
      \expandafter\def\csname LT4\endcsname{\color[rgb]{0,1,1}}%
      \expandafter\def\csname LT5\endcsname{\color[rgb]{1,1,0}}%
      \expandafter\def\csname LT6\endcsname{\color[rgb]{0,0,0}}%
      \expandafter\def\csname LT7\endcsname{\color[rgb]{1,0.3,0}}%
      \expandafter\def\csname LT8\endcsname{\color[rgb]{0.5,0.5,0.5}}%
    \else
      \def\colorrgb#1{\color{black}}%
      \def\colorgray#1{\color[gray]{#1}}%
      \expandafter\def\csname LTw\endcsname{\color{white}}%
      \expandafter\def\csname LTb\endcsname{\color{black}}%
      \expandafter\def\csname LTa\endcsname{\color{black}}%
      \expandafter\def\csname LT0\endcsname{\color{black}}%
      \expandafter\def\csname LT1\endcsname{\color{black}}%
      \expandafter\def\csname LT2\endcsname{\color{black}}%
      \expandafter\def\csname LT3\endcsname{\color{black}}%
      \expandafter\def\csname LT4\endcsname{\color{black}}%
      \expandafter\def\csname LT5\endcsname{\color{black}}%
      \expandafter\def\csname LT6\endcsname{\color{black}}%
      \expandafter\def\csname LT7\endcsname{\color{black}}%
      \expandafter\def\csname LT8\endcsname{\color{black}}%
    \fi
  \fi
    \setlength{\unitlength}{0.0500bp}%
    \ifx\gptboxheight\undefined%
      \newlength{\gptboxheight}%
      \newlength{\gptboxwidth}%
      \newsavebox{\gptboxtext}%
    \fi%
    \setlength{\fboxrule}{0.5pt}%
    \setlength{\fboxsep}{1pt}%
    \definecolor{tbcol}{rgb}{1,1,1}%
\begin{picture}(7370.00,3968.00)%
    \gplgaddtomacro\gplbacktext{%
      \csname LTb\endcsname
      \put(990,3699){\makebox(0,0)[r]{\strut{}\footnotesize{$1$}}}%
      \put(990,484){\makebox(0,0)[r]{\strut{}\footnotesize{$10^{-200}$}}}%
      \put(990,1288){\makebox(0,0)[r]{\strut{}\footnotesize{$10^{-150}$}}}%
      \put(990,2091){\makebox(0,0)[r]{\strut{}\footnotesize{$10^{-100}$}}}%
      \put(990,2895){\makebox(0,0)[r]{\strut{}\footnotesize{$10^{-50}$}}}%
      \put(1056,352){\makebox(0,0){\strut{}\footnotesize{$0$}}}%
      \put(1768,352){\makebox(0,0){\strut{}\footnotesize{$200$}}}%
      \put(2479,352){\makebox(0,0){\strut{}\footnotesize{$400$}}}%
      \put(3191,352){\makebox(0,0){\strut{}\footnotesize{$600$}}}%
      \put(3903,352){\makebox(0,0){\strut{}\footnotesize{$800$}}}%
      \put(4614,352){\makebox(0,0){\strut{}\footnotesize{$1000$}}}%
      \put(5326,352){\makebox(0,0){\strut{}\footnotesize{$1200$}}}%
      \put(6038,352){\makebox(0,0){\strut{}\footnotesize{$1400$}}}%
      \put(6749,352){\makebox(0,0){\strut{}\footnotesize{$1600$}}}%
    }%
    \gplgaddtomacro\gplfronttext{%
      \csname LTb\endcsname
      \put(1198,2175){\makebox(0,0)[l]{\strut{}\footnotesize{$\eps_P \approx 2.220\cdot 10^{-16}$}}}%
      \csname LTb\endcsname
      \put(1198,2021){\makebox(0,0)[l]{\strut{}\footnotesize{$\eps_P \approx 2.242\cdot 10^{-44}$}}}%
      \csname LTb\endcsname
      \put(1198,1867){\makebox(0,0)[l]{\strut{}\footnotesize{$\eps_P \approx 2.264\cdot 10^{-72}$}}}%
      \csname LTb\endcsname
      \put(1198,1713){\makebox(0,0)[l]{\strut{}\footnotesize{$\eps_P \approx 4.572\cdot 10^{-100}$}}}%
      \csname LTb\endcsname
      \put(1198,1559){\makebox(0,0)[l]{\strut{}\footnotesize{$\eps_P \approx 4.616\cdot 10^{-128}$}}}%
      \csname LTb\endcsname
      \put(1198,1405){\makebox(0,0)[l]{\strut{}\footnotesize{$\eps_P \approx 4.661\cdot 10^{-156}$}}}%
      \csname LTb\endcsname
      \put(1198,1251){\makebox(0,0)[l]{\strut{}\footnotesize{$\eps_P \approx 4.707\cdot 10^{-184}$}}}%
      \csname LTb\endcsname
      \put(1198,1097){\makebox(0,0)[l]{\strut{}\footnotesize{$\eps_P \approx 9.505\cdot 10^{-212}$}}}%
      \csname LTb\endcsname
      \put(1198,943){\makebox(0,0)[l]{\strut{}\footnotesize{$\eps_P \approx 9.598\cdot 10^{-240}$}}}%
      \csname LTb\endcsname
      \put(1198,789){\makebox(0,0)[l]{\strut{}\footnotesize{$\eps_P = \eps_{\text{min}}$}}}%
      \csname LTb\endcsname
      \put(187,2115){\rotatebox{-270.00}{\makebox(0,0){\strut{}\footnotesize{$d_H(\,\cdot\,,\,\cdot\,)$}}}}%
      \put(4080,154){\makebox(0,0){\strut{}\footnotesize{N}}}%
    }%
    \gplbacktext
    \put(0,0){\includegraphics[width={368.50bp},height={198.40bp}]{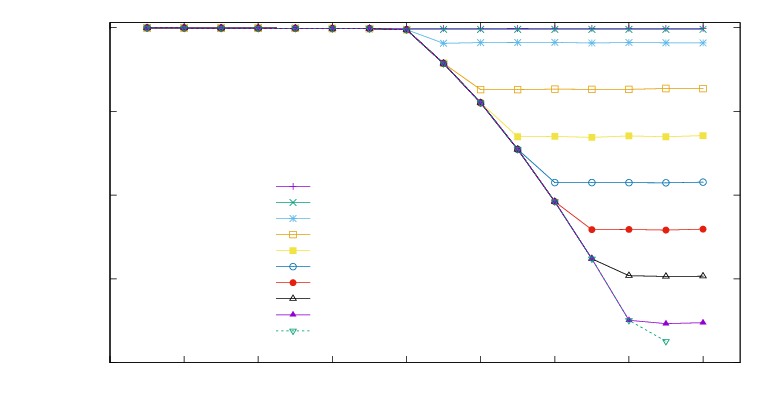}}%
    \gplfronttext
  \end{picture}%
\endgroup

%% file: apcfe-paper-figure-dists_NPmax_Re1e5_P_dH.txt
\begingroup
  \makeatletter
  \providecommand\color[2][]{%
    \GenericError{(gnuplot) \space\space\space\@spaces}{%
      Package color not loaded in conjunction with
      terminal option `colourtext'%
    }{See the gnuplot documentation for explanation.%
    }{Either use 'blacktext' in gnuplot or load the package
      color.sty in LaTeX.}%
    \renewcommand\color[2][]{}%
  }%
  \providecommand\includegraphics[2][]{%
    \GenericError{(gnuplot) \space\space\space\@spaces}{%
      Package graphicx or graphics not loaded%
    }{See the gnuplot documentation for explanation.%
    }{The gnuplot epslatex terminal needs graphicx.sty or graphics.sty.}%
    \renewcommand\includegraphics[2][]{}%
  }%
  \providecommand\rotatebox[2]{#2}%
  \@ifundefined{ifGPcolor}{%
    \newif\ifGPcolor
    \GPcolorfalse
  }{}%
  \@ifundefined{ifGPblacktext}{%
    \newif\ifGPblacktext
    \GPblacktexttrue
  }{}%
  \let\gplgaddtomacro\g@addto@macro
  \gdef\gplbacktext{}%
  \gdef\gplfronttext{}%
  \makeatother
  \ifGPblacktext
    \def\colorrgb#1{}%
    \def\colorgray#1{}%
  \else
    \ifGPcolor
      \def\colorrgb#1{\color[rgb]{#1}}%
      \def\colorgray#1{\color[gray]{#1}}%
      \expandafter\def\csname LTw\endcsname{\color{white}}%
      \expandafter\def\csname LTb\endcsname{\color{black}}%
      \expandafter\def\csname LTa\endcsname{\color{black}}%
      \expandafter\def\csname LT0\endcsname{\color[rgb]{1,0,0}}%
      \expandafter\def\csname LT1\endcsname{\color[rgb]{0,1,0}}%
      \expandafter\def\csname LT2\endcsname{\color[rgb]{0,0,1}}%
      \expandafter\def\csname LT3\endcsname{\color[rgb]{1,0,1}}%
      \expandafter\def\csname LT4\endcsname{\color[rgb]{0,1,1}}%
      \expandafter\def\csname LT5\endcsname{\color[rgb]{1,1,0}}%
      \expandafter\def\csname LT6\endcsname{\color[rgb]{0,0,0}}%
      \expandafter\def\csname LT7\endcsname{\color[rgb]{1,0.3,0}}%
      \expandafter\def\csname LT8\endcsname{\color[rgb]{0.5,0.5,0.5}}%
    \else
      \def\colorrgb#1{\color{black}}%
      \def\colorgray#1{\color[gray]{#1}}%
      \expandafter\def\csname LTw\endcsname{\color{white}}%
      \expandafter\def\csname LTb\endcsname{\color{black}}%
      \expandafter\def\csname LTa\endcsname{\color{black}}%
      \expandafter\def\csname LT0\endcsname{\color{black}}%
      \expandafter\def\csname LT1\endcsname{\color{black}}%
      \expandafter\def\csname LT2\endcsname{\color{black}}%
      \expandafter\def\csname LT3\endcsname{\color{black}}%
      \expandafter\def\csname LT4\endcsname{\color{black}}%
      \expandafter\def\csname LT5\endcsname{\color{black}}%
      \expandafter\def\csname LT6\endcsname{\color{black}}%
      \expandafter\def\csname LT7\endcsname{\color{black}}%
      \expandafter\def\csname LT8\endcsname{\color{black}}%
    \fi
  \fi
    \setlength{\unitlength}{0.0500bp}%
    \ifx\gptboxheight\undefined%
      \newlength{\gptboxheight}%
      \newlength{\gptboxwidth}%
      \newsavebox{\gptboxtext}%
    \fi%
    \setlength{\fboxrule}{0.5pt}%
    \setlength{\fboxsep}{1pt}%
    \definecolor{tbcol}{rgb}{1,1,1}%
\begin{picture}(7370.00,3968.00)%
    \gplgaddtomacro\gplbacktext{%
      \csname LTb\endcsname
      \put(990,3708){\makebox(0,0)[r]{\strut{}\footnotesize{$1$}}}%
      \put(990,484){\makebox(0,0)[r]{\strut{}\footnotesize{$10^{-250}$}}}%
      \put(990,1129){\makebox(0,0)[r]{\strut{}\footnotesize{$10^{-200}$}}}%
      \put(990,1774){\makebox(0,0)[r]{\strut{}\footnotesize{$10^{-150}$}}}%
      \put(990,2419){\makebox(0,0)[r]{\strut{}\footnotesize{$10^{-100}$}}}%
      \put(990,3063){\makebox(0,0)[r]{\strut{}\footnotesize{$10^{-50}$}}}%
      \put(7171,264){\makebox(0,0){\strut{}\footnotesize{$1$}}}%
      \put(1355,264){\makebox(0,0){\strut{}\footnotesize{$10^{-250}$}}}%
      \put(2518,264){\makebox(0,0){\strut{}\footnotesize{$10^{-200}$}}}%
      \put(3681,264){\makebox(0,0){\strut{}\footnotesize{$10^{-150}$}}}%
      \put(4844,264){\makebox(0,0){\strut{}\footnotesize{$10^{-100}$}}}%
      \put(6008,264){\makebox(0,0){\strut{}\footnotesize{$10^{-50}$}}}%
    }%
    \gplgaddtomacro\gplfronttext{%
      \csname LTb\endcsname
      \put(6174,690){\makebox(0,0)[l]{\strut{}\footnotesize{$\mathcal O(\eps_P)$}}}%
      \csname LTb\endcsname
      \put(6174,844){\makebox(0,0)[l]{\strut{}\footnotesize{$N = 1097$}}}%
      \csname LTb\endcsname
      \put(6174,998){\makebox(0,0)[l]{\strut{}\footnotesize{$N = 997$}}}%
      \csname LTb\endcsname
      \put(6174,1152){\makebox(0,0)[l]{\strut{}\footnotesize{$N = 897$}}}%
      \csname LTb\endcsname
      \put(6174,1306){\makebox(0,0)[l]{\strut{}\footnotesize{$N = 797$}}}%
      \csname LTb\endcsname
      \put(6174,1460){\makebox(0,0)[l]{\strut{}\footnotesize{$N = 697$}}}%
      \csname LTb\endcsname
      \put(6174,1614){\makebox(0,0)[l]{\strut{}\footnotesize{$N = 597$}}}%
      \csname LTb\endcsname
      \put(6174,1768){\makebox(0,0)[l]{\strut{}\footnotesize{$N = 497$}}}%
      \csname LTb\endcsname
      \put(6174,1922){\makebox(0,0)[l]{\strut{}\footnotesize{$N = 397$}}}%
      \csname LTb\endcsname
      \put(6174,2076){\makebox(0,0)[l]{\strut{}\footnotesize{$N = 297$}}}%
      \csname LTb\endcsname
      \put(187,2115){\rotatebox{-270.00}{\makebox(0,0){\strut{}\footnotesize{$d_H(\,\cdot\,,\,\cdot\,)$}}}}%
      \put(4080,154){\makebox(0,0){\strut{}\footnotesize{$\eps_P$}}}%
    }%
    \gplbacktext
    \put(0,0){\includegraphics[width={368.50bp},height={198.40bp}]{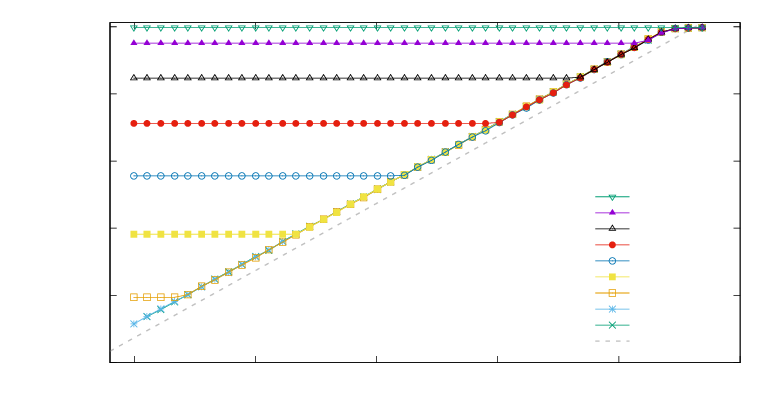}}%
    \gplfronttext
  \end{picture}%
\endgroup

%% file: apcfe-paper-figure-dists_NPmax_Re2e5_P_dH.txt
\begingroup
  \makeatletter
  \providecommand\color[2][]{%
    \GenericError{(gnuplot) \space\space\space\@spaces}{%
      Package color not loaded in conjunction with
      terminal option `colourtext'%
    }{See the gnuplot documentation for explanation.%
    }{Either use 'blacktext' in gnuplot or load the package
      color.sty in LaTeX.}%
    \renewcommand\color[2][]{}%
  }%
  \providecommand\includegraphics[2][]{%
    \GenericError{(gnuplot) \space\space\space\@spaces}{%
      Package graphicx or graphics not loaded%
    }{See the gnuplot documentation for explanation.%
    }{The gnuplot epslatex terminal needs graphicx.sty or graphics.sty.}%
    \renewcommand\includegraphics[2][]{}%
  }%
  \providecommand\rotatebox[2]{#2}%
  \@ifundefined{ifGPcolor}{%
    \newif\ifGPcolor
    \GPcolorfalse
  }{}%
  \@ifundefined{ifGPblacktext}{%
    \newif\ifGPblacktext
    \GPblacktexttrue
  }{}%
  \let\gplgaddtomacro\g@addto@macro
  \gdef\gplbacktext{}%
  \gdef\gplfronttext{}%
  \makeatother
  \ifGPblacktext
    \def\colorrgb#1{}%
    \def\colorgray#1{}%
  \else
    \ifGPcolor
      \def\colorrgb#1{\color[rgb]{#1}}%
      \def\colorgray#1{\color[gray]{#1}}%
      \expandafter\def\csname LTw\endcsname{\color{white}}%
      \expandafter\def\csname LTb\endcsname{\color{black}}%
      \expandafter\def\csname LTa\endcsname{\color{black}}%
      \expandafter\def\csname LT0\endcsname{\color[rgb]{1,0,0}}%
      \expandafter\def\csname LT1\endcsname{\color[rgb]{0,1,0}}%
      \expandafter\def\csname LT2\endcsname{\color[rgb]{0,0,1}}%
      \expandafter\def\csname LT3\endcsname{\color[rgb]{1,0,1}}%
      \expandafter\def\csname LT4\endcsname{\color[rgb]{0,1,1}}%
      \expandafter\def\csname LT5\endcsname{\color[rgb]{1,1,0}}%
      \expandafter\def\csname LT6\endcsname{\color[rgb]{0,0,0}}%
      \expandafter\def\csname LT7\endcsname{\color[rgb]{1,0.3,0}}%
      \expandafter\def\csname LT8\endcsname{\color[rgb]{0.5,0.5,0.5}}%
    \else
      \def\colorrgb#1{\color{black}}%
      \def\colorgray#1{\color[gray]{#1}}%
      \expandafter\def\csname LTw\endcsname{\color{white}}%
      \expandafter\def\csname LTb\endcsname{\color{black}}%
      \expandafter\def\csname LTa\endcsname{\color{black}}%
      \expandafter\def\csname LT0\endcsname{\color{black}}%
      \expandafter\def\csname LT1\endcsname{\color{black}}%
      \expandafter\def\csname LT2\endcsname{\color{black}}%
      \expandafter\def\csname LT3\endcsname{\color{black}}%
      \expandafter\def\csname LT4\endcsname{\color{black}}%
      \expandafter\def\csname LT5\endcsname{\color{black}}%
      \expandafter\def\csname LT6\endcsname{\color{black}}%
      \expandafter\def\csname LT7\endcsname{\color{black}}%
      \expandafter\def\csname LT8\endcsname{\color{black}}%
    \fi
  \fi
    \setlength{\unitlength}{0.0500bp}%
    \ifx\gptboxheight\undefined%
      \newlength{\gptboxheight}%
      \newlength{\gptboxwidth}%
      \newsavebox{\gptboxtext}%
    \fi%
    \setlength{\fboxrule}{0.5pt}%
    \setlength{\fboxsep}{1pt}%
    \definecolor{tbcol}{rgb}{1,1,1}%
\begin{picture}(7370.00,3968.00)%
    \gplgaddtomacro\gplbacktext{%
      \csname LTb\endcsname
      \put(990,3704){\makebox(0,0)[r]{\strut{}\footnotesize{$1$}}}%
      \put(990,842){\makebox(0,0)[r]{\strut{}\footnotesize{$10^{-200}$}}}%
      \put(990,1557){\makebox(0,0)[r]{\strut{}\footnotesize{$10^{-150}$}}}%
      \put(990,2273){\makebox(0,0)[r]{\strut{}\footnotesize{$10^{-100}$}}}%
      \put(990,2988){\makebox(0,0)[r]{\strut{}\footnotesize{$10^{-50}$}}}%
      \put(7171,264){\makebox(0,0){\strut{}\footnotesize{$1$}}}%
      \put(1355,264){\makebox(0,0){\strut{}\footnotesize{$10^{-250}$}}}%
      \put(2518,264){\makebox(0,0){\strut{}\footnotesize{$10^{-200}$}}}%
      \put(3681,264){\makebox(0,0){\strut{}\footnotesize{$10^{-150}$}}}%
      \put(4844,264){\makebox(0,0){\strut{}\footnotesize{$10^{-100}$}}}%
      \put(6008,264){\makebox(0,0){\strut{}\footnotesize{$10^{-50}$}}}%
    }%
    \gplgaddtomacro\gplfronttext{%
      \csname LTb\endcsname
      \put(6174,704){\makebox(0,0)[l]{\strut{}\footnotesize{$\mathcal O(\eps_P)$}}}%
      \csname LTb\endcsname
      \put(6174,858){\makebox(0,0)[l]{\strut{}\footnotesize{$N = 1197$}}}%
      \csname LTb\endcsname
      \put(6174,1012){\makebox(0,0)[l]{\strut{}\footnotesize{$N = 1097$}}}%
      \csname LTb\endcsname
      \put(6174,1166){\makebox(0,0)[l]{\strut{}\footnotesize{$N = 997$}}}%
      \csname LTb\endcsname
      \put(6174,1320){\makebox(0,0)[l]{\strut{}\footnotesize{$N = 897$}}}%
      \csname LTb\endcsname
      \put(6174,1474){\makebox(0,0)[l]{\strut{}\footnotesize{$N = 797$}}}%
      \csname LTb\endcsname
      \put(6174,1628){\makebox(0,0)[l]{\strut{}\footnotesize{$N = 697$}}}%
      \csname LTb\endcsname
      \put(6174,1782){\makebox(0,0)[l]{\strut{}\footnotesize{$N = 597$}}}%
      \csname LTb\endcsname
      \put(6174,1936){\makebox(0,0)[l]{\strut{}\footnotesize{$N = 497$}}}%
      \csname LTb\endcsname
      \put(187,2115){\rotatebox{-270.00}{\makebox(0,0){\strut{}\footnotesize{$d_H(\,\cdot\,,\,\cdot\,)$}}}}%
      \put(4080,154){\makebox(0,0){\strut{}\footnotesize{$\eps_P$}}}%
    }%
    \gplbacktext
    \put(0,0){\includegraphics[width={368.50bp},height={198.40bp}]{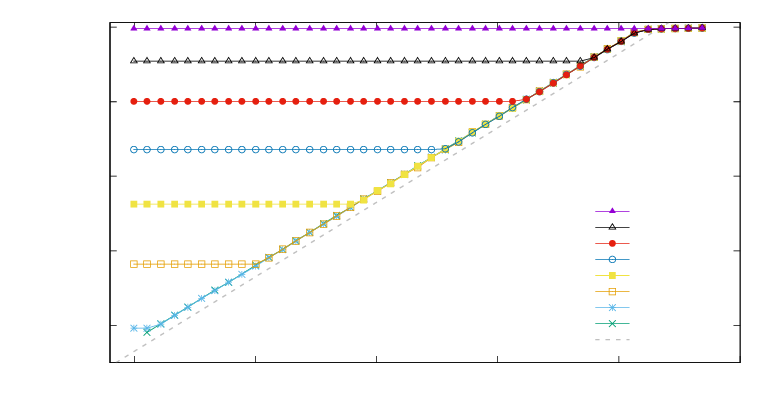}}%
    \gplfronttext
  \end{picture}%
\endgroup

%% file: apcfe-paper-figure-dists_NPmax_Re5e5_P_dH.txt
\begingroup
  \makeatletter
  \providecommand\color[2][]{%
    \GenericError{(gnuplot) \space\space\space\@spaces}{%
      Package color not loaded in conjunction with
      terminal option `colourtext'%
    }{See the gnuplot documentation for explanation.%
    }{Either use 'blacktext' in gnuplot or load the package
      color.sty in LaTeX.}%
    \renewcommand\color[2][]{}%
  }%
  \providecommand\includegraphics[2][]{%
    \GenericError{(gnuplot) \space\space\space\@spaces}{%
      Package graphicx or graphics not loaded%
    }{See the gnuplot documentation for explanation.%
    }{The gnuplot epslatex terminal needs graphicx.sty or graphics.sty.}%
    \renewcommand\includegraphics[2][]{}%
  }%
  \providecommand\rotatebox[2]{#2}%
  \@ifundefined{ifGPcolor}{%
    \newif\ifGPcolor
    \GPcolorfalse
  }{}%
  \@ifundefined{ifGPblacktext}{%
    \newif\ifGPblacktext
    \GPblacktexttrue
  }{}%
  \let\gplgaddtomacro\g@addto@macro
  \gdef\gplbacktext{}%
  \gdef\gplfronttext{}%
  \makeatother
  \ifGPblacktext
    \def\colorrgb#1{}%
    \def\colorgray#1{}%
  \else
    \ifGPcolor
      \def\colorrgb#1{\color[rgb]{#1}}%
      \def\colorgray#1{\color[gray]{#1}}%
      \expandafter\def\csname LTw\endcsname{\color{white}}%
      \expandafter\def\csname LTb\endcsname{\color{black}}%
      \expandafter\def\csname LTa\endcsname{\color{black}}%
      \expandafter\def\csname LT0\endcsname{\color[rgb]{1,0,0}}%
      \expandafter\def\csname LT1\endcsname{\color[rgb]{0,1,0}}%
      \expandafter\def\csname LT2\endcsname{\color[rgb]{0,0,1}}%
      \expandafter\def\csname LT3\endcsname{\color[rgb]{1,0,1}}%
      \expandafter\def\csname LT4\endcsname{\color[rgb]{0,1,1}}%
      \expandafter\def\csname LT5\endcsname{\color[rgb]{1,1,0}}%
      \expandafter\def\csname LT6\endcsname{\color[rgb]{0,0,0}}%
      \expandafter\def\csname LT7\endcsname{\color[rgb]{1,0.3,0}}%
      \expandafter\def\csname LT8\endcsname{\color[rgb]{0.5,0.5,0.5}}%
    \else
      \def\colorrgb#1{\color{black}}%
      \def\colorgray#1{\color[gray]{#1}}%
      \expandafter\def\csname LTw\endcsname{\color{white}}%
      \expandafter\def\csname LTb\endcsname{\color{black}}%
      \expandafter\def\csname LTa\endcsname{\color{black}}%
      \expandafter\def\csname LT0\endcsname{\color{black}}%
      \expandafter\def\csname LT1\endcsname{\color{black}}%
      \expandafter\def\csname LT2\endcsname{\color{black}}%
      \expandafter\def\csname LT3\endcsname{\color{black}}%
      \expandafter\def\csname LT4\endcsname{\color{black}}%
      \expandafter\def\csname LT5\endcsname{\color{black}}%
      \expandafter\def\csname LT6\endcsname{\color{black}}%
      \expandafter\def\csname LT7\endcsname{\color{black}}%
      \expandafter\def\csname LT8\endcsname{\color{black}}%
    \fi
  \fi
    \setlength{\unitlength}{0.0500bp}%
    \ifx\gptboxheight\undefined%
      \newlength{\gptboxheight}%
      \newlength{\gptboxwidth}%
      \newsavebox{\gptboxtext}%
    \fi%
    \setlength{\fboxrule}{0.5pt}%
    \setlength{\fboxsep}{1pt}%
    \definecolor{tbcol}{rgb}{1,1,1}%
\begin{picture}(7370.00,3968.00)%
    \gplgaddtomacro\gplbacktext{%
      \csname LTb\endcsname
      \put(990,3699){\makebox(0,0)[r]{\strut{}\footnotesize{$1$}}}%
      \put(990,484){\makebox(0,0)[r]{\strut{}\footnotesize{$10^{-200}$}}}%
      \put(990,1288){\makebox(0,0)[r]{\strut{}\footnotesize{$10^{-150}$}}}%
      \put(990,2091){\makebox(0,0)[r]{\strut{}\footnotesize{$10^{-100}$}}}%
      \put(990,2895){\makebox(0,0)[r]{\strut{}\footnotesize{$10^{-50}$}}}%
      \put(7171,264){\makebox(0,0){\strut{}\footnotesize{$1$}}}%
      \put(1355,264){\makebox(0,0){\strut{}\footnotesize{$10^{-250}$}}}%
      \put(2518,264){\makebox(0,0){\strut{}\footnotesize{$10^{-200}$}}}%
      \put(3681,264){\makebox(0,0){\strut{}\footnotesize{$10^{-150}$}}}%
      \put(4844,264){\makebox(0,0){\strut{}\footnotesize{$10^{-100}$}}}%
      \put(6008,264){\makebox(0,0){\strut{}\footnotesize{$10^{-50}$}}}%
    }%
    \gplgaddtomacro\gplfronttext{%
      \csname LTb\endcsname
      \put(6174,722){\makebox(0,0)[l]{\strut{}\footnotesize{$\mathcal O(\eps_P)$}}}%
      \csname LTb\endcsname
      \put(6174,876){\makebox(0,0)[l]{\strut{}\footnotesize{$N = 1497$}}}%
      \csname LTb\endcsname
      \put(6174,1030){\makebox(0,0)[l]{\strut{}\footnotesize{$N = 1397$}}}%
      \csname LTb\endcsname
      \put(6174,1184){\makebox(0,0)[l]{\strut{}\footnotesize{$N = 1297$}}}%
      \csname LTb\endcsname
      \put(6174,1338){\makebox(0,0)[l]{\strut{}\footnotesize{$N = 1197$}}}%
      \csname LTb\endcsname
      \put(6174,1492){\makebox(0,0)[l]{\strut{}\footnotesize{$N = 1097$}}}%
      \csname LTb\endcsname
      \put(6174,1646){\makebox(0,0)[l]{\strut{}\footnotesize{$N = 997$}}}%
      \csname LTb\endcsname
      \put(6174,1800){\makebox(0,0)[l]{\strut{}\footnotesize{$N = 897$}}}%
      \csname LTb\endcsname
      \put(6174,1954){\makebox(0,0)[l]{\strut{}\footnotesize{$N = 797$}}}%
      \csname LTb\endcsname
      \put(187,2115){\rotatebox{-270.00}{\makebox(0,0){\strut{}\footnotesize{$d_H(\,\cdot\,,\,\cdot\,)$}}}}%
      \put(4080,154){\makebox(0,0){\strut{}\footnotesize{$\eps_{P}$}}}%
    }%
    \gplbacktext
    \put(0,0){\includegraphics[width={368.50bp},height={198.40bp}]{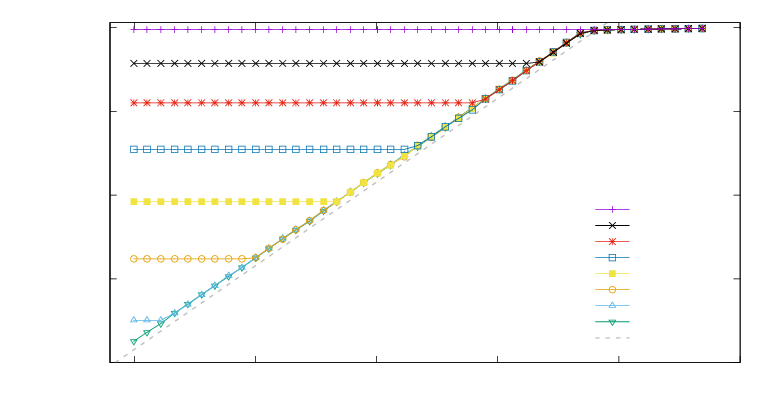}}%
    \gplfronttext
  \end{picture}%
\endgroup

%% file: apcfe-paper-figure-out_pois_2e5_spectra.txt
\begingroup
  \makeatletter
  \providecommand\color[2][]{%
    \GenericError{(gnuplot) \space\space\space\@spaces}{%
      Package color not loaded in conjunction with
      terminal option `colourtext'%
    }{See the gnuplot documentation for explanation.%
    }{Either use 'blacktext' in gnuplot or load the package
      color.sty in LaTeX.}%
    \renewcommand\color[2][]{}%
  }%
  \providecommand\includegraphics[2][]{%
    \GenericError{(gnuplot) \space\space\space\@spaces}{%
      Package graphicx or graphics not loaded%
    }{See the gnuplot documentation for explanation.%
    }{The gnuplot epslatex terminal needs graphicx.sty or graphics.sty.}%
    \renewcommand\includegraphics[2][]{}%
  }%
  \providecommand\rotatebox[2]{#2}%
  \@ifundefined{ifGPcolor}{%
    \newif\ifGPcolor
    \GPcolorfalse
  }{}%
  \@ifundefined{ifGPblacktext}{%
    \newif\ifGPblacktext
    \GPblacktexttrue
  }{}%
  \let\gplgaddtomacro\g@addto@macro
  \gdef\gplbacktext{}%
  \gdef\gplfronttext{}%
  \makeatother
  \ifGPblacktext
    \def\colorrgb#1{}%
    \def\colorgray#1{}%
  \else
    \ifGPcolor
      \def\colorrgb#1{\color[rgb]{#1}}%
      \def\colorgray#1{\color[gray]{#1}}%
      \expandafter\def\csname LTw\endcsname{\color{white}}%
      \expandafter\def\csname LTb\endcsname{\color{black}}%
      \expandafter\def\csname LTa\endcsname{\color{black}}%
      \expandafter\def\csname LT0\endcsname{\color[rgb]{1,0,0}}%
      \expandafter\def\csname LT1\endcsname{\color[rgb]{0,1,0}}%
      \expandafter\def\csname LT2\endcsname{\color[rgb]{0,0,1}}%
      \expandafter\def\csname LT3\endcsname{\color[rgb]{1,0,1}}%
      \expandafter\def\csname LT4\endcsname{\color[rgb]{0,1,1}}%
      \expandafter\def\csname LT5\endcsname{\color[rgb]{1,1,0}}%
      \expandafter\def\csname LT6\endcsname{\color[rgb]{0,0,0}}%
      \expandafter\def\csname LT7\endcsname{\color[rgb]{1,0.3,0}}%
      \expandafter\def\csname LT8\endcsname{\color[rgb]{0.5,0.5,0.5}}%
    \else
      \def\colorrgb#1{\color{black}}%
      \def\colorgray#1{\color[gray]{#1}}%
      \expandafter\def\csname LTw\endcsname{\color{white}}%
      \expandafter\def\csname LTb\endcsname{\color{black}}%
      \expandafter\def\csname LTa\endcsname{\color{black}}%
      \expandafter\def\csname LT0\endcsname{\color{black}}%
      \expandafter\def\csname LT1\endcsname{\color{black}}%
      \expandafter\def\csname LT2\endcsname{\color{black}}%
      \expandafter\def\csname LT3\endcsname{\color{black}}%
      \expandafter\def\csname LT4\endcsname{\color{black}}%
      \expandafter\def\csname LT5\endcsname{\color{black}}%
      \expandafter\def\csname LT6\endcsname{\color{black}}%
      \expandafter\def\csname LT7\endcsname{\color{black}}%
      \expandafter\def\csname LT8\endcsname{\color{black}}%
    \fi
  \fi
    \setlength{\unitlength}{0.0500bp}%
    \ifx\gptboxheight\undefined%
      \newlength{\gptboxheight}%
      \newlength{\gptboxwidth}%
      \newsavebox{\gptboxtext}%
    \fi%
    \setlength{\fboxrule}{0.5pt}%
    \setlength{\fboxsep}{1pt}%
    \definecolor{tbcol}{rgb}{1,1,1}%
\begin{picture}(7936.00,7936.00)%
    \gplgaddtomacro\gplbacktext{%
      \csname LTb\endcsname
      \put(138,5752){\makebox(0,0)[r]{\strut{}\footnotesize{-1}}}%
      \put(138,7737){\makebox(0,0)[r]{\strut{}\footnotesize{0}}}%
      \put(165,5642){\makebox(0,0){\strut{}\footnotesize{0}}}%
      \put(2150,5642){\makebox(0,0){\strut{}\footnotesize{1}}}%
    }%
    \gplgaddtomacro\gplfronttext{%
      \csname LTb\endcsname
      \put(22,6722){\rotatebox{-270}{\makebox(0,0){\strut{}\footnotesize{imag}}}}%
      \put(1157,5554){\makebox(0,0){\strut{}\footnotesize{real}}}%
      \csname LTb\endcsname
      \put(1823,5322){\makebox(0,0)[r]{\strut{}\footnotesize{$P = 53, N = 100$}}}%
    }%
    \gplgaddtomacro\gplbacktext{%
      \csname LTb\endcsname
      \put(2783,5752){\makebox(0,0)[r]{\strut{}\footnotesize{-1}}}%
      \put(2783,7737){\makebox(0,0)[r]{\strut{}\footnotesize{0}}}%
      \put(2810,5642){\makebox(0,0){\strut{}\footnotesize{0}}}%
      \put(4795,5642){\makebox(0,0){\strut{}\footnotesize{1}}}%
    }%
    \gplgaddtomacro\gplfronttext{%
      \csname LTb\endcsname
      \put(2667,6722){\rotatebox{-270}{\makebox(0,0){\strut{}\footnotesize{imag}}}}%
      \put(3802,5554){\makebox(0,0){\strut{}\footnotesize{real}}}%
      \csname LTb\endcsname
      \put(4468,5322){\makebox(0,0)[r]{\strut{}\footnotesize{$P = 53, N = 400$}}}%
    }%
    \gplgaddtomacro\gplbacktext{%
      \csname LTb\endcsname
      \put(5428,5752){\makebox(0,0)[r]{\strut{}\footnotesize{-1}}}%
      \put(5428,7737){\makebox(0,0)[r]{\strut{}\footnotesize{0}}}%
      \put(5455,5642){\makebox(0,0){\strut{}\footnotesize{0}}}%
      \put(7440,5642){\makebox(0,0){\strut{}\footnotesize{1}}}%
    }%
    \gplgaddtomacro\gplfronttext{%
      \csname LTb\endcsname
      \put(5312,6722){\rotatebox{-270}{\makebox(0,0){\strut{}\footnotesize{imag}}}}%
      \put(6447,5554){\makebox(0,0){\strut{}\footnotesize{real}}}%
      \csname LTb\endcsname
      \put(7113,5322){\makebox(0,0)[r]{\strut{}\footnotesize{$P = 53, N = 1200$}}}%
    }%
    \gplgaddtomacro\gplbacktext{%
      \csname LTb\endcsname
      \put(138,3107){\makebox(0,0)[r]{\strut{}\footnotesize{-1}}}%
      \put(138,5092){\makebox(0,0)[r]{\strut{}\footnotesize{0}}}%
      \put(165,2997){\makebox(0,0){\strut{}\footnotesize{0}}}%
      \put(2150,2997){\makebox(0,0){\strut{}\footnotesize{1}}}%
    }%
    \gplgaddtomacro\gplfronttext{%
      \csname LTb\endcsname
      \put(22,4077){\rotatebox{-270}{\makebox(0,0){\strut{}\footnotesize{imag}}}}%
      \put(1157,2909){\makebox(0,0){\strut{}\footnotesize{real}}}%
      \csname LTb\endcsname
      \put(1823,2677){\makebox(0,0)[r]{\strut{}\footnotesize{$P = 424, N = 100$}}}%
    }%
    \gplgaddtomacro\gplbacktext{%
      \csname LTb\endcsname
      \put(2783,3107){\makebox(0,0)[r]{\strut{}\footnotesize{-1}}}%
      \put(2783,5092){\makebox(0,0)[r]{\strut{}\footnotesize{0}}}%
      \put(2810,2997){\makebox(0,0){\strut{}\footnotesize{0}}}%
      \put(4795,2997){\makebox(0,0){\strut{}\footnotesize{1}}}%
    }%
    \gplgaddtomacro\gplfronttext{%
      \csname LTb\endcsname
      \put(2667,4077){\rotatebox{-270}{\makebox(0,0){\strut{}\footnotesize{imag}}}}%
      \put(3802,2909){\makebox(0,0){\strut{}\footnotesize{real}}}%
      \csname LTb\endcsname
      \put(4468,2677){\makebox(0,0)[r]{\strut{}\footnotesize{$P = 424, N = 400$}}}%
    }%
    \gplgaddtomacro\gplbacktext{%
      \csname LTb\endcsname
      \put(5428,3107){\makebox(0,0)[r]{\strut{}\footnotesize{-1}}}%
      \put(5428,5092){\makebox(0,0)[r]{\strut{}\footnotesize{0}}}%
      \put(5455,2997){\makebox(0,0){\strut{}\footnotesize{0}}}%
      \put(7440,2997){\makebox(0,0){\strut{}\footnotesize{1}}}%
    }%
    \gplgaddtomacro\gplfronttext{%
      \csname LTb\endcsname
      \put(5312,4077){\rotatebox{-270}{\makebox(0,0){\strut{}\footnotesize{imag}}}}%
      \put(6447,2909){\makebox(0,0){\strut{}\footnotesize{real}}}%
      \csname LTb\endcsname
      \put(7113,2677){\makebox(0,0)[r]{\strut{}\footnotesize{$P = 424, N = 1200$}}}%
    }%
    \gplgaddtomacro\gplbacktext{%
      \csname LTb\endcsname
      \put(138,462){\makebox(0,0)[r]{\strut{}\footnotesize{-1}}}%
      \put(138,2447){\makebox(0,0)[r]{\strut{}\footnotesize{0}}}%
      \put(165,352){\makebox(0,0){\strut{}\footnotesize{0}}}%
      \put(2150,352){\makebox(0,0){\strut{}\footnotesize{1}}}%
    }%
    \gplgaddtomacro\gplfronttext{%
      \csname LTb\endcsname
      \put(22,1432){\rotatebox{-270}{\makebox(0,0){\strut{}\footnotesize{imag}}}}%
      \put(1157,264){\makebox(0,0){\strut{}\footnotesize{real}}}%
      \csname LTb\endcsname
      \put(1823,32){\makebox(0,0)[r]{\strut{}\footnotesize{$P = 832, N = 100$}}}%
    }%
    \gplgaddtomacro\gplbacktext{%
      \csname LTb\endcsname
      \put(2783,462){\makebox(0,0)[r]{\strut{}\footnotesize{-1}}}%
      \put(2783,2447){\makebox(0,0)[r]{\strut{}\footnotesize{0}}}%
      \put(2810,352){\makebox(0,0){\strut{}\footnotesize{0}}}%
      \put(4795,352){\makebox(0,0){\strut{}\footnotesize{1}}}%
    }%
    \gplgaddtomacro\gplfronttext{%
      \csname LTb\endcsname
      \put(2667,1432){\rotatebox{-270}{\makebox(0,0){\strut{}\footnotesize{imag}}}}%
      \put(3802,264){\makebox(0,0){\strut{}\footnotesize{real}}}%
      \csname LTb\endcsname
      \put(4468,32){\makebox(0,0)[r]{\strut{}\footnotesize{$P = 832, N = 400$}}}%
    }%
    \gplgaddtomacro\gplbacktext{%
      \csname LTb\endcsname
      \put(5428,462){\makebox(0,0)[r]{\strut{}\footnotesize{-1}}}%
      \put(5428,2447){\makebox(0,0)[r]{\strut{}\footnotesize{0}}}%
      \put(5455,352){\makebox(0,0){\strut{}\footnotesize{0}}}%
      \put(7440,352){\makebox(0,0){\strut{}\footnotesize{1}}}%
    }%
    \gplgaddtomacro\gplfronttext{%
      \csname LTb\endcsname
      \put(5312,1432){\rotatebox{-270}{\makebox(0,0){\strut{}\footnotesize{imag}}}}%
      \put(6447,264){\makebox(0,0){\strut{}\footnotesize{real}}}%
      \csname LTb\endcsname
      \put(7113,32){\makebox(0,0)[r]{\strut{}\footnotesize{$P = 832, N = 1200$}}}%
    }%
    \gplbacktext
    \put(0,0){\includegraphics[width={396.80bp},height={396.80bp}]{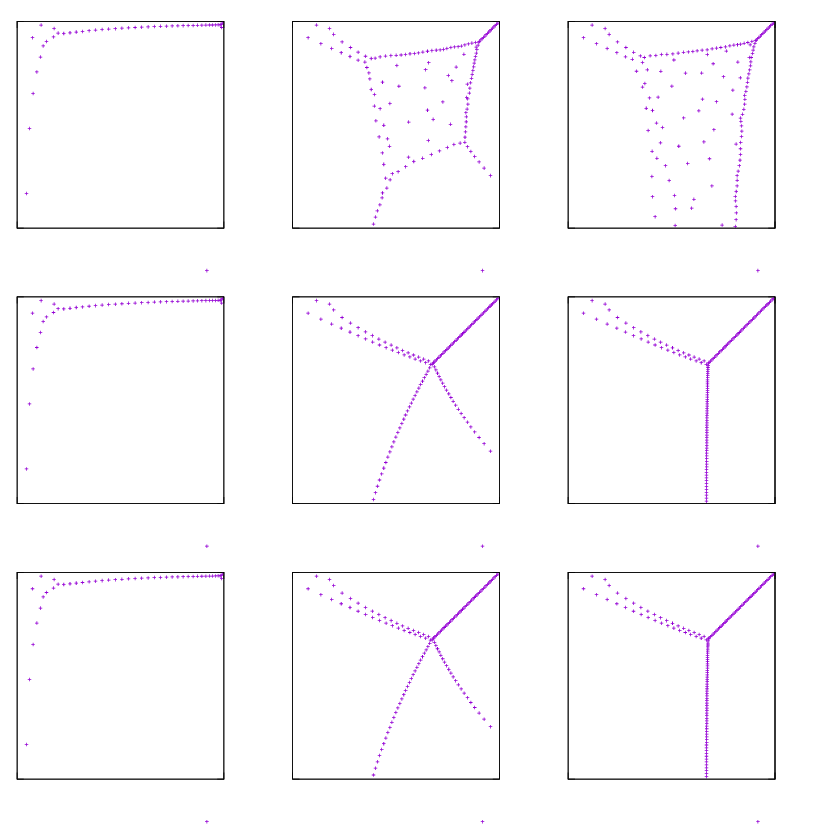}}%
    \gplfronttext
  \end{picture}%
\endgroup

%% file: apcfe-paper-figure-out_pois_5e5_spectra.txt
\begingroup
  \makeatletter
  \providecommand\color[2][]{%
    \GenericError{(gnuplot) \space\space\space\@spaces}{%
      Package color not loaded in conjunction with
      terminal option `colourtext'%
    }{See the gnuplot documentation for explanation.%
    }{Either use 'blacktext' in gnuplot or load the package
      color.sty in LaTeX.}%
    \renewcommand\color[2][]{}%
  }%
  \providecommand\includegraphics[2][]{%
    \GenericError{(gnuplot) \space\space\space\@spaces}{%
      Package graphicx or graphics not loaded%
    }{See the gnuplot documentation for explanation.%
    }{The gnuplot epslatex terminal needs graphicx.sty or graphics.sty.}%
    \renewcommand\includegraphics[2][]{}%
  }%
  \providecommand\rotatebox[2]{#2}%
  \@ifundefined{ifGPcolor}{%
    \newif\ifGPcolor
    \GPcolorfalse
  }{}%
  \@ifundefined{ifGPblacktext}{%
    \newif\ifGPblacktext
    \GPblacktexttrue
  }{}%
  \let\gplgaddtomacro\g@addto@macro
  \gdef\gplbacktext{}%
  \gdef\gplfronttext{}%
  \makeatother
  \ifGPblacktext
    \def\colorrgb#1{}%
    \def\colorgray#1{}%
  \else
    \ifGPcolor
      \def\colorrgb#1{\color[rgb]{#1}}%
      \def\colorgray#1{\color[gray]{#1}}%
      \expandafter\def\csname LTw\endcsname{\color{white}}%
      \expandafter\def\csname LTb\endcsname{\color{black}}%
      \expandafter\def\csname LTa\endcsname{\color{black}}%
      \expandafter\def\csname LT0\endcsname{\color[rgb]{1,0,0}}%
      \expandafter\def\csname LT1\endcsname{\color[rgb]{0,1,0}}%
      \expandafter\def\csname LT2\endcsname{\color[rgb]{0,0,1}}%
      \expandafter\def\csname LT3\endcsname{\color[rgb]{1,0,1}}%
      \expandafter\def\csname LT4\endcsname{\color[rgb]{0,1,1}}%
      \expandafter\def\csname LT5\endcsname{\color[rgb]{1,1,0}}%
      \expandafter\def\csname LT6\endcsname{\color[rgb]{0,0,0}}%
      \expandafter\def\csname LT7\endcsname{\color[rgb]{1,0.3,0}}%
      \expandafter\def\csname LT8\endcsname{\color[rgb]{0.5,0.5,0.5}}%
    \else
      \def\colorrgb#1{\color{black}}%
      \def\colorgray#1{\color[gray]{#1}}%
      \expandafter\def\csname LTw\endcsname{\color{white}}%
      \expandafter\def\csname LTb\endcsname{\color{black}}%
      \expandafter\def\csname LTa\endcsname{\color{black}}%
      \expandafter\def\csname LT0\endcsname{\color{black}}%
      \expandafter\def\csname LT1\endcsname{\color{black}}%
      \expandafter\def\csname LT2\endcsname{\color{black}}%
      \expandafter\def\csname LT3\endcsname{\color{black}}%
      \expandafter\def\csname LT4\endcsname{\color{black}}%
      \expandafter\def\csname LT5\endcsname{\color{black}}%
      \expandafter\def\csname LT6\endcsname{\color{black}}%
      \expandafter\def\csname LT7\endcsname{\color{black}}%
      \expandafter\def\csname LT8\endcsname{\color{black}}%
    \fi
  \fi
    \setlength{\unitlength}{0.0500bp}%
    \ifx\gptboxheight\undefined%
      \newlength{\gptboxheight}%
      \newlength{\gptboxwidth}%
      \newsavebox{\gptboxtext}%
    \fi%
    \setlength{\fboxrule}{0.5pt}%
    \setlength{\fboxsep}{1pt}%
    \definecolor{tbcol}{rgb}{1,1,1}%
\begin{picture}(7936.00,7936.00)%
    \gplgaddtomacro\gplbacktext{%
      \csname LTb\endcsname
      \put(138,5752){\makebox(0,0)[r]{\strut{}\footnotesize{-1}}}%
      \put(138,7737){\makebox(0,0)[r]{\strut{}\footnotesize{0}}}%
      \put(165,5642){\makebox(0,0){\strut{}\footnotesize{0}}}%
      \put(2150,5642){\makebox(0,0){\strut{}\footnotesize{1}}}%
    }%
    \gplgaddtomacro\gplfronttext{%
      \csname LTb\endcsname
      \put(22,6722){\rotatebox{-270}{\makebox(0,0){\strut{}\footnotesize{imag}}}}%
      \put(1157,5554){\makebox(0,0){\strut{}\footnotesize{real}}}%
      \csname LTb\endcsname
      \put(1823,5322){\makebox(0,0)[r]{\strut{}\footnotesize{$P = 53, N = 100$}}}%
    }%
    \gplgaddtomacro\gplbacktext{%
      \csname LTb\endcsname
      \put(2783,5752){\makebox(0,0)[r]{\strut{}\footnotesize{-1}}}%
      \put(2783,7737){\makebox(0,0)[r]{\strut{}\footnotesize{0}}}%
      \put(2810,5642){\makebox(0,0){\strut{}\footnotesize{0}}}%
      \put(4795,5642){\makebox(0,0){\strut{}\footnotesize{1}}}%
    }%
    \gplgaddtomacro\gplfronttext{%
      \csname LTb\endcsname
      \put(2667,6722){\rotatebox{-270}{\makebox(0,0){\strut{}\footnotesize{imag}}}}%
      \put(3802,5554){\makebox(0,0){\strut{}\footnotesize{real}}}%
      \csname LTb\endcsname
      \put(4468,5322){\makebox(0,0)[r]{\strut{}\footnotesize{$P = 53, N = 600$}}}%
    }%
    \gplgaddtomacro\gplbacktext{%
      \csname LTb\endcsname
      \put(5428,5752){\makebox(0,0)[r]{\strut{}\footnotesize{-1}}}%
      \put(5428,7737){\makebox(0,0)[r]{\strut{}\footnotesize{0}}}%
      \put(5455,5642){\makebox(0,0){\strut{}\footnotesize{0}}}%
      \put(7440,5642){\makebox(0,0){\strut{}\footnotesize{1}}}%
    }%
    \gplgaddtomacro\gplfronttext{%
      \csname LTb\endcsname
      \put(5312,6722){\rotatebox{-270}{\makebox(0,0){\strut{}\footnotesize{imag}}}}%
      \put(6447,5554){\makebox(0,0){\strut{}\footnotesize{real}}}%
      \csname LTb\endcsname
      \put(7113,5322){\makebox(0,0)[r]{\strut{}\footnotesize{$P = 53, N = 1500$}}}%
    }%
    \gplgaddtomacro\gplbacktext{%
      \csname LTb\endcsname
      \put(138,3107){\makebox(0,0)[r]{\strut{}\footnotesize{-1}}}%
      \put(138,5092){\makebox(0,0)[r]{\strut{}\footnotesize{0}}}%
      \put(165,2997){\makebox(0,0){\strut{}\footnotesize{0}}}%
      \put(2150,2997){\makebox(0,0){\strut{}\footnotesize{1}}}%
    }%
    \gplgaddtomacro\gplfronttext{%
      \csname LTb\endcsname
      \put(22,4077){\rotatebox{-270}{\makebox(0,0){\strut{}\footnotesize{imag}}}}%
      \put(1157,2909){\makebox(0,0){\strut{}\footnotesize{real}}}%
      \csname LTb\endcsname
      \put(1823,2677){\makebox(0,0)[r]{\strut{}\footnotesize{$P = 109, N = 100$}}}%
    }%
    \gplgaddtomacro\gplbacktext{%
      \csname LTb\endcsname
      \put(2783,3107){\makebox(0,0)[r]{\strut{}\footnotesize{-1}}}%
      \put(2783,5092){\makebox(0,0)[r]{\strut{}\footnotesize{0}}}%
      \put(2810,2997){\makebox(0,0){\strut{}\footnotesize{0}}}%
      \put(4795,2997){\makebox(0,0){\strut{}\footnotesize{1}}}%
    }%
    \gplgaddtomacro\gplfronttext{%
      \csname LTb\endcsname
      \put(2667,4077){\rotatebox{-270}{\makebox(0,0){\strut{}\footnotesize{imag}}}}%
      \put(3802,2909){\makebox(0,0){\strut{}\footnotesize{real}}}%
      \csname LTb\endcsname
      \put(4468,2677){\makebox(0,0)[r]{\strut{}\footnotesize{$P = 109, N = 600$}}}%
    }%
    \gplgaddtomacro\gplbacktext{%
      \csname LTb\endcsname
      \put(5428,3107){\makebox(0,0)[r]{\strut{}\footnotesize{-1}}}%
      \put(5428,5092){\makebox(0,0)[r]{\strut{}\footnotesize{0}}}%
      \put(5455,2997){\makebox(0,0){\strut{}\footnotesize{0}}}%
      \put(7440,2997){\makebox(0,0){\strut{}\footnotesize{1}}}%
    }%
    \gplgaddtomacro\gplfronttext{%
      \csname LTb\endcsname
      \put(5312,4077){\rotatebox{-270}{\makebox(0,0){\strut{}\footnotesize{imag}}}}%
      \put(6447,2909){\makebox(0,0){\strut{}\footnotesize{real}}}%
      \csname LTb\endcsname
      \put(7113,2677){\makebox(0,0)[r]{\strut{}\footnotesize{$P = 109, N = 1500$}}}%
    }%
    \gplgaddtomacro\gplbacktext{%
      \csname LTb\endcsname
      \put(138,462){\makebox(0,0)[r]{\strut{}\footnotesize{-1}}}%
      \put(138,2447){\makebox(0,0)[r]{\strut{}\footnotesize{0}}}%
      \put(165,352){\makebox(0,0){\strut{}\footnotesize{0}}}%
      \put(2150,352){\makebox(0,0){\strut{}\footnotesize{1}}}%
    }%
    \gplgaddtomacro\gplfronttext{%
      \csname LTb\endcsname
      \put(22,1432){\rotatebox{-270}{\makebox(0,0){\strut{}\footnotesize{imag}}}}%
      \put(1157,264){\makebox(0,0){\strut{}\footnotesize{real}}}%
      \csname LTb\endcsname
      \put(1823,32){\makebox(0,0)[r]{\strut{}\footnotesize{$P = 832, N = 100$}}}%
    }%
    \gplgaddtomacro\gplbacktext{%
      \csname LTb\endcsname
      \put(2783,462){\makebox(0,0)[r]{\strut{}\footnotesize{-1}}}%
      \put(2783,2447){\makebox(0,0)[r]{\strut{}\footnotesize{0}}}%
      \put(2810,352){\makebox(0,0){\strut{}\footnotesize{0}}}%
      \put(4795,352){\makebox(0,0){\strut{}\footnotesize{1}}}%
    }%
    \gplgaddtomacro\gplfronttext{%
      \csname LTb\endcsname
      \put(2667,1432){\rotatebox{-270}{\makebox(0,0){\strut{}\footnotesize{imag}}}}%
      \put(3802,264){\makebox(0,0){\strut{}\footnotesize{real}}}%
      \csname LTb\endcsname
      \put(4468,32){\makebox(0,0)[r]{\strut{}\footnotesize{$P = 832, N = 600$}}}%
    }%
    \gplgaddtomacro\gplbacktext{%
      \csname LTb\endcsname
      \put(5428,462){\makebox(0,0)[r]{\strut{}\footnotesize{-1}}}%
      \put(5428,2447){\makebox(0,0)[r]{\strut{}\footnotesize{0}}}%
      \put(5455,352){\makebox(0,0){\strut{}\footnotesize{0}}}%
      \put(7440,352){\makebox(0,0){\strut{}\footnotesize{1}}}%
    }%
    \gplgaddtomacro\gplfronttext{%
      \csname LTb\endcsname
      \put(5312,1432){\rotatebox{-270}{\makebox(0,0){\strut{}\footnotesize{imag}}}}%
      \put(6447,264){\makebox(0,0){\strut{}\footnotesize{real}}}%
      \csname LTb\endcsname
      \put(7113,32){\makebox(0,0)[r]{\strut{}\footnotesize{$P = 832, N = 1500$}}}%
    }%
    \gplbacktext
    \put(0,0){\includegraphics[width={396.80bp},height={396.80bp}]{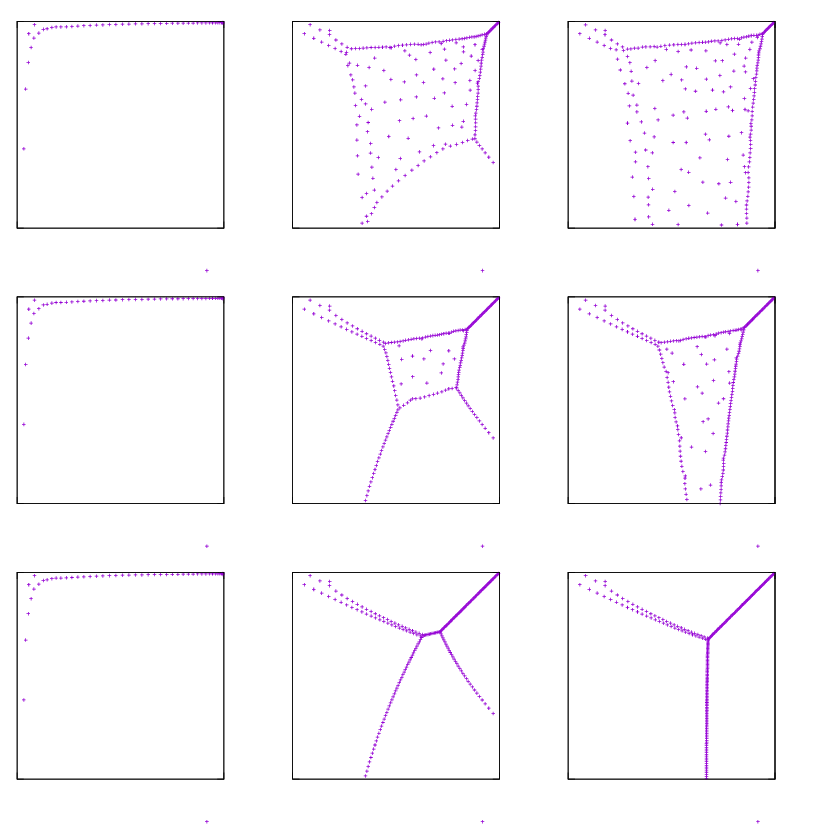}}%
    \gplfronttext
  \end{picture}%
\endgroup

%% file: apcfe-paper-figure-dists_coue_NPmax_Re13000_N_dH.txt
\begingroup
  \makeatletter
  \providecommand\color[2][]{%
    \GenericError{(gnuplot) \space\space\space\@spaces}{%
      Package color not loaded in conjunction with
      terminal option `colourtext'%
    }{See the gnuplot documentation for explanation.%
    }{Either use 'blacktext' in gnuplot or load the package
      color.sty in LaTeX.}%
    \renewcommand\color[2][]{}%
  }%
  \providecommand\includegraphics[2][]{%
    \GenericError{(gnuplot) \space\space\space\@spaces}{%
      Package graphicx or graphics not loaded%
    }{See the gnuplot documentation for explanation.%
    }{The gnuplot epslatex terminal needs graphicx.sty or graphics.sty.}%
    \renewcommand\includegraphics[2][]{}%
  }%
  \providecommand\rotatebox[2]{#2}%
  \@ifundefined{ifGPcolor}{%
    \newif\ifGPcolor
    \GPcolorfalse
  }{}%
  \@ifundefined{ifGPblacktext}{%
    \newif\ifGPblacktext
    \GPblacktexttrue
  }{}%
  \let\gplgaddtomacro\g@addto@macro
  \gdef\gplbacktext{}%
  \gdef\gplfronttext{}%
  \makeatother
  \ifGPblacktext
    \def\colorrgb#1{}%
    \def\colorgray#1{}%
  \else
    \ifGPcolor
      \def\colorrgb#1{\color[rgb]{#1}}%
      \def\colorgray#1{\color[gray]{#1}}%
      \expandafter\def\csname LTw\endcsname{\color{white}}%
      \expandafter\def\csname LTb\endcsname{\color{black}}%
      \expandafter\def\csname LTa\endcsname{\color{black}}%
      \expandafter\def\csname LT0\endcsname{\color[rgb]{1,0,0}}%
      \expandafter\def\csname LT1\endcsname{\color[rgb]{0,1,0}}%
      \expandafter\def\csname LT2\endcsname{\color[rgb]{0,0,1}}%
      \expandafter\def\csname LT3\endcsname{\color[rgb]{1,0,1}}%
      \expandafter\def\csname LT4\endcsname{\color[rgb]{0,1,1}}%
      \expandafter\def\csname LT5\endcsname{\color[rgb]{1,1,0}}%
      \expandafter\def\csname LT6\endcsname{\color[rgb]{0,0,0}}%
      \expandafter\def\csname LT7\endcsname{\color[rgb]{1,0.3,0}}%
      \expandafter\def\csname LT8\endcsname{\color[rgb]{0.5,0.5,0.5}}%
    \else
      \def\colorrgb#1{\color{black}}%
      \def\colorgray#1{\color[gray]{#1}}%
      \expandafter\def\csname LTw\endcsname{\color{white}}%
      \expandafter\def\csname LTb\endcsname{\color{black}}%
      \expandafter\def\csname LTa\endcsname{\color{black}}%
      \expandafter\def\csname LT0\endcsname{\color{black}}%
      \expandafter\def\csname LT1\endcsname{\color{black}}%
      \expandafter\def\csname LT2\endcsname{\color{black}}%
      \expandafter\def\csname LT3\endcsname{\color{black}}%
      \expandafter\def\csname LT4\endcsname{\color{black}}%
      \expandafter\def\csname LT5\endcsname{\color{black}}%
      \expandafter\def\csname LT6\endcsname{\color{black}}%
      \expandafter\def\csname LT7\endcsname{\color{black}}%
      \expandafter\def\csname LT8\endcsname{\color{black}}%
    \fi
  \fi
    \setlength{\unitlength}{0.0500bp}%
    \ifx\gptboxheight\undefined%
      \newlength{\gptboxheight}%
      \newlength{\gptboxwidth}%
      \newsavebox{\gptboxtext}%
    \fi%
    \setlength{\fboxrule}{0.5pt}%
    \setlength{\fboxsep}{1pt}%
    \definecolor{tbcol}{rgb}{1,1,1}%
\begin{picture}(7370.00,3968.00)%
    \gplgaddtomacro\gplbacktext{%
      \csname LTb\endcsname
      \put(990,3701){\makebox(0,0)[r]{\strut{}\footnotesize{$1$}}}%
      \put(990,637){\makebox(0,0)[r]{\strut{}\footnotesize{$10^{-200}$}}}%
      \put(990,1403){\makebox(0,0)[r]{\strut{}\footnotesize{$10^{-150}$}}}%
      \put(990,2169){\makebox(0,0)[r]{\strut{}\footnotesize{$10^{-100}$}}}%
      \put(990,2935){\makebox(0,0)[r]{\strut{}\footnotesize{$10^{-50}$}}}%
      \put(1056,352){\makebox(0,0){\strut{}\footnotesize{$0$}}}%
      \put(1987,352){\makebox(0,0){\strut{}\footnotesize{$100$}}}%
      \put(2917,352){\makebox(0,0){\strut{}\footnotesize{$200$}}}%
      \put(3848,352){\makebox(0,0){\strut{}\footnotesize{$300$}}}%
      \put(4778,352){\makebox(0,0){\strut{}\footnotesize{$400$}}}%
      \put(5709,352){\makebox(0,0){\strut{}\footnotesize{$500$}}}%
      \put(6640,352){\makebox(0,0){\strut{}\footnotesize{$600$}}}%
    }%
    \gplgaddtomacro\gplfronttext{%
      \csname LTb\endcsname
      \put(1428,2016){\makebox(0,0)[l]{\strut{}\footnotesize{$\eps_P \approx 2.220\cdot 10^{-16}$}}}%
      \csname LTb\endcsname
      \put(1428,1862){\makebox(0,0)[l]{\strut{}\footnotesize{$\eps_P \approx 2.242\cdot 10^{-44}$}}}%
      \csname LTb\endcsname
      \put(1428,1708){\makebox(0,0)[l]{\strut{}\footnotesize{$\eps_P \approx 2.264\cdot 10^{-72}$}}}%
      \csname LTb\endcsname
      \put(1428,1554){\makebox(0,0)[l]{\strut{}\footnotesize{$\eps_P \approx 4.572\cdot 10^{-100}$}}}%
      \csname LTb\endcsname
      \put(1428,1400){\makebox(0,0)[l]{\strut{}\footnotesize{$\eps_P \approx 4.616\cdot 10^{-128}$}}}%
      \csname LTb\endcsname
      \put(1428,1246){\makebox(0,0)[l]{\strut{}\footnotesize{$\eps_P \approx 4.661\cdot 10^{-156}$}}}%
      \csname LTb\endcsname
      \put(1428,1092){\makebox(0,0)[l]{\strut{}\footnotesize{$\eps_P \approx 4.707\cdot 10^{-184}$}}}%
      \csname LTb\endcsname
      \put(1428,938){\makebox(0,0)[l]{\strut{}\footnotesize{$\eps_P \approx 9.505\cdot 10^{-212}$}}}%
      \csname LTb\endcsname
      \put(1428,784){\makebox(0,0)[l]{\strut{}\footnotesize{$\eps_P \approx 9.598\cdot 10^{-240}$}}}%
      \csname LTb\endcsname
      \put(1428,630){\makebox(0,0)[l]{\strut{}\footnotesize{$\eps_P = \eps_{\text{min}}$}}}%
      \csname LTb\endcsname
      \put(187,2115){\rotatebox{-270.00}{\makebox(0,0){\strut{}\footnotesize{$d_H(\,\cdot\,,\,\cdot\,)$}}}}%
      \put(4080,154){\makebox(0,0){\strut{}\footnotesize{N}}}%
    }%
    \gplbacktext
    \put(0,0){\includegraphics[width={368.50bp},height={198.40bp}]{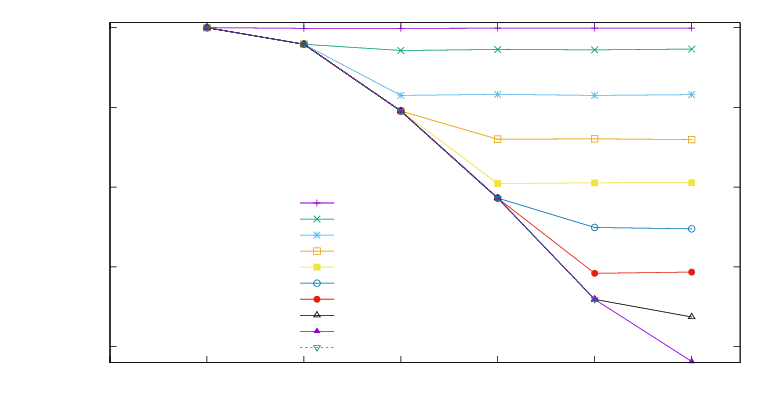}}%
    \gplfronttext
  \end{picture}%
\endgroup

%% file: apcfe-paper-figure-dists_coue_NPmax_Re20000_N_dH.txt
\begingroup
  \makeatletter
  \providecommand\color[2][]{%
    \GenericError{(gnuplot) \space\space\space\@spaces}{%
      Package color not loaded in conjunction with
      terminal option `colourtext'%
    }{See the gnuplot documentation for explanation.%
    }{Either use 'blacktext' in gnuplot or load the package
      color.sty in LaTeX.}%
    \renewcommand\color[2][]{}%
  }%
  \providecommand\includegraphics[2][]{%
    \GenericError{(gnuplot) \space\space\space\@spaces}{%
      Package graphicx or graphics not loaded%
    }{See the gnuplot documentation for explanation.%
    }{The gnuplot epslatex terminal needs graphicx.sty or graphics.sty.}%
    \renewcommand\includegraphics[2][]{}%
  }%
  \providecommand\rotatebox[2]{#2}%
  \@ifundefined{ifGPcolor}{%
    \newif\ifGPcolor
    \GPcolorfalse
  }{}%
  \@ifundefined{ifGPblacktext}{%
    \newif\ifGPblacktext
    \GPblacktexttrue
  }{}%
  \let\gplgaddtomacro\g@addto@macro
  \gdef\gplbacktext{}%
  \gdef\gplfronttext{}%
  \makeatother
  \ifGPblacktext
    \def\colorrgb#1{}%
    \def\colorgray#1{}%
  \else
    \ifGPcolor
      \def\colorrgb#1{\color[rgb]{#1}}%
      \def\colorgray#1{\color[gray]{#1}}%
      \expandafter\def\csname LTw\endcsname{\color{white}}%
      \expandafter\def\csname LTb\endcsname{\color{black}}%
      \expandafter\def\csname LTa\endcsname{\color{black}}%
      \expandafter\def\csname LT0\endcsname{\color[rgb]{1,0,0}}%
      \expandafter\def\csname LT1\endcsname{\color[rgb]{0,1,0}}%
      \expandafter\def\csname LT2\endcsname{\color[rgb]{0,0,1}}%
      \expandafter\def\csname LT3\endcsname{\color[rgb]{1,0,1}}%
      \expandafter\def\csname LT4\endcsname{\color[rgb]{0,1,1}}%
      \expandafter\def\csname LT5\endcsname{\color[rgb]{1,1,0}}%
      \expandafter\def\csname LT6\endcsname{\color[rgb]{0,0,0}}%
      \expandafter\def\csname LT7\endcsname{\color[rgb]{1,0.3,0}}%
      \expandafter\def\csname LT8\endcsname{\color[rgb]{0.5,0.5,0.5}}%
    \else
      \def\colorrgb#1{\color{black}}%
      \def\colorgray#1{\color[gray]{#1}}%
      \expandafter\def\csname LTw\endcsname{\color{white}}%
      \expandafter\def\csname LTb\endcsname{\color{black}}%
      \expandafter\def\csname LTa\endcsname{\color{black}}%
      \expandafter\def\csname LT0\endcsname{\color{black}}%
      \expandafter\def\csname LT1\endcsname{\color{black}}%
      \expandafter\def\csname LT2\endcsname{\color{black}}%
      \expandafter\def\csname LT3\endcsname{\color{black}}%
      \expandafter\def\csname LT4\endcsname{\color{black}}%
      \expandafter\def\csname LT5\endcsname{\color{black}}%
      \expandafter\def\csname LT6\endcsname{\color{black}}%
      \expandafter\def\csname LT7\endcsname{\color{black}}%
      \expandafter\def\csname LT8\endcsname{\color{black}}%
    \fi
  \fi
    \setlength{\unitlength}{0.0500bp}%
    \ifx\gptboxheight\undefined%
      \newlength{\gptboxheight}%
      \newlength{\gptboxwidth}%
      \newsavebox{\gptboxtext}%
    \fi%
    \setlength{\fboxrule}{0.5pt}%
    \setlength{\fboxsep}{1pt}%
    \definecolor{tbcol}{rgb}{1,1,1}%
\begin{picture}(7370.00,3968.00)%
    \gplgaddtomacro\gplbacktext{%
      \csname LTb\endcsname
      \put(990,3711){\makebox(0,0)[r]{\strut{}\footnotesize{$1$}}}%
      \put(990,1338){\makebox(0,0)[r]{\strut{}\footnotesize{$10^{-200}$}}}%
      \put(990,1932){\makebox(0,0)[r]{\strut{}\footnotesize{$10^{-150}$}}}%
      \put(990,2525){\makebox(0,0)[r]{\strut{}\footnotesize{$10^{-100}$}}}%
      \put(990,3118){\makebox(0,0)[r]{\strut{}\footnotesize{$10^{-50}$}}}%
      \put(1056,352){\makebox(0,0){\strut{}\footnotesize{$0$}}}%
      \put(1693,352){\makebox(0,0){\strut{}\footnotesize{$100$}}}%
      \put(2329,352){\makebox(0,0){\strut{}\footnotesize{$200$}}}%
      \put(2966,352){\makebox(0,0){\strut{}\footnotesize{$300$}}}%
      \put(3603,352){\makebox(0,0){\strut{}\footnotesize{$400$}}}%
      \put(4240,352){\makebox(0,0){\strut{}\footnotesize{$500$}}}%
      \put(4876,352){\makebox(0,0){\strut{}\footnotesize{$600$}}}%
    }%
    \gplgaddtomacro\gplfronttext{%
      \csname LTb\endcsname
      \put(1311,2389){\makebox(0,0)[l]{\strut{}\footnotesize{$\eps_P \approx 2.220\cdot 10^{-16}$}}}%
      \csname LTb\endcsname
      \put(1311,2235){\makebox(0,0)[l]{\strut{}\footnotesize{$\eps_P \approx 2.242\cdot 10^{-44}$}}}%
      \csname LTb\endcsname
      \put(1311,2081){\makebox(0,0)[l]{\strut{}\footnotesize{$\eps_P \approx 2.264\cdot 10^{-72}$}}}%
      \csname LTb\endcsname
      \put(1311,1927){\makebox(0,0)[l]{\strut{}\footnotesize{$\eps_P \approx 4.572\cdot 10^{-100}$}}}%
      \csname LTb\endcsname
      \put(1311,1773){\makebox(0,0)[l]{\strut{}\footnotesize{$\eps_P \approx 4.616\cdot 10^{-128}$}}}%
      \csname LTb\endcsname
      \put(1311,1619){\makebox(0,0)[l]{\strut{}\footnotesize{$\eps_P \approx 4.661\cdot 10^{-156}$}}}%
      \csname LTb\endcsname
      \put(1311,1465){\makebox(0,0)[l]{\strut{}\footnotesize{$\eps_P \approx 4.707\cdot 10^{-184}$}}}%
      \csname LTb\endcsname
      \put(1311,1311){\makebox(0,0)[l]{\strut{}\footnotesize{$\eps_P \approx 9.505\cdot 10^{-212}$}}}%
      \csname LTb\endcsname
      \put(1311,1157){\makebox(0,0)[l]{\strut{}\footnotesize{$\eps_P \approx 9.598\cdot 10^{-240}$}}}%
      \csname LTb\endcsname
      \put(1311,1003){\makebox(0,0)[l]{\strut{}\footnotesize{$\eps_P \approx 2.423\cdot 10^{-268}$}}}%
      \csname LTb\endcsname
      \put(1311,849){\makebox(0,0)[l]{\strut{}\footnotesize{$\eps_P \approx 7.645\cdot 10^{-298}$}}}%
      \csname LTb\endcsname
      \put(1311,695){\makebox(0,0)[l]{\strut{}\footnotesize{$\eps_P = \eps_{\text{min}}$}}}%
      \csname LTb\endcsname
      \put(187,2115){\rotatebox{-270.00}{\makebox(0,0){\strut{}\footnotesize{$d_H(\,\cdot\,,\,\cdot\,)$}}}}%
      \put(4080,154){\makebox(0,0){\strut{}\footnotesize{N}}}%
    }%
    \gplbacktext
    \put(0,0){\includegraphics[width={368.50bp},height={198.40bp}]{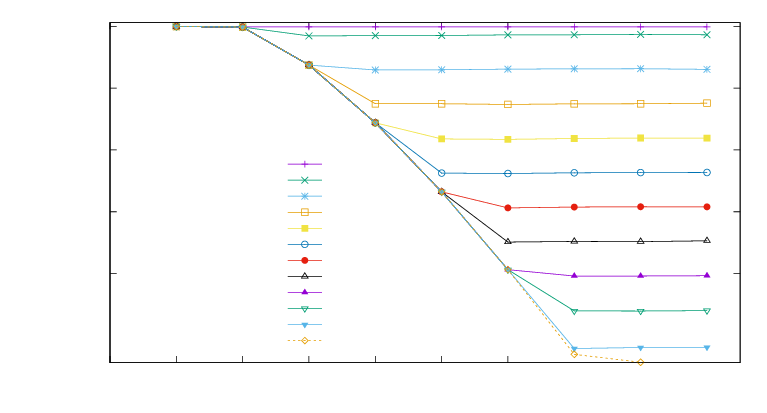}}%
    \gplfronttext
  \end{picture}%
\endgroup

%% file: apcfe-paper-figure-dists_coue_NPmax_Re13000_P_dH.txt
\begingroup
  \makeatletter
  \providecommand\color[2][]{%
    \GenericError{(gnuplot) \space\space\space\@spaces}{%
      Package color not loaded in conjunction with
      terminal option `colourtext'%
    }{See the gnuplot documentation for explanation.%
    }{Either use 'blacktext' in gnuplot or load the package
      color.sty in LaTeX.}%
    \renewcommand\color[2][]{}%
  }%
  \providecommand\includegraphics[2][]{%
    \GenericError{(gnuplot) \space\space\space\@spaces}{%
      Package graphicx or graphics not loaded%
    }{See the gnuplot documentation for explanation.%
    }{The gnuplot epslatex terminal needs graphicx.sty or graphics.sty.}%
    \renewcommand\includegraphics[2][]{}%
  }%
  \providecommand\rotatebox[2]{#2}%
  \@ifundefined{ifGPcolor}{%
    \newif\ifGPcolor
    \GPcolorfalse
  }{}%
  \@ifundefined{ifGPblacktext}{%
    \newif\ifGPblacktext
    \GPblacktexttrue
  }{}%
  \let\gplgaddtomacro\g@addto@macro
  \gdef\gplbacktext{}%
  \gdef\gplfronttext{}%
  \makeatother
  \ifGPblacktext
    \def\colorrgb#1{}%
    \def\colorgray#1{}%
  \else
    \ifGPcolor
      \def\colorrgb#1{\color[rgb]{#1}}%
      \def\colorgray#1{\color[gray]{#1}}%
      \expandafter\def\csname LTw\endcsname{\color{white}}%
      \expandafter\def\csname LTb\endcsname{\color{black}}%
      \expandafter\def\csname LTa\endcsname{\color{black}}%
      \expandafter\def\csname LT0\endcsname{\color[rgb]{1,0,0}}%
      \expandafter\def\csname LT1\endcsname{\color[rgb]{0,1,0}}%
      \expandafter\def\csname LT2\endcsname{\color[rgb]{0,0,1}}%
      \expandafter\def\csname LT3\endcsname{\color[rgb]{1,0,1}}%
      \expandafter\def\csname LT4\endcsname{\color[rgb]{0,1,1}}%
      \expandafter\def\csname LT5\endcsname{\color[rgb]{1,1,0}}%
      \expandafter\def\csname LT6\endcsname{\color[rgb]{0,0,0}}%
      \expandafter\def\csname LT7\endcsname{\color[rgb]{1,0.3,0}}%
      \expandafter\def\csname LT8\endcsname{\color[rgb]{0.5,0.5,0.5}}%
    \else
      \def\colorrgb#1{\color{black}}%
      \def\colorgray#1{\color[gray]{#1}}%
      \expandafter\def\csname LTw\endcsname{\color{white}}%
      \expandafter\def\csname LTb\endcsname{\color{black}}%
      \expandafter\def\csname LTa\endcsname{\color{black}}%
      \expandafter\def\csname LT0\endcsname{\color{black}}%
      \expandafter\def\csname LT1\endcsname{\color{black}}%
      \expandafter\def\csname LT2\endcsname{\color{black}}%
      \expandafter\def\csname LT3\endcsname{\color{black}}%
      \expandafter\def\csname LT4\endcsname{\color{black}}%
      \expandafter\def\csname LT5\endcsname{\color{black}}%
      \expandafter\def\csname LT6\endcsname{\color{black}}%
      \expandafter\def\csname LT7\endcsname{\color{black}}%
      \expandafter\def\csname LT8\endcsname{\color{black}}%
    \fi
  \fi
    \setlength{\unitlength}{0.0500bp}%
    \ifx\gptboxheight\undefined%
      \newlength{\gptboxheight}%
      \newlength{\gptboxwidth}%
      \newsavebox{\gptboxtext}%
    \fi%
    \setlength{\fboxrule}{0.5pt}%
    \setlength{\fboxsep}{1pt}%
    \definecolor{tbcol}{rgb}{1,1,1}%
\begin{picture}(7370.00,3968.00)%
    \gplgaddtomacro\gplbacktext{%
      \csname LTb\endcsname
      \put(990,3708){\makebox(0,0)[r]{\strut{}\footnotesize{$1$}}}%
      \put(990,484){\makebox(0,0)[r]{\strut{}\footnotesize{$10^{-250}$}}}%
      \put(990,1129){\makebox(0,0)[r]{\strut{}\footnotesize{$10^{-200}$}}}%
      \put(990,1774){\makebox(0,0)[r]{\strut{}\footnotesize{$10^{-150}$}}}%
      \put(990,2419){\makebox(0,0)[r]{\strut{}\footnotesize{$10^{-100}$}}}%
      \put(990,3063){\makebox(0,0)[r]{\strut{}\footnotesize{$10^{-50}$}}}%
      \put(7171,264){\makebox(0,0){\strut{}\footnotesize{$1$}}}%
      \put(1355,264){\makebox(0,0){\strut{}\footnotesize{$10^{-250}$}}}%
      \put(2518,264){\makebox(0,0){\strut{}\footnotesize{$10^{-200}$}}}%
      \put(3681,264){\makebox(0,0){\strut{}\footnotesize{$10^{-150}$}}}%
      \put(4844,264){\makebox(0,0){\strut{}\footnotesize{$10^{-100}$}}}%
      \put(6008,264){\makebox(0,0){\strut{}\footnotesize{$10^{-50}$}}}%
    }%
    \gplgaddtomacro\gplfronttext{%
      \csname LTb\endcsname
      \put(6306,690){\makebox(0,0)[l]{\strut{}\footnotesize{$\mathcal O(\eps_P)$}}}%
      \csname LTb\endcsname
      \put(6306,844){\makebox(0,0)[l]{\strut{}\footnotesize{$N = 597$}}}%
      \csname LTb\endcsname
      \put(6306,998){\makebox(0,0)[l]{\strut{}\footnotesize{$N = 497$}}}%
      \csname LTb\endcsname
      \put(6306,1152){\makebox(0,0)[l]{\strut{}\footnotesize{$N = 397$}}}%
      \csname LTb\endcsname
      \put(6306,1306){\makebox(0,0)[l]{\strut{}\footnotesize{$N = 297$}}}%
      \csname LTb\endcsname
      \put(6306,1460){\makebox(0,0)[l]{\strut{}\footnotesize{$N = 197$}}}%
      \csname LTb\endcsname
      \put(187,2115){\rotatebox{-270.00}{\makebox(0,0){\strut{}\footnotesize{$d_H(\,\cdot\,,\,\cdot\,)$}}}}%
      \put(4080,154){\makebox(0,0){\strut{}\footnotesize{$\eps_P$}}}%
    }%
    \gplbacktext
    \put(0,0){\includegraphics[width={368.50bp},height={198.40bp}]{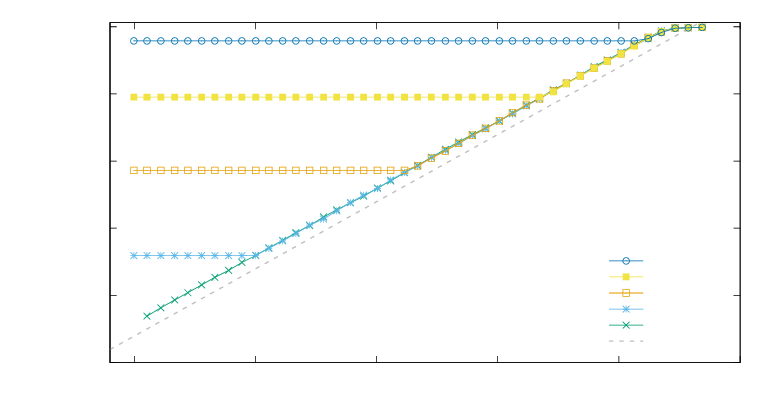}}%
    \gplfronttext
  \end{picture}%
\endgroup

%% file: apcfe-paper-figure-dists_coue_NPmax_Re20000_P_dH.txt
\begingroup
  \makeatletter
  \providecommand\color[2][]{%
    \GenericError{(gnuplot) \space\space\space\@spaces}{%
      Package color not loaded in conjunction with
      terminal option `colourtext'%
    }{See the gnuplot documentation for explanation.%
    }{Either use 'blacktext' in gnuplot or load the package
      color.sty in LaTeX.}%
    \renewcommand\color[2][]{}%
  }%
  \providecommand\includegraphics[2][]{%
    \GenericError{(gnuplot) \space\space\space\@spaces}{%
      Package graphicx or graphics not loaded%
    }{See the gnuplot documentation for explanation.%
    }{The gnuplot epslatex terminal needs graphicx.sty or graphics.sty.}%
    \renewcommand\includegraphics[2][]{}%
  }%
  \providecommand\rotatebox[2]{#2}%
  \@ifundefined{ifGPcolor}{%
    \newif\ifGPcolor
    \GPcolorfalse
  }{}%
  \@ifundefined{ifGPblacktext}{%
    \newif\ifGPblacktext
    \GPblacktexttrue
  }{}%
  \let\gplgaddtomacro\g@addto@macro
  \gdef\gplbacktext{}%
  \gdef\gplfronttext{}%
  \makeatother
  \ifGPblacktext
    \def\colorrgb#1{}%
    \def\colorgray#1{}%
  \else
    \ifGPcolor
      \def\colorrgb#1{\color[rgb]{#1}}%
      \def\colorgray#1{\color[gray]{#1}}%
      \expandafter\def\csname LTw\endcsname{\color{white}}%
      \expandafter\def\csname LTb\endcsname{\color{black}}%
      \expandafter\def\csname LTa\endcsname{\color{black}}%
      \expandafter\def\csname LT0\endcsname{\color[rgb]{1,0,0}}%
      \expandafter\def\csname LT1\endcsname{\color[rgb]{0,1,0}}%
      \expandafter\def\csname LT2\endcsname{\color[rgb]{0,0,1}}%
      \expandafter\def\csname LT3\endcsname{\color[rgb]{1,0,1}}%
      \expandafter\def\csname LT4\endcsname{\color[rgb]{0,1,1}}%
      \expandafter\def\csname LT5\endcsname{\color[rgb]{1,1,0}}%
      \expandafter\def\csname LT6\endcsname{\color[rgb]{0,0,0}}%
      \expandafter\def\csname LT7\endcsname{\color[rgb]{1,0.3,0}}%
      \expandafter\def\csname LT8\endcsname{\color[rgb]{0.5,0.5,0.5}}%
    \else
      \def\colorrgb#1{\color{black}}%
      \def\colorgray#1{\color[gray]{#1}}%
      \expandafter\def\csname LTw\endcsname{\color{white}}%
      \expandafter\def\csname LTb\endcsname{\color{black}}%
      \expandafter\def\csname LTa\endcsname{\color{black}}%
      \expandafter\def\csname LT0\endcsname{\color{black}}%
      \expandafter\def\csname LT1\endcsname{\color{black}}%
      \expandafter\def\csname LT2\endcsname{\color{black}}%
      \expandafter\def\csname LT3\endcsname{\color{black}}%
      \expandafter\def\csname LT4\endcsname{\color{black}}%
      \expandafter\def\csname LT5\endcsname{\color{black}}%
      \expandafter\def\csname LT6\endcsname{\color{black}}%
      \expandafter\def\csname LT7\endcsname{\color{black}}%
      \expandafter\def\csname LT8\endcsname{\color{black}}%
    \fi
  \fi
    \setlength{\unitlength}{0.0500bp}%
    \ifx\gptboxheight\undefined%
      \newlength{\gptboxheight}%
      \newlength{\gptboxwidth}%
      \newsavebox{\gptboxtext}%
    \fi%
    \setlength{\fboxrule}{0.5pt}%
    \setlength{\fboxsep}{1pt}%
    \definecolor{tbcol}{rgb}{1,1,1}%
\begin{picture}(7370.00,3968.00)%
    \gplgaddtomacro\gplbacktext{%
      \csname LTb\endcsname
      \put(990,3716){\makebox(0,0)[r]{\strut{}\footnotesize{$1$}}}%
      \put(990,588){\makebox(0,0)[r]{\strut{}\footnotesize{$10^{-300}$}}}%
      \put(990,1109){\makebox(0,0)[r]{\strut{}\footnotesize{$10^{-250}$}}}%
      \put(990,1631){\makebox(0,0)[r]{\strut{}\footnotesize{$10^{-200}$}}}%
      \put(990,2152){\makebox(0,0)[r]{\strut{}\footnotesize{$10^{-150}$}}}%
      \put(990,2673){\makebox(0,0)[r]{\strut{}\footnotesize{$10^{-100}$}}}%
      \put(990,3194){\makebox(0,0)[r]{\strut{}\footnotesize{$10^{-50}$}}}%
      \put(7171,264){\makebox(0,0){\strut{}\footnotesize{$1$}}}%
      \put(2293,264){\makebox(0,0){\strut{}\footnotesize{$10^{-250}$}}}%
      \put(3268,264){\makebox(0,0){\strut{}\footnotesize{$10^{-200}$}}}%
      \put(4244,264){\makebox(0,0){\strut{}\footnotesize{$10^{-150}$}}}%
      \put(5220,264){\makebox(0,0){\strut{}\footnotesize{$10^{-100}$}}}%
      \put(6195,264){\makebox(0,0){\strut{}\footnotesize{$10^{-50}$}}}%
    }%
    \gplgaddtomacro\gplfronttext{%
      \csname LTb\endcsname
      \put(6261,665){\makebox(0,0)[l]{\strut{}\footnotesize{$\mathcal O(\eps_P)$}}}%
      \csname LTb\endcsname
      \put(6261,819){\makebox(0,0)[l]{\strut{}\footnotesize{$N = 897$}}}%
      \csname LTb\endcsname
      \put(6261,973){\makebox(0,0)[l]{\strut{}\footnotesize{$N = 797$}}}%
      \csname LTb\endcsname
      \put(6261,1127){\makebox(0,0)[l]{\strut{}\footnotesize{$N = 697$}}}%
      \csname LTb\endcsname
      \put(6261,1281){\makebox(0,0)[l]{\strut{}\footnotesize{$N = 597$}}}%
      \csname LTb\endcsname
      \put(6261,1435){\makebox(0,0)[l]{\strut{}\footnotesize{$N = 497$}}}%
      \csname LTb\endcsname
      \put(6261,1589){\makebox(0,0)[l]{\strut{}\footnotesize{$N = 397$}}}%
      \csname LTb\endcsname
      \put(6261,1743){\makebox(0,0)[l]{\strut{}\footnotesize{$N = 297$}}}%
      \csname LTb\endcsname
      \put(6261,1897){\makebox(0,0)[l]{\strut{}\footnotesize{$N = 197$}}}%
      \csname LTb\endcsname
      \put(187,2115){\rotatebox{-270.00}{\makebox(0,0){\strut{}\footnotesize{$d_H(\,\cdot\,,\,\cdot\,)$}}}}%
      \put(4080,154){\makebox(0,0){\strut{}\footnotesize{$\eps_P$}}}%
    }%
    \gplbacktext
    \put(0,0){\includegraphics[width={368.50bp},height={198.40bp}]{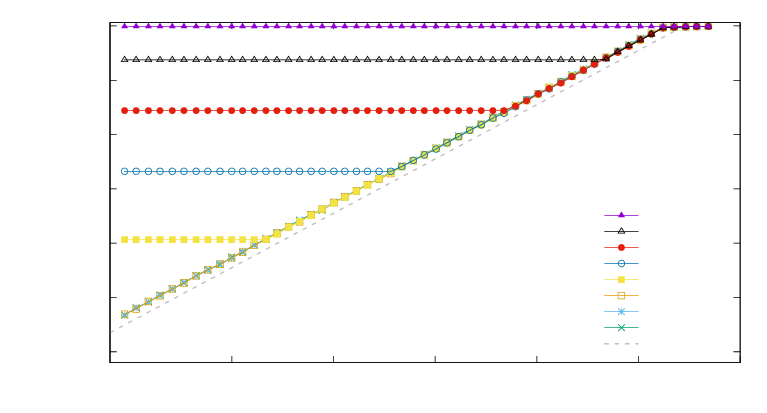}}%
    \gplfronttext
  \end{picture}%
\endgroup

%% file: apcfe-paper-figure-out_coue_20000_spectra.txt
\begingroup
  \makeatletter
  \providecommand\color[2][]{%
    \GenericError{(gnuplot) \space\space\space\@spaces}{%
      Package color not loaded in conjunction with
      terminal option `colourtext'%
    }{See the gnuplot documentation for explanation.%
    }{Either use 'blacktext' in gnuplot or load the package
      color.sty in LaTeX.}%
    \renewcommand\color[2][]{}%
  }%
  \providecommand\includegraphics[2][]{%
    \GenericError{(gnuplot) \space\space\space\@spaces}{%
      Package graphicx or graphics not loaded%
    }{See the gnuplot documentation for explanation.%
    }{The gnuplot epslatex terminal needs graphicx.sty or graphics.sty.}%
    \renewcommand\includegraphics[2][]{}%
  }%
  \providecommand\rotatebox[2]{#2}%
  \@ifundefined{ifGPcolor}{%
    \newif\ifGPcolor
    \GPcolorfalse
  }{}%
  \@ifundefined{ifGPblacktext}{%
    \newif\ifGPblacktext
    \GPblacktexttrue
  }{}%
  \let\gplgaddtomacro\g@addto@macro
  \gdef\gplbacktext{}%
  \gdef\gplfronttext{}%
  \makeatother
  \ifGPblacktext
    \def\colorrgb#1{}%
    \def\colorgray#1{}%
  \else
    \ifGPcolor
      \def\colorrgb#1{\color[rgb]{#1}}%
      \def\colorgray#1{\color[gray]{#1}}%
      \expandafter\def\csname LTw\endcsname{\color{white}}%
      \expandafter\def\csname LTb\endcsname{\color{black}}%
      \expandafter\def\csname LTa\endcsname{\color{black}}%
      \expandafter\def\csname LT0\endcsname{\color[rgb]{1,0,0}}%
      \expandafter\def\csname LT1\endcsname{\color[rgb]{0,1,0}}%
      \expandafter\def\csname LT2\endcsname{\color[rgb]{0,0,1}}%
      \expandafter\def\csname LT3\endcsname{\color[rgb]{1,0,1}}%
      \expandafter\def\csname LT4\endcsname{\color[rgb]{0,1,1}}%
      \expandafter\def\csname LT5\endcsname{\color[rgb]{1,1,0}}%
      \expandafter\def\csname LT6\endcsname{\color[rgb]{0,0,0}}%
      \expandafter\def\csname LT7\endcsname{\color[rgb]{1,0.3,0}}%
      \expandafter\def\csname LT8\endcsname{\color[rgb]{0.5,0.5,0.5}}%
    \else
      \def\colorrgb#1{\color{black}}%
      \def\colorgray#1{\color[gray]{#1}}%
      \expandafter\def\csname LTw\endcsname{\color{white}}%
      \expandafter\def\csname LTb\endcsname{\color{black}}%
      \expandafter\def\csname LTa\endcsname{\color{black}}%
      \expandafter\def\csname LT0\endcsname{\color{black}}%
      \expandafter\def\csname LT1\endcsname{\color{black}}%
      \expandafter\def\csname LT2\endcsname{\color{black}}%
      \expandafter\def\csname LT3\endcsname{\color{black}}%
      \expandafter\def\csname LT4\endcsname{\color{black}}%
      \expandafter\def\csname LT5\endcsname{\color{black}}%
      \expandafter\def\csname LT6\endcsname{\color{black}}%
      \expandafter\def\csname LT7\endcsname{\color{black}}%
      \expandafter\def\csname LT8\endcsname{\color{black}}%
    \fi
  \fi
    \setlength{\unitlength}{0.0500bp}%
    \ifx\gptboxheight\undefined%
      \newlength{\gptboxheight}%
      \newlength{\gptboxwidth}%
      \newsavebox{\gptboxtext}%
    \fi%
    \setlength{\fboxrule}{0.5pt}%
    \setlength{\fboxsep}{1pt}%
    \definecolor{tbcol}{rgb}{1,1,1}%
\begin{picture}(7936.00,5102.00)%
    \gplgaddtomacro\gplbacktext{%
      \csname LTb\endcsname
      \put(185,3013){\makebox(0,0)[r]{\strut{}\footnotesize{-1}}}%
      \put(185,4903){\makebox(0,0)[r]{\strut{}\footnotesize{0}}}%
      \put(1158,2903){\makebox(0,0){\strut{}\footnotesize{0}}}%
      \put(2103,2903){\makebox(0,0){\strut{}\footnotesize{1}}}%
    }%
    \gplgaddtomacro\gplfronttext{%
      \csname LTb\endcsname
      \put(1776,2598){\makebox(0,0)[r]{\strut{}\footnotesize{$P = 53, N = 100$}}}%
      \csname LTb\endcsname
      \put(69,3936){\rotatebox{-270.00}{\makebox(0,0){\strut{}\footnotesize{imag}}}}%
      \put(1157,2815){\makebox(0,0){\strut{}\footnotesize{real}}}%
    }%
    \gplgaddtomacro\gplbacktext{%
      \csname LTb\endcsname
      \put(2830,3013){\makebox(0,0)[r]{\strut{}\footnotesize{-1}}}%
      \put(2830,4903){\makebox(0,0)[r]{\strut{}\footnotesize{0}}}%
      \put(3803,2903){\makebox(0,0){\strut{}\footnotesize{0}}}%
      \put(4748,2903){\makebox(0,0){\strut{}\footnotesize{1}}}%
    }%
    \gplgaddtomacro\gplfronttext{%
      \csname LTb\endcsname
      \put(4421,2598){\makebox(0,0)[r]{\strut{}\footnotesize{$P = 53, N = 200$}}}%
      \csname LTb\endcsname
      \put(2714,3936){\rotatebox{-270.00}{\makebox(0,0){\strut{}\footnotesize{imag}}}}%
      \put(3802,2815){\makebox(0,0){\strut{}\footnotesize{real}}}%
    }%
    \gplgaddtomacro\gplbacktext{%
      \csname LTb\endcsname
      \put(5475,3013){\makebox(0,0)[r]{\strut{}\footnotesize{-1}}}%
      \put(5475,4903){\makebox(0,0)[r]{\strut{}\footnotesize{0}}}%
      \put(6448,2903){\makebox(0,0){\strut{}\footnotesize{0}}}%
      \put(7393,2903){\makebox(0,0){\strut{}\footnotesize{1}}}%
    }%
    \gplgaddtomacro\gplfronttext{%
      \csname LTb\endcsname
      \put(7066,2598){\makebox(0,0)[r]{\strut{}\footnotesize{$P = 53, N = 700$}}}%
      \csname LTb\endcsname
      \put(5359,3936){\rotatebox{-270.00}{\makebox(0,0){\strut{}\footnotesize{imag}}}}%
      \put(6447,2815){\makebox(0,0){\strut{}\footnotesize{real}}}%
    }%
    \gplgaddtomacro\gplbacktext{%
      \csname LTb\endcsname
      \put(185,462){\makebox(0,0)[r]{\strut{}\footnotesize{-1}}}%
      \put(185,2353){\makebox(0,0)[r]{\strut{}\footnotesize{0}}}%
      \put(1158,352){\makebox(0,0){\strut{}\footnotesize{0}}}%
      \put(2103,352){\makebox(0,0){\strut{}\footnotesize{1}}}%
    }%
    \gplgaddtomacro\gplfronttext{%
      \csname LTb\endcsname
      \put(1776,46){\makebox(0,0)[r]{\strut{}\footnotesize{$P = 424, N = 100$}}}%
      \csname LTb\endcsname
      \put(69,1385){\rotatebox{-270.00}{\makebox(0,0){\strut{}\footnotesize{imag}}}}%
      \put(1157,264){\makebox(0,0){\strut{}\footnotesize{real}}}%
    }%
    \gplgaddtomacro\gplbacktext{%
      \csname LTb\endcsname
      \put(2830,462){\makebox(0,0)[r]{\strut{}\footnotesize{-1}}}%
      \put(2830,2353){\makebox(0,0)[r]{\strut{}\footnotesize{0}}}%
      \put(3803,352){\makebox(0,0){\strut{}\footnotesize{0}}}%
      \put(4748,352){\makebox(0,0){\strut{}\footnotesize{1}}}%
    }%
    \gplgaddtomacro\gplfronttext{%
      \csname LTb\endcsname
      \put(4421,46){\makebox(0,0)[r]{\strut{}\footnotesize{$P = 424, N = 200$}}}%
      \csname LTb\endcsname
      \put(2714,1385){\rotatebox{-270.00}{\makebox(0,0){\strut{}\footnotesize{imag}}}}%
      \put(3802,264){\makebox(0,0){\strut{}\footnotesize{real}}}%
    }%
    \gplgaddtomacro\gplbacktext{%
      \csname LTb\endcsname
      \put(5475,462){\makebox(0,0)[r]{\strut{}\footnotesize{-1}}}%
      \put(5475,2353){\makebox(0,0)[r]{\strut{}\footnotesize{0}}}%
      \put(6448,352){\makebox(0,0){\strut{}\footnotesize{0}}}%
      \put(7393,352){\makebox(0,0){\strut{}\footnotesize{1}}}%
    }%
    \gplgaddtomacro\gplfronttext{%
      \csname LTb\endcsname
      \put(7066,46){\makebox(0,0)[r]{\strut{}\footnotesize{$P = 424, N = 700$}}}%
      \csname LTb\endcsname
      \put(5359,1385){\rotatebox{-270.00}{\makebox(0,0){\strut{}\footnotesize{imag}}}}%
      \put(6447,264){\makebox(0,0){\strut{}\footnotesize{real}}}%
    }%
    \gplbacktext
    \put(0,0){\includegraphics[width={396.80bp},height={255.10bp}]{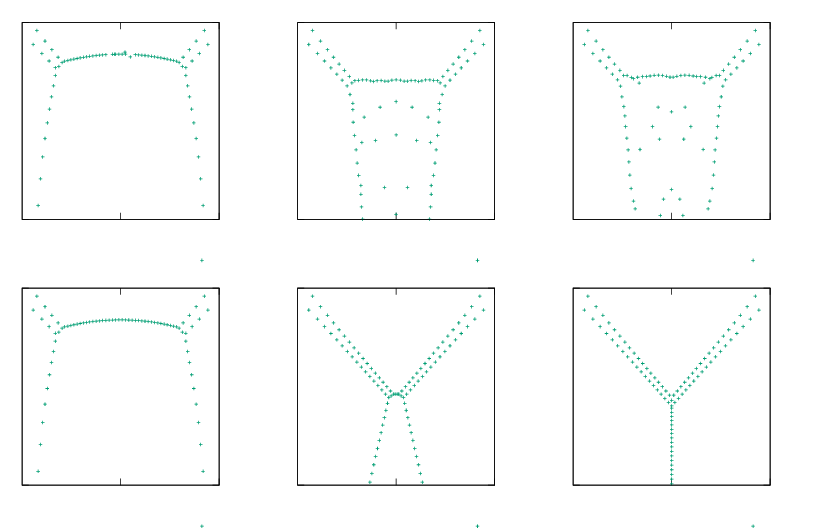}}%
    \gplfronttext
  \end{picture}%
\endgroup

%% file: apcfe-paper-figure-dists_comp_D2D4_P_53.txt
\begingroup
  \makeatletter
  \providecommand\color[2][]{%
    \GenericError{(gnuplot) \space\space\space\@spaces}{%
      Package color not loaded in conjunction with
      terminal option `colourtext'%
    }{See the gnuplot documentation for explanation.%
    }{Either use 'blacktext' in gnuplot or load the package
      color.sty in LaTeX.}%
    \renewcommand\color[2][]{}%
  }%
  \providecommand\includegraphics[2][]{%
    \GenericError{(gnuplot) \space\space\space\@spaces}{%
      Package graphicx or graphics not loaded%
    }{See the gnuplot documentation for explanation.%
    }{The gnuplot epslatex terminal needs graphicx.sty or graphics.sty.}%
    \renewcommand\includegraphics[2][]{}%
  }%
  \providecommand\rotatebox[2]{#2}%
  \@ifundefined{ifGPcolor}{%
    \newif\ifGPcolor
    \GPcolorfalse
  }{}%
  \@ifundefined{ifGPblacktext}{%
    \newif\ifGPblacktext
    \GPblacktexttrue
  }{}%
  \let\gplgaddtomacro\g@addto@macro
  \gdef\gplbacktext{}%
  \gdef\gplfronttext{}%
  \makeatother
  \ifGPblacktext
    \def\colorrgb#1{}%
    \def\colorgray#1{}%
  \else
    \ifGPcolor
      \def\colorrgb#1{\color[rgb]{#1}}%
      \def\colorgray#1{\color[gray]{#1}}%
      \expandafter\def\csname LTw\endcsname{\color{white}}%
      \expandafter\def\csname LTb\endcsname{\color{black}}%
      \expandafter\def\csname LTa\endcsname{\color{black}}%
      \expandafter\def\csname LT0\endcsname{\color[rgb]{1,0,0}}%
      \expandafter\def\csname LT1\endcsname{\color[rgb]{0,1,0}}%
      \expandafter\def\csname LT2\endcsname{\color[rgb]{0,0,1}}%
      \expandafter\def\csname LT3\endcsname{\color[rgb]{1,0,1}}%
      \expandafter\def\csname LT4\endcsname{\color[rgb]{0,1,1}}%
      \expandafter\def\csname LT5\endcsname{\color[rgb]{1,1,0}}%
      \expandafter\def\csname LT6\endcsname{\color[rgb]{0,0,0}}%
      \expandafter\def\csname LT7\endcsname{\color[rgb]{1,0.3,0}}%
      \expandafter\def\csname LT8\endcsname{\color[rgb]{0.5,0.5,0.5}}%
    \else
      \def\colorrgb#1{\color{black}}%
      \def\colorgray#1{\color[gray]{#1}}%
      \expandafter\def\csname LTw\endcsname{\color{white}}%
      \expandafter\def\csname LTb\endcsname{\color{black}}%
      \expandafter\def\csname LTa\endcsname{\color{black}}%
      \expandafter\def\csname LT0\endcsname{\color{black}}%
      \expandafter\def\csname LT1\endcsname{\color{black}}%
      \expandafter\def\csname LT2\endcsname{\color{black}}%
      \expandafter\def\csname LT3\endcsname{\color{black}}%
      \expandafter\def\csname LT4\endcsname{\color{black}}%
      \expandafter\def\csname LT5\endcsname{\color{black}}%
      \expandafter\def\csname LT6\endcsname{\color{black}}%
      \expandafter\def\csname LT7\endcsname{\color{black}}%
      \expandafter\def\csname LT8\endcsname{\color{black}}%
    \fi
  \fi
    \setlength{\unitlength}{0.0500bp}%
    \ifx\gptboxheight\undefined%
      \newlength{\gptboxheight}%
      \newlength{\gptboxwidth}%
      \newsavebox{\gptboxtext}%
    \fi%
    \setlength{\fboxrule}{0.5pt}%
    \setlength{\fboxsep}{1pt}%
    \definecolor{tbcol}{rgb}{1,1,1}%
\begin{picture}(7370.00,2266.00)%
    \gplgaddtomacro\gplbacktext{%
      \csname LTb\endcsname
      \put(678,453){\makebox(0,0)[r]{\strut{}\footnotesize{$0.1$}}}%
      \put(678,985){\makebox(0,0)[r]{\strut{}\footnotesize{$0.2$}}}%
      \put(678,1516){\makebox(0,0)[r]{\strut{}\footnotesize{$0.4$}}}%
      \put(678,2048){\makebox(0,0)[r]{\strut{}\footnotesize{$0.8$}}}%
      \put(906,343){\makebox(0,0){\strut{}\footnotesize{$100$}}}%
      \put(1536,343){\makebox(0,0){\strut{}\footnotesize{$200$}}}%
      \put(2166,343){\makebox(0,0){\strut{}\footnotesize{$400$}}}%
      \put(2796,343){\makebox(0,0){\strut{}\footnotesize{$800$}}}%
      \put(3426,343){\makebox(0,0){\strut{}\footnotesize{$1600$}}}%
    }%
    \gplgaddtomacro\gplfronttext{%
      \csname LTb\endcsname
      \put(1165,856){\makebox(0,0)[r]{\strut{}\footnotesize{$D^2$}}}%
      \csname LTb\endcsname
      \put(1165,636){\makebox(0,0)[r]{\strut{}\footnotesize{$D^4$}}}%
      \csname LTb\endcsname
      \put(271,1336){\rotatebox{-270.00}{\makebox(0,0){\strut{}\footnotesize{$d_H(\,\cdot\,,\,\cdot\,)$}}}}%
      \put(2118,123){\makebox(0,0){\strut{}\footnotesize{N}}}%
    }%
    \gplgaddtomacro\gplbacktext{%
      \csname LTb\endcsname
      \put(4540,806){\makebox(0,0)[r]{\strut{}\footnotesize{$10^{1}$}}}%
      \put(4540,1159){\makebox(0,0)[r]{\strut{}\footnotesize{$10^{2}$}}}%
      \put(4540,1513){\makebox(0,0)[r]{\strut{}\footnotesize{$10^{3}$}}}%
      \put(4540,1866){\makebox(0,0)[r]{\strut{}\footnotesize{$10^{4}$}}}%
      \put(4540,2219){\makebox(0,0)[r]{\strut{}\footnotesize{$10^{5}$}}}%
      \put(4702,343){\makebox(0,0){\strut{}\footnotesize{$100$}}}%
      \put(5332,343){\makebox(0,0){\strut{}\footnotesize{$200$}}}%
      \put(5961,343){\makebox(0,0){\strut{}\footnotesize{$400$}}}%
      \put(6591,343){\makebox(0,0){\strut{}\footnotesize{$800$}}}%
      \put(7221,343){\makebox(0,0){\strut{}\footnotesize{$1600$}}}%
    }%
    \gplgaddtomacro\gplfronttext{%
      \csname LTb\endcsname
      \put(5200,2046){\makebox(0,0)[r]{\strut{}\footnotesize{$D^2$}}}%
      \csname LTb\endcsname
      \put(5200,1826){\makebox(0,0)[r]{\strut{}\footnotesize{$D^4$}}}%
      \csname LTb\endcsname
      \put(4067,1336){\rotatebox{-270.00}{\makebox(0,0){\strut{}\footnotesize{time} [$\mathrm s$]}}}%
      \put(5913,123){\makebox(0,0){\strut{}\footnotesize{N}}}%
    }%
    \gplbacktext
    \put(0,0){\includegraphics[width={368.50bp},height={113.30bp}]{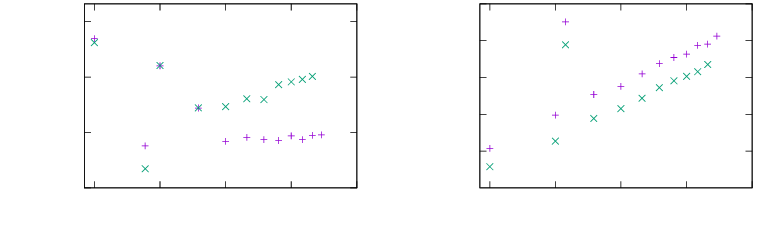}}%
    \gplfronttext
  \end{picture}%
\endgroup

%% file: apcfe-paper-figure-dists_comp_D2D4_P_832.txt
\begingroup
  \makeatletter
  \providecommand\color[2][]{%
    \GenericError{(gnuplot) \space\space\space\@spaces}{%
      Package color not loaded in conjunction with
      terminal option `colourtext'%
    }{See the gnuplot documentation for explanation.%
    }{Either use 'blacktext' in gnuplot or load the package
      color.sty in LaTeX.}%
    \renewcommand\color[2][]{}%
  }%
  \providecommand\includegraphics[2][]{%
    \GenericError{(gnuplot) \space\space\space\@spaces}{%
      Package graphicx or graphics not loaded%
    }{See the gnuplot documentation for explanation.%
    }{The gnuplot epslatex terminal needs graphicx.sty or graphics.sty.}%
    \renewcommand\includegraphics[2][]{}%
  }%
  \providecommand\rotatebox[2]{#2}%
  \@ifundefined{ifGPcolor}{%
    \newif\ifGPcolor
    \GPcolorfalse
  }{}%
  \@ifundefined{ifGPblacktext}{%
    \newif\ifGPblacktext
    \GPblacktexttrue
  }{}%
  \let\gplgaddtomacro\g@addto@macro
  \gdef\gplbacktext{}%
  \gdef\gplfronttext{}%
  \makeatother
  \ifGPblacktext
    \def\colorrgb#1{}%
    \def\colorgray#1{}%
  \else
    \ifGPcolor
      \def\colorrgb#1{\color[rgb]{#1}}%
      \def\colorgray#1{\color[gray]{#1}}%
      \expandafter\def\csname LTw\endcsname{\color{white}}%
      \expandafter\def\csname LTb\endcsname{\color{black}}%
      \expandafter\def\csname LTa\endcsname{\color{black}}%
      \expandafter\def\csname LT0\endcsname{\color[rgb]{1,0,0}}%
      \expandafter\def\csname LT1\endcsname{\color[rgb]{0,1,0}}%
      \expandafter\def\csname LT2\endcsname{\color[rgb]{0,0,1}}%
      \expandafter\def\csname LT3\endcsname{\color[rgb]{1,0,1}}%
      \expandafter\def\csname LT4\endcsname{\color[rgb]{0,1,1}}%
      \expandafter\def\csname LT5\endcsname{\color[rgb]{1,1,0}}%
      \expandafter\def\csname LT6\endcsname{\color[rgb]{0,0,0}}%
      \expandafter\def\csname LT7\endcsname{\color[rgb]{1,0.3,0}}%
      \expandafter\def\csname LT8\endcsname{\color[rgb]{0.5,0.5,0.5}}%
    \else
      \def\colorrgb#1{\color{black}}%
      \def\colorgray#1{\color[gray]{#1}}%
      \expandafter\def\csname LTw\endcsname{\color{white}}%
      \expandafter\def\csname LTb\endcsname{\color{black}}%
      \expandafter\def\csname LTa\endcsname{\color{black}}%
      \expandafter\def\csname LT0\endcsname{\color{black}}%
      \expandafter\def\csname LT1\endcsname{\color{black}}%
      \expandafter\def\csname LT2\endcsname{\color{black}}%
      \expandafter\def\csname LT3\endcsname{\color{black}}%
      \expandafter\def\csname LT4\endcsname{\color{black}}%
      \expandafter\def\csname LT5\endcsname{\color{black}}%
      \expandafter\def\csname LT6\endcsname{\color{black}}%
      \expandafter\def\csname LT7\endcsname{\color{black}}%
      \expandafter\def\csname LT8\endcsname{\color{black}}%
    \fi
  \fi
    \setlength{\unitlength}{0.0500bp}%
    \ifx\gptboxheight\undefined%
      \newlength{\gptboxheight}%
      \newlength{\gptboxwidth}%
      \newsavebox{\gptboxtext}%
    \fi%
    \setlength{\fboxrule}{0.5pt}%
    \setlength{\fboxsep}{1pt}%
    \definecolor{tbcol}{rgb}{1,1,1}%
\begin{picture}(7370.00,2266.00)%
    \gplgaddtomacro\gplbacktext{%
      \csname LTb\endcsname
      \put(678,2195){\makebox(0,0)[r]{\strut{}\footnotesize{$1$}}}%
      \put(678,626){\makebox(0,0)[r]{\strut{}\footnotesize{$10^{-200}$}}}%
      \put(678,1018){\makebox(0,0)[r]{\strut{}\footnotesize{$10^{-150}$}}}%
      \put(678,1411){\makebox(0,0)[r]{\strut{}\footnotesize{$10^{-100}$}}}%
      \put(678,1803){\makebox(0,0)[r]{\strut{}\footnotesize{$10^{-50}$}}}%
      \put(906,343){\makebox(0,0){\strut{}\footnotesize{$100$}}}%
      \put(1536,343){\makebox(0,0){\strut{}\footnotesize{$200$}}}%
      \put(2166,343){\makebox(0,0){\strut{}\footnotesize{$400$}}}%
      \put(2796,343){\makebox(0,0){\strut{}\footnotesize{$800$}}}%
      \put(3426,343){\makebox(0,0){\strut{}\footnotesize{$1600$}}}%
    }%
    \gplgaddtomacro\gplfronttext{%
      \csname LTb\endcsname
      \put(1404,846){\makebox(0,0)[r]{\strut{}\footnotesize{$D^2$}}}%
      \csname LTb\endcsname
      \put(1404,626){\makebox(0,0)[r]{\strut{}\footnotesize{$D^4$}}}%
      \csname LTb\endcsname
      \put(-125,1336){\rotatebox{-270.00}{\makebox(0,0){\strut{}\footnotesize{$d_H(\,\cdot\,,\,\cdot\,)$}}}}%
      \put(2118,123){\makebox(0,0){\strut{}\footnotesize{N}}}%
    }%
    \gplgaddtomacro\gplbacktext{%
      \csname LTb\endcsname
      \put(4540,806){\makebox(0,0)[r]{\strut{}\footnotesize{$10^{1}$}}}%
      \put(4540,1159){\makebox(0,0)[r]{\strut{}\footnotesize{$10^{2}$}}}%
      \put(4540,1513){\makebox(0,0)[r]{\strut{}\footnotesize{$10^{3}$}}}%
      \put(4540,1866){\makebox(0,0)[r]{\strut{}\footnotesize{$10^{4}$}}}%
      \put(4540,2219){\makebox(0,0)[r]{\strut{}\footnotesize{$10^{5}$}}}%
      \put(4702,343){\makebox(0,0){\strut{}\footnotesize{$100$}}}%
      \put(5332,343){\makebox(0,0){\strut{}\footnotesize{$200$}}}%
      \put(5961,343){\makebox(0,0){\strut{}\footnotesize{$400$}}}%
      \put(6591,343){\makebox(0,0){\strut{}\footnotesize{$800$}}}%
      \put(7221,343){\makebox(0,0){\strut{}\footnotesize{$1600$}}}%
    }%
    \gplgaddtomacro\gplfronttext{%
      \csname LTb\endcsname
      \put(5200,2046){\makebox(0,0)[r]{\strut{}\footnotesize{$D^2$}}}%
      \csname LTb\endcsname
      \put(5200,1826){\makebox(0,0)[r]{\strut{}\footnotesize{$D^4$}}}%
      \csname LTb\endcsname
      \put(4067,1336){\rotatebox{-270.00}{\makebox(0,0){\strut{}\footnotesize{time} [$\mathrm s$]}}}%
      \put(5913,123){\makebox(0,0){\strut{}\footnotesize{N}}}%
    }%
    \gplbacktext
    \put(0,0){\includegraphics[width={368.50bp},height={113.30bp}]{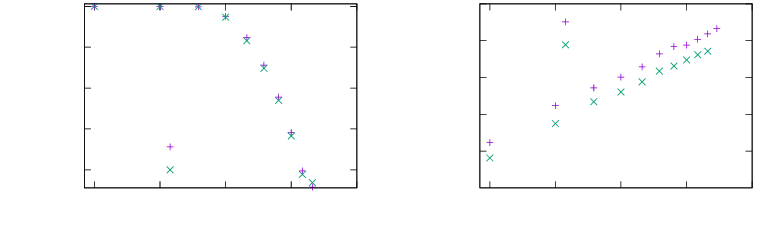}}%
    \gplfronttext
  \end{picture}%
\endgroup